\documentclass{amsart}

\usepackage{amsmath, amsthm}
\usepackage{amssymb, latexsym}
\usepackage{amsfonts}

\newtheorem{Theorem}[equation]{Theorem}
\newtheorem{theorem}[equation]{Theorem}

\newtheorem{proposition}[equation]{Proposition}
\newtheorem{lemma}[equation]{Lemma}
\newtheorem{Lemma}[equation]{Lemma}

\newtheorem{corollary}[equation]{Corollary}

\newtheorem{Remark}[equation]{Remark}
\newtheorem{Remarks}[equation]{Remarks}

\newtheorem{prop}[equation]{Proposition}

\newtheorem{Examples}[equation]{Examples}

\newcounter{com}

\newtheorem{comments}{}

\newtheorem*{A'}{Corollary A$'$}

\newtheorem*{C}{Theorem C} 

\newtheorem*{C'}{Theorem C'}
\newtheorem*{C''}{Theorem B$''$ (Holt's Conjecture)}

\newtheorem*{CC}{Conjecture C}

\newcommand{\eqn}[1]{(\ref{#1})}
\newcommand{\beql}[1]{\begin{equation}\label{#1}}
\newcommand{\eeq} {\end{equation}}

\font\Aaa=msam10

\def\noreq{\hbox{ {\Aaa  E} }}

\def\qed{\hbox{~~\Aaa\char'003}}

\font\Bbb=msbm10
\def\semi{\hbox{\Bbb o}}
\def\R{\hbox{\Bbb R}}

\def\Z{\hbox{\Bbb Z}}

\def\F{\hbox{\Bbb F}}

\numberwithin{equation}{section}

\def\Proof {\noindent{\em Proof.} }

\let\define=\def
\let\redefine=\def

\def\Sym{{\rm Sym}}
\def\Alt{{\rm Alt}}

\define\a{{ \alpha }}
\redefine\b{{ \beta }}
\redefine\c{{ \gamma }}

\redefine\g{\gamma}

\define\z{{ \zeta }}

\define\({ \left( }
\define\){ \right) }
\define\[{ \left[ }
\define\]{ \right] }
\define\<{ \langle }
\define\>{ \rangle }
\let\ljunk=\{
\let\rjunk=\}
\redefine\{{\left\ljunk}
\redefine\}{\right\rjunk}

\redefine\.{ \diamomd }
\redefine\-{ \ominus }
\redefine\+{\oplus}
\redefine\.{ \cdot}
\redefine\#{\sharp}

\redefine\.{ \cdot }

        \def\PSL{{\rm PSL}}
        \def\SL{{\rm SL}}

        \def\SU{{\rm SU}}
        
        \def\PSU{{\rm PSU}}
        
        \def\PGL{{\rm PGL}}
        
        \def\PSp{{\rm PSp}}
        \def\PO{\mbox{{\rm P}$\Omega$}}
        \def\Sp{{\rm Sp}}
        \def\GF{{\rm GF}}

        \def\AGL{{\rm AGL}}

     \def\diag{{\rm diag}}

        \def\Sz{{\rm Sz}}

        \font\Aaa=msam10

        \def\noreq{\mathop{\hbox{\Aaa E}}}

\DeclareRobustCommand{\SkipTocEntry}[4]{}
\begin{document}

\title[Presentations of finite simple groups]
{Presentations of finite simple groups:  \\a computational  approach}

\thanks{The authors were partially supported by
        NSF grants DMS 0140578,
          DMS 0242983, DMS 0600244  and   DMS~0354731. The research by the third author also was 
         supported by Centennial Fellowship from the American Mathematics Society.}

       \author{R. M. Guralnick}
       \address{Department of Mathematics, University of Southern California,
       Los Angeles, CA 90089-2532 USA}
       \email{guralnic@usc.edu}

       \author{W. M. Kantor}
       \address{Department of Mathematics, University of Oregon,
       Eugene, OR 97403 USA}
       \email{kantor@math.uoregon.edu}

    \author{M. Kassabov}
       \address{Department of Mathematics, Cornell University,
Ithaca, NY 14853-4201  USA}
       \email{kassabov@math.cornell.edu}

    \author{A. Lubotzky}
       \address{Department of Mathematics, Hebrew University, Givat Ram,
Jerusalem 91904 Israel}
       \email{alexlub@math.huji.ac.il}

\subjclass[2000]{Primary 20D06, 20F05 Secondary 20J06}

{\abstract 

All nonabelian finite
simple groups of rank $n$ over a field of size $q$,
with the possible exception of the Ree groups $^2G_2(3^{2e+1})$,  have
presentations with  at most  $80 $  relations
and bit-length $O(\log n +\log q)$. 
Moreover,  $A_n$ and $S_n$  have  
presentations with
$3$ generators$,$  $7$ relations and bit-length $O(\log n)$,
while
$\SL(n,q)$   has a  presentation with   $7$ generators, $2 5$ relations
and bit-length $O(\log n +\log q)$. 
 
}

\maketitle


\tableofcontents


\section {Introduction}
\label {Introduction}
  
In \cite{GKKL} we provided short presentations for all alternating
groups, and all finite simple groups of  Lie type  other than  the
Ree groups $^2G_2(q),$ using at most 1000 generators and relations.  
In \cite{GKKL2} we proved the existence of  profinite presentations for the
same groups using fewer than 20 relations.  The goal of the present paper is
similar:  we will provide  presentations for the  same simple groups using 2
generators and at most 80 relations.  
%
%
These and other  new  presentations   have the
potential advantage that they are simpler than those in \cite{GKKL},  at
least in the sense of requiring fewer relations; we hope that both types
of presentations will turn out  to be useful in Computational Group
Theory.

The fundamental difference between this paper and \cite{GKKL} 
is that here we achieve a smaller number of  relations 
at the cost of relinquishing some   control over the length of the
presentations.
Our first result does not deal with lengths at all:

{
\renewcommand{\theequation}{A}
\begin{theorem}
\label{A}
All  nonabelian finite simple 
groups of Lie type$,$  with the possible exception of the
Ree groups $^2G_2(q),$  have  presentations with $2$ generators 
and  at most
$80 $ relations$.$ 

All symmetric and alternating groups have presentations 
with $2$ generators 
and   $8$ relations$.$  
\end{theorem}
\addtocounter{equation}{-1}
}

In fact, a similar result holds for   \emph{all} finite simple groups,
except perhaps $^2G_2(q) $
(the sporadic groups are surveyed in \cite{Soi}).
Both the bounds of $20$ relations in
\cite{GKKL2} and   
$80$ here are  not optimal -- in all cases we will provide much
better bounds, though usually with more generators.  Possibly  $4$ is the
correct upper bound for both standard and profinite  presentations.
Wilson
\cite{Wil} has even conjectured that 2 relations suffice for the
universal covers of all finite  simple groups.

Although we are giving up the  requirement of length  used  in \cite{GKKL}, we
can still retain some control over a weaker notion of length  
used in
\cite{BGKLP,Con} and    especially suited for Computer Science complexity
considerations:  \emph{bit-length}.  This is the total number of  bits
required  to write the presentation, so that  all exponents are encoded as
binary strings,  the sum of whose lengths enters into the bit-length.
(The presentation length that had to be kept small in \cite{GKKL} 
involves at least   the
much larger sum of the actual exponents; cf. Section~\ref{Preliminaries}.)

{
\renewcommand{\theequation}{B}
\begin{theorem}
\label{B}
All nonabelian finite
simple groups of rank $n$ over a field of size $q,$
with the possible exception of the Ree groups  
$^2G_2(q),$   have presentations   with
at most $1 4$ generators$,$     $7 8$ relations  and   bit-length  $O(\log
n+\log q)
$.%
\footnote{Logarithms will 
 be to the base 2.}%
\end{theorem}
\addtocounter{equation}{-1}
}

Here we view alternating  (and symmetric)  groups as
having rank $n-1$ over  ``the  field of size 1'' \cite{Tits_order_1}.
The   bounds in both of the
preceding theorems  also hold   for many of the almost simple groups.

By \cite[Lemma~2.1]{GKKL}, 
if we have any presentation of a finite simple group  $G$  with 
at most  $7 8$
 relations, we obtain a presentation with $2$ generators and
 at most $80$  relations (cf.~Lemma~\ref{d generators} below). 
Moreover, the proof of that lemma shows that {\em any} pair of generators
of $G$ can be used for such a presentation (cf. Corollary~\ref{any
generators of Sn}).  ``Almost all'' pairs of elements of a   finite 
simple group generate the group
\cite{Di,KLu,LiSh}; some pairs will probably force the length or even the
bit-length to be fairly large.  Indeed,
\cite[Lemma~2.1]{GKKL} is so general that it allows us to
cheat somewhat: 
the resulting presentations are not even slightly explicit,  and we have no
information concerning their bit-lengths.  In particular, we are  unable to
prove  Theorem~\ref{B} using   2 instead of 14 generators. 

In view of \cite[Lemma~2.1]{GKKL}  or  Lemma~\ref{d generators}, our
goal will be to prove Theorem~\ref{B}. Much better bounds are obtained in
various cases.  For example, Theorems~\ref{Alt(n) and Sym(n) 100}  and
\ref{Alt(n) and Sym(n) 100 down to 7}   deal  with the alternating and
symmetric groups, where we go further than any previous types    of
presentations for these groups   in terms of the small number of
relations   used (cf. \cite{GKKL,Con}): 

{
\renewcommand{\theequation}{C}
\begin{theorem}
\label{C}
For   each   $n\ge 5,$ $A_n$ and $S_n$  have 
presentations with
$3$ generators$,$  $7$ relations and  bit-length $O(\log n),$
using   a bounded number of exponents. 
\end{theorem}
\addtocounter{equation}{-1}
}
\noindent 
Moreover, for the preceding groups,  in addition to the second part of
Theorem~\ref{A} we show that,  {\em    if 
$a$ and
$b$ are}  any  {\em   generators of $G= A_n$   or $S_n,$    then 
there is a presentation of
$G$ using $2$ generators that map onto 
$a$ and $b,$  with     $ 9$ relations} (Corollary~\ref
{any generators of Sn}).
There are similar results in  all cases of Theorem~\ref{A}
(cf. Remark~\ref{one less generator} in Section~\ref{Concluding
remarks}). 
However, as already  noted, we  do not know   if it is possible
to choose $a$ and $b$  in order to obtain
a presentation with bit-length $O(\log n)$;
nor   if it is possible
to choose  them    in order to obtain
a bounded number of  exponents for  groups of Lie type over 
arbitrarily large degree extensions of a prime field   (cf.
Remark~\ref{Concluding remarks expo-length} in Section~\ref{Concluding
remarks}). 

In order to obtain all of the   presentations in the preceding theorems, 
although \cite{GKKL} was a starting point we need significantly new
methods for unbounded rank $n$;  these ideas  may prove to be
more practical for actual group computation  than some of those in
\cite{GKKL}.
Moreover, while  a few  of the arguments used here are  
streamlined, often simpler, and
occasionally improved  versions of  ones in \cite{GKKL}, they are still
rather involved. 
As in \cite{GKKL} we do not use the 
 classification of the finite simple groups
in any proofs.

For groups of bounded rank, our presentations can be made to be  short in
the sense used in \cite{GKKL}, at the cost of adding
a small number of additional generators and relations (so that 
\cite[(3.3)~and~(4.16)]{GKKL} 
will apply; cf.    (\ref{Horner})  and 
(\ref{Horner for polynomials}) below).  It is our treatment of  unbounded rank  that contains new
ideas to decrease  the
number of relations in \cite{GKKL} at the expense of the length of the
presentation.
We provide more than one approach for this purpose:  
 some classical groups  are handled in different ways in
Sections~\ref {General case}  and \ref {More presentations of classical
groups}.  The unitary groups are dealt with separately in
Section~\ref{Presentations of unitary groups } by using an idea  of
Phan~\cite{Ph} as improved in \cite{BeS}.%

In Sections~\ref{Alternating groups}--\ref{More presentations of classical
groups}  we will consider
various types of simple groups in order to prove the above theorems,
  providing better bounds for the number  of generators and
relations in various cases.
For many cases  we only  give hints regarding
the final assertion in Theorem~\ref{B}.

One of our original motivations for work on presentations was   the
following 

{
\renewcommand{\theequation}{D}
\begin{corollary}[Holt's Conjecture \cite{Holt} for simple groups]
\label{HOLT}
\label{D}
There is a  constant $C$\break
 such that$,$
for every  
finite simple  group $G$$,$ every prime
$p$ and every  irreducible  $\F_p[G]$-module $M ,$ $\dim
H^2(G,M)\leq  C\dim  M$.
\end{corollary} 
\addtocounter{equation}{-1}
}
This conjecture has already been proven twice, in   
\cite[Theorem~B$'$]{GKKL}  and \cite[Theorem~B]{GKKL2}. 
  As in \cite{GKKL} it is an immediate consequence of
Theorem~\ref{A}, except for the Ree groups (and these also had to be handled 
separately in \cite{GKKL}).  The proof based on the present paper
(using the elementary result  \cite[Lemma~7.1]{GKKL}) is simpler
than the previous proofs, although a smaller, explicit constant $C$ is
given   in
 \cite{GKKL2}.
See  \cite[Theorem~C]{GKKL2} for a generalization of the preceding
result to all finite groups.

Section~\ref{Concluding remarks} contains further remarks concerning these 
results.    For now we note one further unexpected direction:
 
\medskip 
 
 \noindent
 {\bf Efficient presentations. } 
\label{Efficient presentations}
 If  $\<X\mid R\>$ is a presentation of a
finite group 
$G$, then $|R|-|X|$ is at least the smallest number $d(M)$ of generators of
the Schur multiplier $M$ of $G$; and $\<X\mid R\>$ is called an
\emph{efficient} presentation if  $|R|-|X|= d(M)~$ \cite{CRKMW, CHRR,
CHRR2, CHR}.  The only infinite families of nonsolvable groups known to
have efficient presentations appear to be groups having $\PSL(2,p)$ 
as a composition factor  when
$p$ is prime
\cite{Sun,CR3} (cf.~(\ref{efficient SL2}).  Therefore  it may be of some interest   that
Corollaries~\ref{Ap+2 4 relations}(i) and
\ref{Sp+2 4 relations}(ii) contain  examples of families of groups  having   efficient 
presentations with alternating groups as composition factors.   For
example,
\emph{for any prime  $p \equiv 11$} (mod 12), \emph{there is a 
a presentation of  $A_{p+2}\times T$ with $2$ generators and $3$
relations$,$ where $T$ is the subgroup of index $2$ in $\AGL(1,p)$.}  It
seems plausible that this can be used to   obtain   efficient
presentations  of 
$A_{p+2}$   
with 3 generators and 4 relations  for these primes.  

Examples~\ref{summary of special n} and \ref{more AGL1 examples}, together
with  Table~\ref{Small degrees}
and Remark~\ref{degree q+2}, deal with presentations for various groups  
$A_n$ and
$S_n$ when
$n$ has a special form.   Section~\ref {An explicit presentation for
$S_n$}  contains explicit presentations  of $S_n$   for all   $n\ge
50$. For general $n$ it would be desirable to have even fewer relations
than in Theorem~\ref{C}, with the goal of approaching efficiency for
alternating groups.

\section{Preliminaries}
\label{Preliminaries}
\noindent
 {\bf Presentation  lengths.}
In \cite[Section~1.2]{GKKL} there is a long discussion of
various notions of ``lengths''  of a presentation
$\langle  X\mid R\rangle$ and some of the relationships among them.  Here
we only summarize what is needed for the present purposes.
 
\begin {itemize}
\item[]
\begin {itemize}
\item [] \hspace{-30pt} \emph{length $=$ word length}: $  |X| \, + $ sum of
the lengths of the words in $R$  within the free group on $X$.
Thus, length refers to strings in the alphabet $X\cup X^{-1}$.  This is
the notion of length used in \cite{GKKL}, and seems the most natural
notion from a purely mathematical point of view.  We reserve the term \emph{short
presentation} for one having small length.  Achieving this was one of the goals 
in
\cite{GKKL}, though not of the present paper. 
\item [] \hspace{-30pt} \emph{bit-length}: the total
number of  bits required  to write the presentation,  used in  \cite{BGKLP} 
and    \cite{Con}. 
All exponents are encoded as binary strings, the sum of whose
lengths enters into the bit-length, as does the space required  to
enumerate the list of  generators and relations. 

\item [] \hspace{-30pt} \emph{expo-length}:  the
total number of 
exponents used in   the presentation. 
\end {itemize}
\end {itemize}
 
\smallskip
\noindent
  By comparing the present paper with
\cite{GKKL} it becomes clear that small bit-length
is much easier to achieve than small length.  The properties required  of the bit-length
$bl(w)$ of a word $w$ are
$$
\begin{array} {llll}
bl(x)=1, x\in X; \,\,\,\,   bl(w^n)\le bl(w)+\log |n| \, \mbox{ if } \,n\in \Z\backslash\{0,1,-1\};  \, \mbox{ and }   
\\
bl(ww')\le bl(w)+bl(w') 
\end{array}
$$
for any words  $w,w'.$

\smallskip
  
\noindent{\bf Subgroups.} %
We will use  the elementary fact
\cite[Lemma~2.3]{GKKL}    that a group $G$ 
that has   a presentation   based
on  a  known group,  using  presentations of  
subgroups of that group, has the subgroups automatically built into 
$G$:
 
\begin{lemma}
\label{It's a subgroup}
\label{it's a subgroup}
\label{It's a group}
Let $\pi\colon F_{X\cup Y} \to G=\langle  X,Y\mid R,S\rangle$ and
$\lambda\colon F_X\to H=\langle  X\mid R\rangle $ be the natural
surjections$,$ where $H$ is finite.~Assume that $\alpha \colon G\to G_0$ is
a homomorphism such that $ \alpha \langle \pi(X )  \rangle\cong H$.~Then
$\langle \pi(X ) \rangle \cong H$.%
\end{lemma}

  In the present paper we will  use this freely, often
  without    comments.

\smallskip
\smallskip
\noindent{\bf Curtis-Steinberg-Tits presentation.}
This is a standard presentation for groups of Lie type; see
\cite{Cu}, \cite{St1},  \cite[Theorem~13.32]{Tits
book} and
\cite[Theorem~2.9.3]{GLS}. We will generally refer to 
\cite[Sections 5.1 and 5.2]{GKKL} for a discussion of the versions we will
use.
 
\section{Symmetric and alternating groups}
\label{Alternating groups}


We will use  a presentation for alternating groups, due to
Carmichael~\cite[p.~169]{Carm2}, that is  more symmetrical than
a presentation due to Burnside and Miller 
(\cite[p.~464]{Bur}, \cite[p.~366]{Mil}) in 1911
that was 
used in
\cite[(2.6)]{GKKL}.  Moreover, Carmichael's 
  presentation requires less data (i.e., fewer relations):
\begin{equation}
   \label{Carmichael An}
   A_{n+2}=
\langle  x_1, \ldots, x_{n } \mid x_i^3= (x_ix_{j})^2= 1~
    \mbox{~  whenever  $~i\ne j$}\rangle ,
   \end{equation}
based on the 3-cycles $(i,n+1,n+2) ,$ $ 1 \le i \le n$.
We will also use  the
standard  presentations 
\begin{equation}
   \label{A4 A5}
\mbox{$A_4=\langle x, y \mid x^3=y^2=(xy)^3=1 \rangle~ $   and 
$	~A_5=\langle x, y \mid x^5=y^2=(xy)^3=1 \rangle $}
   \end{equation} 
\cite[p.~137]{CoMo}.
 Presentations are known for the universal central extension of 
$A_n$,
$4\le n\le9$,  with 2 generators and 2 relations
\cite{CHRR,CHRR2}; and  for $A_{10}$ using 2 generators and 3 relations
\cite{Hav} (cf. Example~\ref{more AGL1 examples}(\ref{A10})).   These  can be used in some of our 
presentations in Sections~\ref{PSL(n,q)} and \ref{Presentations of unitary
groups }  in order to  save  several relations.

In Section~\ref{Using 2-transitive groups} 
we make crucial use of  the symmetry of    \eqn{Carmichael An}, as follows. 
Let $T=\<X\mid R\>$ be a group
 acting transitively on $\{1,\dots,n\}$,  viewed as acting on
$\{1,\dots,n,n+1,n+2\}$.  Introduce a new generator $z$ corresponding to
the 3-cycle $(1,n+1,n+2)$.  We use additional relations in order  to
guarantee that $|z^T|=n$ in our presented group, as well as    a very small
number of   relations of the form $(zz^t)^2=1$ for suitable $t\in
T$, in order to use
\eqn{Carmichael An}.  

This idea is reworked several times  in
order to handle various special degrees $n$.  We glue two such
presentations  in Section~\ref{All alternating and symmetric groups} in order
to deal with   symmetric and alternating groups of arbitrary degrees.

\subsection{Using 2-transitive groups for special degrees}
\label{Using 2-transitive groups}
The results of this section are summarized in Examples~\ref{summary of
special n}.  
 We begin with 
an integer $n\ge 3$, together with
\begin{itemize} 
\item a  group $T$  acting  transitively (though not necessarily faithfully)
on the unordered pairs of distinct 
points in $\{1,\dots,n\}$ (i.e., $T$ is {\em $2$-homogeneous}),  
\item  
a presentation $\langle X \mid R\rangle$  of
$T$,
\item a subset $X_1$ of $T$ such that $T_1 =\<X_1\>$ is  the
stabilizer of
$1$ (we usually have $X_1\subset X$), and
\item  a word $w $ in 
$ X$ that moves $1$ and lies in $A_n$ (when  $w$ is  viewed inside $T$).
\end{itemize} 
\smallskip
The last requirement implies  that the permutation group $\bar T$ induced by
$T$ may not be inside
$A_n$; when it is inside then the  
following lemma and its proof are somewhat simpler.  Note that, if $T$ is
not 2-transitive, then $\bar T$ has odd order, and hence lies in $A_n$.
(Here the order is odd because an involution  in  $\bar T$   would allow some
\emph{ordered} 2-set to be moved to any  other one.)

\begin{lemma}  
\label{using Carmichael}
If $J = \langle X, z \mid  R, z^3=1, (zz^w)^2=1, 
z^u=z^{{\rm sign}(u) } $ for $u\in X_1
\rangle,$ then $J \cong A_{n+2} \times T$.
\end{lemma}

\begin{proof} View $T$ as a subgroup of $A_{n+2} = \Alt\{1,
\ldots, n, n+1,n+2\}$,  with  each $t\in T$
inducing $( n+1,n+2)^{{\rm sign}(t) } $ on $\{ n+1,n+2 \} $.   Define 
$\varphi\colon X\cup \{z\}
\rightarrow A_{n+2}\times T$ by 
\begin{equation}
\label{def phi}
\begin{array}{lll}
\mbox{$\varphi(x)=(x( n+1,n+2)^{{\rm sign}(x) }  ,x)~$ for  $x\in X$ } 
\vspace{2pt}\\
\mbox{$ \varphi(z)=(z',1)~$
where $~z'=(1 , n+1 , n+2)  $.} 
\end{array}
\end{equation}
  Then the image of $\varphi$ satisfies the
defining relations for $J$,  
and we obtain a\break 
 homomorphism $\varphi\colon J \rightarrow A_{n+2}\times T$.
We claim  that  $\varphi$  is a surjection.
For, since $T$ is 2-homogeneous  on $\{1,\ldots,n\}$,  we have  
$\langle\varphi(z)^ {\langle\varphi(X)\rangle } \rangle =A_{n+2} \times
1$. Here, $A_{n+2}$  contains $\langle  X \rangle =T $,
while  $\langle\varphi(X)\rangle$ projects onto $T$ in the second
component, so that $\varphi(J)$ also contains $1\times T$. Hence, 
$\varphi$ is, indeed, surjective.

Then there is  also   
a surjection $\pi\colon  J \rightarrow A_{n+2}$. 
By Lemma~\ref{It's a subgroup}, $J$ has  a subgroup  we identify
with 
$T=\<X\>  $.  We also view $z$ as contained in $J$.

 Since $\<z\> ^{T_1}=\<z\>$ by our relations, we have  $| \<z\>^T|\le n$;
but
$|\pi(\<z\>^T)|=n$   and so
$|\<z\>^T|=n$. Similarly, $|z ^{T\cap A_n}|=n$.
Consequently,  $T$ acts   on $\<z\>^T$ 
and  $T\cap A_n$ acts   on $z ^{T\cap A_n}$
as they do  on $\{1,\ldots, n\}$.

Moreover, if $T$ is not inside $A_n$, then  our sign condition in
the presentation implies that $|z^T|=2n$, and
$T\cap A_n$ has 2 orbits on $z^T$, namely,  $z^{T\cap A_n}$ and
$(z^{-1})^{T\cap A_n}$.

By   $2$-homogeneity, any   unordered pair of distinct  
members of
$\<z\>^T$ is $T$-conjugate to  $\{\<z\>,\<z\>^ w\}$.  If $T\cap A_n$ is 
$2$-homogeneous, then  any   unordered pair of distinct  members of $  z^T$ is
$T\cap A_n$--conjugate to  $\{z , z ^ w \}$.  Since the  relation 
$(zz^ w)^2=1$ in the lemma implies that  $(z^ wz)^2=1$,
it follows that  $z^T$ satisfies 
(\ref{Carmichael An}), so that   $N:=\< z ^T \>\cong A_{n+2}$.

 If $T\cap A_n$ is not $2$-homogeneous then, by hypothesis,   $T$ is
2-transitive but 
$T\cap A_n$ is not.  We claim that   we still have  $N:=\< z ^T \>\cong
A_{n+2}$.  For,  any   ordered pair of distinct  members of $\<z\>^T$ is
$T\cap A_n$--conjugate  to  $(\<z\>,\<z\>^ w )$ or to one other pair,
$(\<z\>,\<z\>^y )$, say, where $y\in T\cap A_n$.  \vspace{2pt}
 Some $g\in T\backslash A_n$ satisfies
  $(\<z\>,\<z\>^ w )^g=(\<z\>,\<z\>^y )$.  
Clearly, $z, z^w,  z^y\in z^{T\cap A_n}$. Since $g\notin A_n$, it follows
that both  $z^g$ and $(z^w)^g$ lie  in the other $T\cap A_n$--class 
$(z^{-1})^{T\cap A_n}$ of $z$. Thus,  $z^g=z^ {-1}$ and
$(z^w)^g=(z^y)^{-1}$, so that $(z   z ^y)^2=([z ^{-1}  (z ^ w)^{-1}]^2)^g 
=1$ by our relations,  and  we again have $N \cong A_{n+2}$ by
(\ref{Carmichael An}).  

Clearly  $N \unlhd J$ and 
  $J/N \cong T$.  Then  $|J|=|A_{n+2}\times T |$, so that   $J
\cong A_{n+2}\times T$.   
\end{proof}

\begin{Examples} \rm
\label{AGL1 example}
(1) Let $p$ be an odd prime, $n=p+1$
and
$T=\SL(2,p)$. Then $T$ has a presentation with $2$ generators and $2$
relations
\cite{CR2},  while 
 $T_1$ can be generated by $2$ elements.  Thus, by the Lemma, 
$A_{p+3} \times \SL(2,p)$ has a presentation with $3$ generators
and $6$ relations (cf.  Examples~\ref{summary of special n}(\ref{three
Ap+3}) and \ref{more AGL1 examples}(\ref{Ap+3})).   
\smallskip

(2)
We can do somewhat better by taking $T$
to be   $\AGL(1, p)=P\rtimes T_1$, with $P$   cyclic of  odd  prime
order $p$ and $T_1$  cyclic of order $ p-1 $.  
By \cite{bn},  
if $r$ and $s$ are integers such that $ \F_p^* =\<r\>$
and $s(r-1)\equiv -1$ (mod $p$), then  
\begin{equation}
\label{Neumann}
T =\AGL(1, p) = \< a,b \mid a^p=b^{ p-1 }, (a^s)^b=a^{s-1} \>.
\end{equation}
This time  $|X_1|=1$,  and hence 
\emph{$A_{p+2} \times T$ has a presentation with
$3$ generators and $5$ relations}. 
\smallskip

(3) An example of a 2-homogeneous group that is not 2--transitive 
is the subgroup $T= \AGL(1, p)^{(2)}$ of index 2 in 
$\AGL(1, p)$ for a prime $p\equiv 3$ (mod 4), $p>3$.   This time 
$T=P\rtimes T_1$  with
$P$   cyclic of  order $p$ and $T_1$  cyclic of order $ (p-1 )/2$. 
By \cite{bn}, 
$$T= \AGL(1, p)^{(2)} = \< a,b \mid a^{p}=b^{( p-1 )/2},
(a^s)^b=a^{s-1}
\>  ,$$
 where this time  
 $ \F_p^*{}^2 =\<r\>$
and $s(r-1)\equiv - 1$ (mod $p$).
Once again  \emph{$A_{p+2} \times T$ has a presentation with
$3$ generators and $5$ relations}. 
\smallskip

(4) For future reference we note that, 
for any prime $p\equiv 3$ (mod 4) with  $p>3$,
$$
 \AGL(1, p)^{(2)}\times\Z_2 = \< a,b \mid a^{2p}=b^{( p-1 )/2},
(a^s)^b=a^{s-2} 
\>  ,
$$ 
where  ${s(r-1)\equiv -2}$ (mod $p$), and 
 $ \F_p^*{}^2 =\<r\>$;
since both $s$ and  $s+p$ both satisfy  the preceding congruence, we may
assume that  with
$s$ odd.  For, the presented group   $T$  surjects onto 
$ \AGL(1, p)^{(2)}\times\Z_2 $  by sending 
$a \to (a',z)$ and $b\to  (b',1),$ with $a',b'$ playing the roles of $a,b$
in (3),  and  $|z|=2$. 
Since $a^{2ps} = (a^{2ps})^b=(a^{2p})^{s-2} $, we have $a^{4 p}=1$. 
Then $(a^{2p})^s=  a^{2p}  $, so that $(a^{ps})^b  =a^{p(s-2)}=a^{-ps}$,
and hence 
$a^{ps}=(a^{ps})^{b ^{(p-1)/2}} = a^{-ps}$ 
since $(p-1)/2$ is odd. 
Now $a^{4p}=1=a^{2ps}$. 
Since $s$ is also odd, it follows that   $a^{2p}=1$, and  also that 
$\<a\>
\unlhd  T$.  Then
$a^{ p} \in Z(T)$, and
$T$ is as claimed.

\end{Examples}

\begin{corollary} 
\label{Ap+2 with 6}
If $p$ is prime then   
$A_{p+2}$ has a presentation with $3$ generators$,$ $6$ relations
and bit-length $O(\log p)$.  
\end{corollary}

\Proof
By Example~\ref{AGL1 example}(2), the lemma provides  a presentation for
$A_{p+2}\times
\AGL(1,p)$ with 3 generators
${\bf a} ,{\bf b },{\bf z}$. 
%
%
Under the isomorphism \eqn{def phi} in the lemma, 
in $ A_{p+2} \times T$ we have    $  {\bf b } =\big (b (p+1, p+2), b\big )  $
for a   $(p-1)$-cycle $b\in T _p $ (the point 1 in that lemma is now the
point $p\equiv 0$ fixed by $b$).  Then 
$b^a$ moves
$0$   (and fixes $p+1$ and
$p+2$), so that  $  {\bf b^a  z}    =\big (b^a (p+1, p+2), b^a\big )(z',1) =
(c ,b^a)$ for a 
$p$-cycle $c $. Thus, $(  {\bf b^a  z} )^p=(1,b^a  )$.  Since $T$ is the
normal  closure of $b^a$ in $T$,  imposing the relation $(  {\bf b^a  z} 
)^p=1$ gives a presentation for
$A_{p+2}.$~\qed
 
\medskip
Note that this is not, however,  a short presentation 
(cf.~Section~\ref{Preliminaries}).

    For some choices of
$p$ we can improve the preceding result:

\begin{corollary} 
\label{Ap+2 4 relations}
For any prime    \mbox{\rm $p\equiv 11$~$ ($mod $12),$}   
\begin{itemize}
\item[\rm(i)]
$A_{p+2}\times \AGL(1,p)^{(2)}$ has a presentation with $2$ generators$,$ $3$
relations and bit-length $O(\log p);$ and
\item[\rm(ii)] 
$A_{p+2}$ has a presentation with $2$ generators$,$ $4$ relations
and bit-length $O(\log p)$.  
\end{itemize}
\end{corollary} 

\Proof (i) Since $p\equiv 3$~$ ($mod $4),$ we can let $T=\AGL(1,p)^{(2)}\le
A_p,$ 
$r$ and
$s$ be as in Example~\ref{AGL1 example}(3).  Consider the group $J$ defined by
the presentation 
$$
  \langle a, g \mid   a^p=b^{ (p-1)/2 }, (a^s)^b=a^{s-1},
(zz^a)^2=1  \rangle,
$$
where $b:=g^{3}$ and $z:=g^{(p-1)/2}$. 

First note that $A_{p+2}\times T$ satisfies this presentation.  Namely, let
$T$ act on  $\{1,\ldots, p, p+1,p+2\}$, fixing
$ p+1$  and $p+2  $.  Let $g$  be the product of the 3-cycle $ (1 , p+1 ,
p+2)  $  and   
an element of $T$ having two  cycles of length
$(p-1)/2$ on $\{2,\ldots, p \}$.  Since 
 $p\equiv 2$~$($mod~$3)$, $g$~has order $3(p-1)$, so that 
$T=\langle a,g^3 \rangle$ satisfies the presentation in  Example~\ref{AGL1
example}(3).  The final relation $(zz^a)^2=1 $ holds as in Lemma~\ref{using
Carmichael} using $w=a$.

By Example~\ref{AGL1 example}(3) and Lemma~\ref{It's a subgroup},  
$J$ has a subgroup we can identify with $T=\langle a,   b \rangle.$
Clearly, $T_p =T_0 =\langle b\rangle$. 
The remaining relations   $z^3=1 $ and  $z^b =z$  in Lemma~\ref{using
Carmichael} are 
automatic: they hold in  $\< g\>$.  This proves (i).

(ii) This is similar to   the preceding
corollary.  This time   $b$ has two cycles of length $(p-1)/2$, and we obtain
\begin{equation}
\label{Ap+2} 
A_{p+2}    \cong  \langle a, g\! \mid \!   a^p\! =\! g^{ 3 \kappa  },
(a^s)^{g^{3}} \!\!  =a^{s-1},    \big (g^{\kappa } (g^{ \kappa  })^a\big)^2
 \hspace{-2pt}=\!  
 \big( g^{3}   (g^{\kappa })^a (g^{\kappa })^{a^{-1}}\big)^{\kappa+1}
\! =\! 1 \>%
\end{equation}
with  
$\kappa =(p-1)/2$,  $s(r-1)\equiv -1$ (mod $p$) and   $ \F_p^*{}^2
=\<r\>$, since  ${-1}$ is a non-square mod $p$.  For,   we view (i) as  a
presentation for
$A_{p+2}\times
\AGL(1,p)^{(2)}$ with  generators
${\bf a}$ and ${\bf g}$, and we use $ {\bf b }:={\bf g }^3$  and  ${\bf z} :={\bf
g}^\kappa  $. By  \eqn{def phi},
as in the proof of Corollary~\ref{Ap+2 with 6} we have
$  {\bf b  z^az^{a^{\mbox{$\scriptscriptstyle -1$}}}}    =\big (b, b\big )(
z'{}^az'{}^{a^{-1}},1) = (c_1c_2 ,b )$ for disjoint cycles $c_1$ and $c_2$
of length $\kappa  +1$, 
since
$a(1)$  and $a(-1 )$ are in
different 
$b$-cycles.  Then  $(  {\bf b  z^az^{a^{\mbox{$\scriptscriptstyle-1$}}}} )^{\kappa+1}=(1,b  )$, and    imposing the relation $(  {\bf b 
z^az^{a^{\mbox{$\scriptscriptstyle -1$}}}})^{\kappa+1}=1$ on the presentation
in (i) gives 
\eqn{Ap+2}.\qed 

\vspace{2pt}
\Remark  
\label{cycles for p+2}
 In the   presentations in the preceding corollaries$,$  every cycle
$(k,k+1,\dots,l)$ 
$($with
$k-l$ even$)$ can be written as a   word of bit-length $O(\log p)$  in the
generators. Any even permutation  with bounded support can
also be expressed as   word of bit-length $O(\log p)$  in the generators.
For all of these elements$,$  
the indicated words  use a bounded number of exponents.
\rm
 
Namely,  all $3$-cycles $(k,p+1,p+2) = (1,p+1,p+2)^{a^{k-1}}$  have
 bit-length $O(\log p)$,  therefore the same is true of any permutation of
bounded  support (because it is a product of 
a bounded number of such 3-cycles).
In particular,  if $x=(1,p)(p+1,p+2)$ then  
$xa=(1,\dots,p-1)(p+1,p+2)$ has  bit-length $O(\log p)$. Since
$$
(xa)^{-k}a^{k} = (p,1,\dots, k)(p+1,p+2)^{k},
$$
if  $k<p$ is even then the $(k+1)$-cycle $(p,1,\dots ,k)$ has  bit-length
$O(\log p)$, hence the same is true of all even cycles of the form
$(k,\dots,l)$, $l < p$.  The remaining cycles arise as, for example,
$(k,\dots,p)(k,p,p+1)=(k+1,\dots,p, p+1)$.  
 
 \medskip
\noindent
 {\bf Symmetric groups.}
 There  are   analogous results for 
symmetric groups.  This time we assume that  our group  $T=\langle X \mid
R\rangle$  acts   2-transitively on     
   $\{1,\ldots, n\}$ and  \emph{does not lie in}  $A_n$; let $w $   be
a   word in  $ X$ that moves $1$  when  $w$ is  viewed inside $T$. 
As above, the obvious examples   are $\AGL(1,p)$ and
$\PGL(2,p)$.   

\begin{lemma}  
\label{using Carmichael2-}
If $J = \langle X, z \mid  R, z^3=1, (zz^w)^2=1, [z,T_1]=1 
\rangle,$ then $J$ has a normal subgroup $T\cap A_{n+2}$ modulo which it is
$S_{n+2}$.  Moreover$,$ $J$ is isomorphic to a 
subgroup of index $2$ in $S_{p+2}\times T $ that projects onto
each factor.
\end{lemma}
 
\Proof
This time view $T$ as a subgroup of $S_{n+2}$ acting on    $\{1,
\ldots, n, n+1,n+2\}$ but fixing  $  n+1$  and  $n+2 $. 
Then  $S_{n+2} $ acts  on a set of size $n+2$,  while  
$ T$  also acts  on a set of size  $n$. 
View $S_{n+2}\times  T$ as acting on the disjoint union of these sets, 
and  let $H$ be the subgroup $ A_{2n+2}\cap (  S_{n+2}\times  T)$
of index 2 (recall that $T$ is not inside $A_n$).

This time we map $J$ into $H$ using  a simpler version of \eqn{def phi}:
$\varphi(x)=(x ,x)~$ for  $x\in X$.  We identify $T=\< X \>$  with a subgroup
of
$J$ and $z$ with an element  of $J$.  As before, 
$T$ acts on $z^T$ as it does on 
$\{1,\ldots, n\}$, and hence is 2-transitive.  Then   the relation 
$(zz^w)^2=1$ and  (\ref{Carmichael An}) imply that $N:=\<z^T\>\cong
A_{n+2}$. Since
$J/ N\cong T$, we have 
$|J| = |A_{n+2} | | T  |=|H|$, so that $J\cong H$.  
\qed

\begin{corollary}  
\label{Sp+2 4 relations}
Let $p$ be a prime.  
\begin{itemize}
\item[\rm(i)]
 $S_{p+2}$ has a
presentation with 
$3$ generators$,$ $6$  relations  and bit-length $O(\log p)$.  
\item[\rm(ii)]
If \mbox{\rm $p\equiv 2$ (mod 3)} then
the subgroup of index $2$ in $S_{p+2}\times \AGL(1,p) $ that projects onto
each factor  has a presentation with
$2$ generators$,$
$3$ relations and bit-length $O(\log p).$   
\item[\rm(iii)]
If \mbox{\rm $p\equiv 2$ (mod 3)} then
$S_{p+2}$ has a presentation with $2$ generators$,$ $4$ relations
and bit-length $O(\log p)$.  
\end{itemize} 
\end{corollary}
  
\Proof
Part (i)  follows from the preceding lemma together with Example~\ref{AGL1
example}(2),  while
 (ii)  is proved exactly as in Corollary~\ref{Ap+2 4 relations} by using 
that  example.  This time in  (iii)  we obtain
\begin{equation}
\label{Sp+2}
\begin{array}{lll}
 S_{p+2}  \cong  
 \langle a, g   \mid  &\hspace{-6pt} 
a^p = (g^{3})^{p-1 }, 
(a^s)^{g^{3}} 
\!=\!a^{s-1} , \vspace{2pt}
\\
  &\hspace{-6pt}  \big (g^{ p-1 } (g^{ p-1})^a\big)^2  
\hspace{-2pt} =
 \big(   g^{ 6  } 
(g^{p-1 })^{a^{-1}}(g^{p-1})^{a^{-r}}\big)^{(p+1)/2}
 \hspace{-2pt}=1 \>  
\end{array}  
\end{equation} 
with   
$s(r-1)\equiv -1$ (mod $p$  and  $ \F_p^* =\<r\>$.   
 For, once again  we view (ii) as  a presentation   with  generators
${\bf a}$ and ${\bf g}$,  we use $ {\bf b }:={\bf g }^3$  and  ${\bf z}:={\bf
g}^{p-1}$, and then 
$  {\bf b } ^2 {\bf z^{a^{-1}}   z^{a^{\mbox{$\scriptscriptstyle -r$}} } } =\big (b^2, b^2 \big )(
z'{}^{a^{-1}}z'{}^{a^{-r}},1) = (c  ,b^2 )$ for  a $(p+1)/2$-cycle  $c $.
Factoring out the 
normal  closure of $b^{p+1} =b^2$ in $T$, we obtain \eqn{Sp+2}.\qed 
  
\vspace{2pt}
\Remark  
\label{cycles for p+2 again}  \hspace{-5,5pt}
For the presentation in {\rm Corollary~\ref{Sp+2 4 relations},}\hspace{-.5pt}
the assertions in  {\rm Remark~\ref{cycles for p+2}}
 hold  once again for even permutations.  \hspace{-6pt} They also hold for odd permutations     if $p\equiv 11 $ \rm (mod 12). 

Even permutations are handled as before.
Odd permutations are slightly more delicate: they require constructing a transposition of  the required bit-length, and we are only able to achieve this when $p\equiv 11 $ \rm (mod 12).   
First we recall  a group-theoretic version of  ``Horner's Rule''   \cite[(3.3)] {GKKL}
for  elements $v,f$ in any group:
 \begin{equation}
\label{Horner}
vv^fv^{f^2}\cdots v^{f^n}=(vf^{-1})^nvf^n
\end{equation}
for any positive integer $n$.

Note that 
\vspace{2pt}

$
b_2:={\bf b} ^{(p-1)/2 } =  (1,p-1)(2,p-2)\cdots \big(( p -1)/2,p-(p-1)/2
\big)
$

\vspace{2pt}
\noindent
is an odd permutation since $p\equiv 3$   (mod 4).  We will use several additional permutations:

$
\mbox{$c({i, j}):=    {\bf z} ^{{\bf a}^{-i}}(  {\bf z}^{{\bf a} ^{-j}}  ) ^{-1}  {\bf z}  ^{{\bf a}^{-i}} = (i,j)(p+1,p+2)  $  for $1\le i,j\le p$, 
}$

\vspace{2pt}
$v_\bullet :=   c({1,p-1})c({2, p-2})  = (1, p-1)(2,p-2),$

\vspace{2pt}
$c_{(p-1)/2}:= \big (c({1,2}){\bf a} \big)^{{(p-1)/2}-2}  c(1,2) {\bf a} ^{-({(p-1)/2}- 2 )} = \big (1,2,\dots,{(p-1)/2} \big ) $

\vspace{2pt}
\noindent
 since $(p-1)/2$ is odd (cf. (\ref{Horner})),
 
\vspace{2pt}
$c_{\bullet}:= c_{(p-1)/2}^{~{}^{~}}\,  c_{(p-1)/2}^{-{\bf a}^{(p+1)/2}}  =  \big(1,2,\dots,(p-1)/2 \big ) \big (p-1,p-2,\dots, p-(p-1)/2 \big)$, and 

  $v : = \big (c(1,p-1) c_{\bullet}  ^{-1} \big )^{(p-1)/2}c(1,p-1)c_{\bullet}^{(p-1)/2}  \!
 = \big (c(1,p-1) c_{\bullet}  ^{-1} \big )^{(p-1)/2}c(1,p-1)$%
\vspace{2pt} 

\hspace{8pt} $ \equiv 
(1,p-1)(2,p-2)\cdots \big ((p-1)/2,p-(p-1)/2\big ) (p+1,p+2)~$   (cf.
(\ref{Horner})).
\vspace{2pt} 

\noindent
Then 
 $v b_2 $ is the transposition $(p+1,p+2)$,  as required 
 (compare Section~\ref{An explicit presentation for $S_n$}). 
\vspace{2pt} 
\vspace{2pt} 
\vspace{2pt} 

We do not know if the final assertion in the preceding remark holds when $p\equiv 1$ (mod 4).
§
\vspace{2pt}
\Examples  
\label{summary of special n}  
We summarize versions of the previous presentations for
special values of $n$. \rm  Additional  explicit  presentations appear
in Examples~\ref{more AGL1 examples}.  

Consider  a prime  $p>3$. 
{
\newcounter{primes}
\begin{list}
{\rmfamily\upshape (\arabic{primes}) }
{\usecounter{primes}
\setlength{\labelwidth}{0cm}\setlength{\leftmargin}{0cm}
\setlength{\labelsep}{0cm}\setlength{\rightmargin}{0cm}
\setlength{\parsep}{0.5ex plus0.2ex}
\setlength{\itemsep}{0.5ex plus0.2ex} \upshape}
\vspace{2pt}
\item  $A_{p+2} \! =  \!\< a,b,z   \!\mid \! a^p\!=b^{ p-1 },
(a^s)^b=a^{s-1},
 z^3\!= (zz^a)^2\!=1, 
z^b=z^{-1} ,
(  { b^a  z})^p\!=1\>$,

where  $s(r-1)\equiv -1$ (mod $p$) with  $ \F_p^* =\<r\>$ 
(by Corollary~\ref{Ap+2 with 6}).

\item   $A_{p+2}$ for $ p\equiv 11$  (mod 12): see  (\ref{Ap+2}).

\item   $S_{p+2}$ for $ p\equiv 2$  (mod 3): see  (\ref{Sp+2}).

\item   $S_{p+2}\!=\!\< a,b,z   \!\mid \! a^p\!=b^{ p-1 },
(a^s)^b=a^{s-1},
 z^3\!=  \hspace{-1pt} (zz^a)^2\!=\![z,b]\!=1,  
 (b^2   z^{a^{-1}} \!z^{a^{-r}}   \hspace{-1pt}   )^{p}\!
=\!1\>$

for $ p\equiv 1$  (mod 3), 
where  $s(r-1)\equiv -1$ (mod $p$) and  $ \F_p^* =\<r\>$ 
(using Corollary~\ref{Sp+2 4 relations}).  
 
\item  
\label{three Ap+3}
  We will give several presentations of
$A_{p+3}\hspace{-1pt}$ both here and   in Example~\ref{more AGL1
examples}(\ref{Ap+3}).  Let $ \F_p^* =\<j\>$  and  $j\bar
\jmath\equiv1\hspace{-1pt}$  (mod
$p$).~Then
$$
\begin{array}{lll}
A_{p+3} = \< x,y, z  \mid   &\hspace{-8pt} 
 x^2 = (xy)^3, (xy^4xy^{(p+1)/2  })^2y^px^{2[p/3]} =
1, \vspace{2pt}
\\
  &\hspace{-8pt}   
z^3\hspace{.8pt}=  (zz^x)^2=  [y,z]=[h ,z]=1,
  (hz^{(xy)^{-1}}z^{{(xy^j)^{-1}}})^{(p+1)/2}  =1 \> ,
\end{array}  
$$
where we have abbreviated  
$h:= y^{\bar \jmath}  (y^{  j})^x y^{\bar
\jmath}x^{-1}. $  
This uses the following presentation for $T:=\SL(2,p)$,   obtained in
\cite{CR2} using \cite{Sun}:
\begin{equation}
\label{efficient SL2}
\SL(2,p)=\<x,y\mid x^2 = (xy)^3, (xy^4xy^{(p+1)/2  })^2y^px^{2[p/3]} =
1\> ,
\end{equation}
where   $x$ and $y$ arise from  elements of order 4 and $p$, respectively.
(These correspond to  the matrices $t$ and $u$ given later in
(\ref{generators}).) Then $T_1=\<X_1\>$ with $X_1:=\{y,h  \}$ in  the
notation used in Lemma~\ref{using Carmichael};  the final relation
in the presentation for $A_ { p+3}$   is
obtained as in  the proof of 
Corollary~\ref{Ap+2 4 relations}(ii).  
$$
\begin{array}{lll}
A_{p+3} = \< u,  h, t,z  \mid   &\hspace{-8pt} 
u^p= t^2 =1, u^h=u^{j^2} , h^t = h^{-1},t =u u^t u,h t=
u^{\bar \jmath}  (u^{  j})^t u^{\bar \jmath}, \vspace{2pt}
\\
  &\hspace{-8pt}   
z^3\hspace{.8pt}=  (zz^t)^2=  [u,z]=[h,z]=1,
  (hz^{(tu)^{-1}}z^{{(tu^j)^{-1}}})^{(p+1)/2}  =1 \> .
\end{array}  
$$
 This uses
Lemma~\ref{using Carmichael} together with the presentation for
$T:=\PSL(2,p)$ given in
\cite[Theorem~3.6]{GKKL}.  
A similar presentation can be obtained  using the presentation for $\PSL(2,p)$  in \cite{To}.

\item  Once again let  $ \F_p^* =\<j\>$.   Then
$$ 
\begin{array}{lll}
S_{p+3} = \< u,  h, t,z  \mid   &\hspace{-8pt} 
u^p= t^2 =1, u^h=u^{j} , h^t = h^{-1},t =u u^t u, \vspace{2pt}
\\
  &\hspace{-8pt}   
z^3\hspace{.8pt} =  (zz^t)^2=  [u,z]=[h,z]=1,
  (hz^ {tu})^{p +1}  =1 \>.
\end{array}  
$$
   This uses
Lemma~\ref{using Carmichael} together with a presentation
$\< u,    h, t   \mid 
u^p= t^2 =1, u^h=u^{j} , h^t = h^{-1},t =u u^t u\>
$
 for
$T:=\PGL(2,p)$ analogous to \cite[Theorem 3.6]{GKKL}. The final
relation is obtained as in  the proof of 
Corollary~\ref{Ap+2 4 relations}(iii).  
Once again one could also  use \cite{To} for a presentation of $T$ .
\end{list}
   }


\subsection{Small $n$}
\label{Small n}
In order to handle a few   degrees  $n<50$ 
we will need further variations on the idea used  
in Corollaries~\ref{Ap+2 4 relations}  and  \ref{Sp+2 4
relations}.
All of the general presentations below have bit-length $O(\log n)$, but this is not significant since
our goal involves bounded $n$.  
We suspect that most readers will wish to skip this
section. 
 
 In  Table~\ref{Small degrees}   we summarize the cases needed later.  
For this table  and our variation on Corollaries~\ref{Ap+2 4 relations} 
and  \ref{Sp+2 4 relations}    we use the following
notation:

 \vspace{5pt} 
%
 \begin{table}[h]
 \caption{Some small $n$}
\label{Small degrees}
\begin{tabular}{|l|l|l|c|c|c|c|c|c|l|}
\hline 
$G$ 				& \hspace{8pt} $n$  &\hspace{17pt}  $T  $       & $|R|$ & $\rho$ &
$|X_1|$ & $|x_1|$ &gens&  rels& in Ex.$\, \#$  \\
\hline \hline
$ S_{11}$& $9+2$           &$\AGL(1,9) $ &4  & 1 & 1 & 8  &2 &6  &     
\mbox{\rm \ref{more AGL1 examples}(\ref{AGL(1,9)})}
\\ 
\hline 
$ A_{11}$&$ 9+2 $   &$\PSL(2,8)$%
&2  &
2 & 1	& 4 & 2 &  5    &    
\mbox{\rm \ref{more AGL1 examples}(\ref{PSL(2,8)})}
\\ 
\hline 
$ A_{11}$&$ 9+2 $      \raisebox{2.5ex}{\hspace{-1pt}}    &$\AGL(1,9)^{(2)}$%
&4  &
2 & 1	& 4 & 2 &  7    &    
\mbox{\rm  \ref{more AGL1 examples}(\ref{AGL(1,9)/2})}
\\ 
\hline 
$ S_{12} $& $10+2 $          &$\PGL(2,9) $	&3  & 1 & 2 & 8 &2 & 6 &     
\mbox{\rm \ref{more AGL1 examples}(\ref{PGL(2,9)})}
\\

\hline
$A_{12}$& $10+2 $            &$6.\PSL(2,9)$ 	&2  & 1 & 2 & 8 &2 &  5
&\mbox{\rm \ref{more AGL1 examples}(\ref{6PSL(2,9)})} 
\\ 
\hline 
$ A_{23}$ &$21+2$ 
\raisebox{2.5ex}{\hspace{-1pt}}\raisebox{-1.2ex}{\hspace{-1pt}} 
& 12.\PSL(3,4)
& 2 &1 & 2 & 5  &  2 & 5
&\mbox{\rm  \ref{more AGL1 examples}(\ref{PSL(3,4)})}
\\
\hline 
$ A_{23}$ &$\binom{7}{2}+2$  
\raisebox{2.7ex}{\hspace{-1pt}}\raisebox{-1.4ex}{\hspace{-1pt}} 
&  $6.A_{7} $   &2 &2 & 2 & 5  &  2 & 6  &   
\mbox{\rm \ref{more AGL1 examples}(\ref{6A7})}
\\
\hline   
 $A_{24}$& $22+2  $  &$12.M_{22}  $
& 2 & 1 & 2 &7 &2 &
5   &
\mbox{\rm \ref{more AGL1 examples}(\ref{M22})}
\\
\hline   
 $A_{24}$& $2\cdot11+2  $  &$\AGL(1,11)^{(2)}\times \Z_2  $
& 2 & 3 & 1 &5 &2 &
6   &
\mbox{\rm  \ref{more AGL1 examples}(\ref{A2p+2})}
\\
\hline 
 $A_{24}$&$ 2\cdot11+2   $ &$\AGL(1,11)  $%
 & 2 & 4 & 1 &5  &2 & 7    &
\mbox{\rm \ref{more AGL1 examples}(\ref{A2p+2 again})}
\\
\hline 
$ A_{47}$&$ \binom{10}{2} +2$
\raisebox{2.7ex}{\hspace{-1pt}}\raisebox{-1.4ex}{\hspace{-1pt}}
& \!$A_{10}$  \!\!\!& 2 & 2 & 2 & 8 &2 & 6    &   
\mbox{\rm \ref{more AGL1 examples}(\ref{A10})} 
\\
\hline 
$ A_{47}$&$ \binom{10}{2} +2$
\raisebox{2.7ex}{\hspace{-1pt}}\raisebox{-1.4ex}{\hspace{-1pt}}
& \!$A_{10}\!\times\!\SL(2,7)$ \!\!\!& 4 & 2 & 2 & 7  &2 & 8   &    
\mbox{\rm \ref{more AGL1 examples}(\ref{p+3 choose 2})}
\\
\hline 
 $A_{48}$&$ 2\cdot23+2   $ &$\AGL(1,23)^{(2)}\times \Z_2 $%
 & 2 & 3 & 1 &11  &2 & 6    &
\mbox{\rm \ref{more AGL1 examples}(\ref{A2p+2})}
\\
\hline 
 $A_{48}$&$ 2\cdot23+2   $ &$\AGL(1,23)  $%
 & 2 & 4 & 1 &11  &2 & 7    &
\mbox{\rm \ref{more AGL1 examples}(\ref{A2p+2 again})}
\\
\hline 
\end{tabular}
 \end{table}
%
%
%


\begin{itemize}
\item
$T$ is   a group acting transitively
$($though not necessarily faithfully$)$ on $\{1,\dots ,n\}$.
\item
$T$ has  exactly
$\rho$ orbits of unordered pairs of distinct   points.
\item
$T=\<X \mid R\>$. 
\item
$T_1=\<  X_1\>$, where $T_1 $ is again the stabilizer of 1.
\item
 $x_1\in
X_1\cap X$ has order
$k$  not divisible by $3$.
\item
  $T$ is also viewed as a subgroup of $\Alt\{1,\dots ,n+2\}$.
\vspace{2pt}
\end{itemize}

Thus,  $\<T^{S_{n+2}} \>$ is $A_{n+2}$ if $T$ is   in  $A_{n+2}$, and 
$S_{n+2}$ otherwise.

\begin{proposition}  
\label{general number of relations}  
If $T$ is the normal closure of one of its elements$,$ then
$\<T^{S_{n+2}} 
\>$ has a presentation  with $|X| $ generators and $|R|+\rho+|X_1| $
relations.

If $T\le A_n$ then $ A_{n+2} \times T$ has a
presentation  with $|X| $ generators and $|R|+\rho+|X_1| -1$ relations.
\end{proposition}

{\noindent \em Proof sketch.}
Let   $w_1,\dots,w_\rho $ be words in $X$ such   that  the $\rho$ pairs  
 $\{1, w_i^{-1}(1)\}$ are  in different $T$-orbits.
Let $X=\{x_1,\dots\}$ and $R=\{r_1,\dots\}$ with each
$r_i$ a word $r_i(x_1,\dots)$ in $X$.   
Let $X ':=X \backslash\{x_1\}$  and  $X_1':=X_1\backslash\{x_1\}$. We claim
that  
$$
\begin{array}{llll}
J:=\<X', g  \mid & \hspace{-7pt} 
r_i(g^3,\dots)=1 \mbox{ for all $i$}, \vspace{2pt}
\\
&  \hspace{-7pt}
[g^{k},X_1']= \big(g^{k}(g^{k})^{w_j}  \big)^2=1
 \mbox{ for all $j$}\> 
\end{array}
$$  
is isomorphic to  $ H:=A_{2n+2}\cap (  S_{n+2}\times  T)$.
For, view $S_{n+2}\times  T$ as
acting on the disjoint union
$\{1,\dots ,n,n+1,n+2\}\dot \cup  \{1,\dots ,n\}$  with $T$ fixing $n+1$ and $n+2$.
Let $z'\! : =(1,n+1,n+2)$ and
let the integer $\nu $ satisfy  $3\nu \equiv 1$ (mod $k$), so that 
 $g:=x_1^\nu z' $ satisfies $g^3=x_1$.  Then  $J$ surjects onto
$H$ as in the proof of Lemma~\ref{using Carmichael}.

We can identify  $\tilde T:
=\<g^3, X' \>$ with a subgroup of $J$ that is  a
homomorphic image of our original $T$. 
Our hypotheses guarantee that $z:=g^k$ commutes with $X_1$.  
Then $|z^{\tilde T}|=n$, and we can
use \eqn{Carmichael An}  as before.
\qed
 
\medskip
Examples~\ref{more AGL1
examples}(\ref{A2p+2}) and \ref{more AGL1
examples}(\ref{A2p+2 again})   contain  further  
variations on the idea behind the
Proposition.

\begin{Examples} 
\label{more AGL1 examples}  \rm

{
\newcounter{Ex}
\begin{list}
{\rmfamily\upshape (\arabic{Ex}) }
{\usecounter{Ex}
\setlength{\labelwidth}{0cm}\setlength{\leftmargin}{0cm}
\setlength{\labelsep}{0cm}\setlength{\rightmargin}{0cm}
\setlength{\parsep}{0.5ex plus0.2ex}
\setlength{\itemsep}{0.5ex plus0.2ex} \upshape}
\item  
\label{AGL(1,9)}
$n=11$:~ $\,T=\AGL(1,9)$ has the presentation   
 $ \< a,b \mid a^3=b^8=1,  a^{b^2}=aa^{-b }, [a,a^b ]=1
\>$, 
$\rho= |X_1|=1$ and $|x_1|=8$, so that 
\emph{$S_{11}$ has a presentation with $2$ generators and $ 4+1+1 $
relations.}   However, for use in Theorem~\ref{C}  it is  easier simply 
to use the presentation of $S_{11}$ with 2 generators  and 6 relations in
\cite[p.~54]{Ar} (cf.   \cite[p.~64]{CoMo}).
See Remark~\ref{degree q+2} for a generalization.  

\item \label{PSL(2,8)}
 $n=11$:~  $T=\PSL(2,8)$  
has a presentation with 2 generators and 2
relations \cite{CHRR}, $\rho=1,|X_1|=2$  and  $ |x_1|=7$, so that  
Proposition~\ref{general number of relations} produces a
\emph{presentation of
$A_{11}$ with $2$ generators and $2+2+1$ relations.}
Once again, it has long been known that   $A_{11}$  has a presentation  with 2
generators  and 6 relations  \cite[p.~67]{CoMo}.

\item  
\label{AGL(1,9)/2}
$n=11$:~ $T$ has index 2 in $\AGL(1,9)$, $T$  has the presentation   
 $ \< a,b \mid a^3=b^4=1, a^{b^2}=a^{-1},
 [a,a^b ]=1 \>$,   
$\rho=2$, $|X_1|=1$ and $|x_1|=4$, so that 
\emph{$A_{11}$ has a presentation with $2$ generators and $4+2+1 $
relations.}
See Remark~\ref{degree q+2} for a generalization.  
\item  
\label{PGL(2,9)}
$n=12$:~ $T=\PGL(2,9)$ has   presentations  with 2 generators and $3$
relations provided by G. Havas \cite{Hav}.  The following are some of his many
presentations related to $A_6$:
$$
\label{Havas}
\begin{array}{rlll}
\PGL(2,9)& \hspace{-6pt} =&\hspace{-6pt}\<  a, b \mid b^5=(aba)^2= 
ab^{-1}a^4ab^{-1}a^3b^{-1}ab^2=1  
\>%
 \vspace{2pt}
\\ 
\PGL(2,9)& \hspace{-6pt} =&\hspace{-6pt}\<  a, b \mid  b ^8 = baba^3bab^2=
b^{-1}  a^{-1}  b^{-2} a^2ba^2b^{-1}  =1\> 
 \vspace{2pt}
\\
S_6&\hspace{-6pt}=&\hspace{-6pt}\< a, b \mid
a^{-1}b^{-1}a^3b^{-1}a^{-2}=(ab)^5= b^{-3}a^{-1}b^{-1}a^2b^{-2}aba =1\> 
 \vspace{2pt}
\\
S_6&\hspace{-6pt}=&\hspace{-6pt}\< a, b \mid
 a ^4= b^{-1}ab^{-2}a^2b^2a^{-1} b^{-1}= aba^{-1} b^{-1}a^{-2} b^{-1}
aba^{-1} b=1\>.
\end{array}
$$  
This time $\rho=1$, $|X_1|=2$  and  we may assume that    $|x_1|=8$, so that
  \emph{$S_{12}$  has a presentation with $2$ generators and  $3+1+2$ 
relations.}  Once again, it has long been known that  $S_{12}$  has a
presentation  with 2 generators   and 7 relations in
\cite[p.~54]{Ar} (cf.   \cite[p.~64]{CoMo}).

\item  
\label{6PSL(2,9)}
$n=12$, $T \cong 6.\PSL(2,9)\cong 6.A_6$ has a presentation with 2 generators
and 2 relations \cite{Ro},   
$\rho=1$, $|X_1|=2$ and   $|x_1|=8$, so that 
\emph{$A_{12}$ has a presentation with  $2$  generators and $2+1+2  $
relations.}   Once again, it has long been known that   $A_{12}$ has a
presentation  with 2 generators  and 7 relations  \cite[p.~67]{CoMo}.

\item  
\label{PSL(3,4)}
 $n=23$:  $ T=12.\PSL(3,4)\!$ has a presentation with 2 generators and 2
relations \cite{CHRR}, $\rho=1,|X_1|=2$  and  $ |x_1|=5$, so that  
Proposition~\ref{general number of relations} produces a
\emph{presentation of
$A_{23}$ with $2$ generators and $2+2+1$ relations.}

\item  
\label{6A7}
$n\!=\!23$:~$T=6.A_7\!$ has  a presentation with 2 generators and 2  relations
\cite{CRKMW},
$n=   \binom{7}{2} $, $\rho=2$, $|X_1|=2$   
and $|x_1|=5$, so that 
\emph{$A_{23}$ has a presentation with $2$ generators and $2+2+2  $ 
relations.}

\item  
\label{M22}
$n=24$:~  $ T=12.M_{22}$ has a presentation with 2 generators and 2
relations \cite{CHRR}, $\rho=1,|X_1|=2$  and  $ |x_1|=7$, so that  
Proposition~\ref{general number of relations} produces a
presentation of
\emph{$A_{24}$ with $2$ generators and $2+2+1$ relations.}

\item  
\label{Ap+3}
{\em If   $p>3$ is prime then  $A_{p+3}\times \SL(2,p)$ has a
presentation with 
$2$ generators and $4$  relations.}  Namely,  
apply Proposition~\ref{general number of relations} 
using  (\ref{efficient SL2})  and $|R|=2 =|X_1|,$
$\rho=1$,  $x_1=y$. 
It follows that {\em  $A_{p+3} $
has a presentation with 
$2$ generators and $5$  relations}.  
We will need this below in (\ref{p+3 choose 2}).

\quad
Explicitly, as in 
Example~\ref{summary of special n}(\ref{three Ap+3})  we have 
$$
\begin{array}{llll}
A_{p+3} =\<x, g  \mid & \hspace{-7pt} 
x^2 = (x g^3)^3, (xg^{12}xg^{3(p+1)/2 
})^2 g^{3p} x^{2[p/3]} = 1, \vspace{2pt}
\\
&  \hspace{-7pt}
[g^{p},h]= \big(g^{p}(g^{p})^{x}  \big)^{\!2}=1,
   (h( g^p)^{x{g^3}}
   )^{p+1}  =1 \> ,
\end{array}
$$  
where 
$h:= {g^3}^{\bar \jmath}  ({g^3}^{  j})^x {g^3}^{\bar\jmath}x^{-1}$   with
$ \F_p^* =\<j\>$ and 
 $j\bar
\jmath\equiv1\hspace{-1pt}$  (mod
$p$).

\item  
\label{A10}

$n=47$: $T= A_{10} $ has  the following  presentation with 2 generators and 3 
relations
\cite{Hav}:
$$
A_{10} \!=  \! \<a, b \!\mid \!  a^{3}b^{-1 }a b^{-1 }\hspace{-.2pt}a^3ba^2b
=\hspace{-.5pt}  a^2 b^{-1}\hspace{-.2pt}  a^5 b^{-3} a^3\hspace{-1pt} =
\hspace{-1pt} a^{-2}  ba  b^{-1}\hspace{-.2pt} aba^3 ba  b^{-1} ab a^{-2}  
b^{-1}\hspace{-1pt}  =1\>,
$$
with   $|a| =15,  |b|= 12$ and     $|ab|=8$; hence we modify this
presentation so that the generators  are $a$ and $ab$.  View  $T$ as acting
on $ \binom{10}{2} $ unordered pairs with 
$\rho=2$,
$|X_1|=2$  and $  x_1  = ab$, so that 
Proposition~\ref{general number of relations}  produces  a
\emph{presentation of $A_{45+2}$  with  $2$
generators and $2+2+2  $ relations.}
 
\item  
\label{p+3 choose 2}
{\em If   $p>3$ is prime then $ A_{\binom{p+3}{2}+2}$ has a presentation 
\vspace{2pt}
with $2$ generators and $8$ relations.} For,   let $T=A_{p+3}\times
\SL(2,p)$ act on $ \binom{p+3}{2} $ unordered pairs of a set of size $p+3$,
with $ \SL(2,p)$ acting trivially. Apply Proposition~\ref{general number of
relations} using (\ref{Ap+3}), with $|X|=3
$,
$|R|=4$,   $\rho=2,$ 
$|X_1|=2$  and  
$x_1=y$.  (Note that $x_1$ fixes $p+2$ and $p+3$ in (\ref{Ap+3}).)

\quad
There is a similar presentation for  
$ A_{\binom{p+2}{2}+2}$.


\item  
\label{A2p+2}             
{\em For any prime}  $p \equiv 11$ (mod 12), $A_{2p+2}$ {\em  has
a  presentation with $2 $ generators and $6 $ relations}.   We will vary
the argument in Proposition~\ref{general number of relations} (and
Lemma~\ref{using Carmichael}), using  the transitive subgroup 
$T:=\AGL(1,p)^{(2)} \times  \<t\>$  of   the transitive group
$\AGL(1,p)\times  \<t\>$ of degree $2p$, where    $t$  is  an involution 
interchanging two blocks of size $p$.  Note that 
    the stabilizer of a point is cyclic of odd order $(p-1)/2$. 
        Moreover,  $T$ has $\rho=3$  orbits of unordered pairs of the
$2p$-set,   with orbit-representatives as follows: contained in a block,
or of the form $\{  p,t(p)\}$, 
    or of the form $\{  p,t(1)\}$  (since $t$  interchanges    $\{ 
p,t(1)\}$    and $\{  1,t(p)\}$).  
           
\quad      
    We view $T$ as a subgroup of $A_{2p+2}$ preserving the two  new
points $2p+1$ and $2p+2$, with  $t$ interchanging these points.
   
\quad
 Let 
$$
J : =\langle a  , g \mid   a^{2p }=b^{(p-1)/2 }, (a^s)^b=a^{s-2},
  (zz^{ {\rm sign} (w_i) w_i})^2=1 ~( i=1,2,3)\rangle 
$$
where   ${s(r-1)\equiv -2}$ (mod $p$)
with $s$ odd  and $r$  of order $(p-1)/2$ (mod $p$),
$b:=g^3$,  
$z:=g^{(p-1)/2 }$,  and suitable words
$w_1,w_2,w_3\in T$ 
(such that  the pairs $\{z,z^{w_i}\}$ are  in different $T$-orbits  on the $2p$-set $z^T$);
here sign refers to the behavior on the $2p$ points. 
For example,    $\{z,z^{t}\}$, $\{{}z,z^{a^{2}}{} \}$ and $\{ {}
z,z^{a^{2}t} {}
\}$  are  representatives of these  $T$-orbits. 
As in the proof of Proposition~\ref{general
number of relations}, using  Example~\ref{AGL1 example}(4) we see that 
$A_{2p+2}\times  T$ satisfies our presentation.
 
\quad  
As usual, using   Example~\ref{AGL1 example}(4) we can view 
$T$  as the subgroup $\<a,b\>$  of $J$. 
Exactly as in Lemma~\ref{using Carmichael} (and Proposition~\ref{general
number of relations}),  $|z|=3$, $|z^T|=2p$,  and 
 $(zz^g)^2$ for all $g\in  \AGL(1,p)^{(2)} $  with $z^g\ne z $
 while
 $(zz^{-gt})^2=1$ for all $g\in  \AGL(1,p)^{(2)} $.

\quad
Then   $N:=\! \<z^T\>\cong A_{2p+2}$  by  \eqn{Carmichael An},\!    and
hence  
 $ J\! =NT\cong A_{2p+2}\times \AGL(1,p)^{(2)} \times  \<t\>$.   One
further relation gives  a presentation of
$A_{2p+2}$ with 3 generators and 
$2+ 3 +1$ relations.  

\quad   \mbox{Explicitly: } 
$$
\begin{array}{llll}
A_{2p+2} =\<a,  g  \mid & \hspace{-7pt} 
  a^{2p }=b^{(p-1)/2 }, (a^s)^b=a^{s-2}, \vspace{2pt}
\\ 
&  \hspace{-7pt}
(zz^{-t})^2= (zz^{a^{2}})^2= (zz^{-a^{2}t})^2=1,  
(t b  z^{a^{-1}}z^{a^{}}z^{ ta^{-1}}z^{ta^{}}  )^{p}=1   \> , 
\end{array}
$$
with  $z$, $r$ and $s$ as above and once again
$b:=g^3$, $z:=g^{(p-1)/2 }$ and
$t : =a^{ p} $,  where the 
last relation  is obtained as in  the proof of 
Corollary~\ref{Ap+2 4 relations}(ii).

\item  
\label{A2p+2 again}
 {\em For any odd prime}  $p \equiv 2$ (mod 3), $A_{2p+2}$ {\em  has
a  presentation with $2 $ generators and  $7$ relations.}  One difference 
between this example  and the preceding one is that we now handle the
case $p \equiv 1$ (mod 4)   using       $T=\AGL(1,p)$.  

\quad
First note that $T$ acts transitively on a set
of size $2p$, with cyclic  point stabilizer of order $ (p-1)/2$ and $\rho = 4 $ orbits on 
unordered pairs of points.  We again  view $T$ as a subgroup of
$A_{2p+2}$ preserving  the additional  points $2p+1$ and $2p+2$.  Once
again signs will refer to the actions of elements of $T$ on the $2p$-set.

\quad
We replace the presentation of $T$ in Example~\ref{AGL1 example}(1)  by
$$
T=\<x,b \mid  (xb^{-2})^p=b^{p-1},  ~((xb^{-2})^s)^b=(xb^{-2})^{s-2}\>,
$$
with $x:=ab^2\in T$ of order  $(p-1)/2$   fixing  the point $r':=(1-r^2)^{-1} $, 
 and  $r$ and $s$ as in Example~\ref{AGL1 example}(1).

\quad
Consider the group 
      $$J:=\!\<g,b \!\mid\! (xb^{-2})^p \hspace{-1pt} =
\hspace{-1pt} b^{p-1},  
((xb^{-2})^s)^b=(xb^{-2})^{s-2} ,    (zz^{ {\rm sign }(w_i ) w_i})^2
\!=1 \, (i=1,2,3,4) \>, $$ 
      where
$x:=g^3$, $z: = g^{(p-1)/2}$, and $\{r', w_i ^{-1} (r')\}$, $1\le i\le 4$,   
are representatives for the orbits of $T$ on pairs of the $2p$ points. 
Then
$ J$ surjects onto $A_{2p+2 } \times T$, and we can view $T=\< x,b  \>\le
J$. In particular,  $ |g|=3(p-1)/2$ and so $z^3=1.$

\quad
Since $x $  centralizes $z$, as usual we obtain   $ |z^T| = 2p $ and 
$N:=\< z^T\>\cong A_{2p+2}$ by (\ref{Carmichael An}),  and then   
$J\cong A_{2p+2 } \times T$. One further relation
produces the  desired presentation.

\end{list}
   }
\end{Examples}


\subsection{Gluing alternating and symmetric groups}
\label{Gluing alternating and symmetric groups}  
We now turn to  the case of all alternating and symmetric groups, starting
with a general gluing lemma:

\begin{lemma} 
\label{Sym m+n-k}
Let $G =\< X\mid R \>$ and  
 $\bar G =\<\bar X\mid  \bar R \>$
  be presentations of   $S_n$ and $S_m,$ respectively$,$ and let  
$m, n > k \ge l + 2 \ge 4$. Consider embeddings $\pi\colon G\to S_{m+n-k}$
and $\bar \pi\colon \bar G\to S_{m+n-k}$ 
into $ \Sym \{ -m+k+1, \dots ,n \}$ such that  
$$ \mbox{$\pi(G ) =\Sym(\{1,\dots , n\} )~$   and  
$~\bar\pi(\bar G ) =\Sym(\{-m + 1 + k, \dots , k\}).$}
$$
Suppose that $a, b, c, d\in G$   and   $ \bar a, \bar b, \bar c, \bar e\in
\bar G,$  viewed as  words in $X$ or $\bar X,$ respectively$,$
are nontrivial permutations such that 
\smallskip
\begin{itemize}  
\item $ \pi(a )= \bar \pi(\bar a )\in 
\Sym(\{1, 
\dots, l\})     <  \pi( G ) \cap  \bar\pi( \bar G) ,$
\item   
$ \pi(b )= \bar \pi(\bar b )\in 
\Sym({l + 1, \dots , k})   <  \pi( G ) \cap  \bar\pi( \bar G) ,$
\item $ \pi(c )= \bar \pi(\bar c )\in  \Sym(\{1,  \dots ,
k\})   <  \pi( G ) \cap  \bar\pi( \bar G) ,$
\item  $\pi ( d ) \in  \Sym(\{l + 1, \dots, n\})    <  \pi( G ) ,$
\item  $\bar\pi ( \bar e)  \in  \Sym(\{-m + 1 + k,  \dots , l\})    <  
\bar\pi( \bar G) ,$
\item  $\< \pi(a), \pi(c)\>  = \Sym(\{1,  \dots , k\})    <  \pi( G ) \cap 
\bar\pi(
\bar G) ,$
\item  $\< \pi(b), \pi(d)\>  = \Sym(\{l + 1, \dots  , n\})    <  \pi( G ) ,$
and  
\item  $\<  \bar\pi(\bar a),  \bar\pi(\bar e)\>  = \Sym(\{-m + 1 + k, \dots ,
l\} )
  <  \bar\pi( \bar G) .$
\smallskip
\end{itemize}  
Then 
\begin{equation}
\label{general symmetric}
J= \< X, \bar X \mid R, \bar R, a = \bar a, c = \bar c, [d, \bar e] =1 \> 
\end{equation}
is a presentation of $S_{m+n -k} =\Sym \{ -m+k+1, \dots ,n \}$$,$
where  
$\<X\>  =S_n $  acts on $\{1,\dots,n\}$ and   $\<\bar X\>  =S_m $
acts on
$\{-m + 1 + k, \dots , k\}$.

\end{lemma}  

The following picture might be helpful.
\smallskip
$$
\begin{array}{cccc}
\multicolumn{4}{c}
{\underline{\quad
-m+k+1, \dots, 0\quad\quad\quad  1, \dots, l\quad \quad\quad \quad\quad
l+1,\dots, k
\quad\quad\quad  k+1,\dots ,n \quad}}\\
&\  \quad\quad\quad\quad\quad\quad\quad\quad \underline{\quad\quad a=\bar a
\quad\quad} & 
\raisebox{2.5ex} {~}  
\quad\quad\underline{\quad\quad b=\bar b \quad\quad} \quad\quad\quad& \\
& \multicolumn{2}{c}
{\quad\quad\quad\quad\quad\quad\quad\quad\underline
         {\quad\quad\quad\quad\quad\quad\quad
c =\bar c \quad\quad\quad\quad\quad\quad\quad}} \quad\quad \ & \\
\multicolumn{2}{c}{\underline{\quad\quad\quad\quad\quad\quad\quad\quad 
\bar e \quad\quad\quad\quad\quad\quad\quad\quad  }} &
\multicolumn{2}{c}{\underline{\quad\quad\quad\quad\quad\quad\quad\quad
 d \quad\quad\quad\quad\quad\quad\quad\quad}}
\end{array}
$$
 
\medskip

\Proof  
The restrictions on $m,n,k$ and  $l$ are designed to guarantee that the
desired permutations exist.   
There is a surjection   
$J\to S_{m+n -k}$ (note that our
  3 extra relations are satisfied).
By  Lemma~\ref{It's a subgroup}, $J$ has subgroups we identify
with 
$G=\<X\> =S_n$ and   $\bar G=\<\bar X\> =S_m$.

By our relations,   $ a =\bar a  $ and $ c =\bar c  $.   
Then the assumption $ \pi(b )= \bar \pi(\bar b )$ states that $b $ 
and $\bar b$  represent 
the same element of  $
\<  a,  c  \> =\<\bar a,\bar c\>$, so that
  the additional relation
$b =
\bar b$ is forced to hold  in $J$.
  Then  we also  have the following relations: \vspace{2pt}
\begin{itemize}  
\item $   [d, \bar a] = [d, a] = 1$
because $d,a \in G=S_n$ have disjoint supports, 
\item$  [d, \bar e] = 1$  by   the last relation in the presentation
\eqn{general symmetric},
 
\item $ [b, \bar a] = [b, a] = 1$ because $b,a\in G$ have disjoint
supports,  and
\item $ [b, \bar e] = [\bar b, \bar e]\, = 1$ because $\bar b,\bar e 
\in\bar  G$ have disjoint supports. 
 \vspace{2pt}
\end{itemize} 
Therefore    
\begin{equation}
 \label {commuting symmetric}
[\< b, d \>, \< \bar a, \bar e  \>] =1,
\end{equation}
where $\< b, d \> =\Sym(\{l + 1, \dots  , n\} ) $ and  $ \< \bar a, \bar e 
\>=
\Sym(\{-m + 1 + k, \dots , l\} ).$

The symmetric groups  $G$  and  
$\bar  G$   are generated, 
respectively,   by 
\vspace{1pt}
the $n-1$ and $m -1$  
 transpositions    
$x_i: =(i,i+1)$, $1\le i<n$,   and
$x_i: =\overline{(i,i+1)}$, $-m+1+k\le i<k$.   The identification of the
two copies of $S_k  =  \<  a,  c  \> = \< \bar a, \bar c \>$
in \eqn{general symmetric}  
identifies  the transpositions $x_i$, $1\le i<k$, common to  these
generating sets,  producing  a generating set of  $J$ consisting of $m+ n
- k- 1$ involutions. These involutions satisfy the relations in the
Coxeter  presentation \cite{Moo}
$$
\begin{array}{lll}
 S_{m+n -k}=\langle x_i,  -m+1+k\le i < n
& \hspace{-8pt} |& \hspace{-12pt}
 x_i^2= (x_ i x_{i+1})^3=
(x_i x_j)^2=1  \vspace{2pt}
\\
  &&\hspace{-12pt}   \mbox{for all possible $i,j$ with  
$ j-i\ge 2\rangle$$: $}   
\end{array}   \quad \quad
$$
   any two $x_i$  either both lie in
$G$, or   both lie in  $\bar G$, or  they  commute  by \eqn {commuting
symmetric}  since one is in $\<b, d \>$ and the other is  in
$\<\bar a, \bar e  \>$.  \qed 
 
\medskip
There is a great deal of flexibility in the choice of  the elements $a, b, c,
d, \bar a, \bar b, \bar c, \bar e  $.
In the proof of  Theorem~\ref{Alt(n) and Sym(n) 100 down to 7}   we will need
to require that $|a|=3$, which is  possible provided that $l\ge3$.

The last three relations in \eqn{general symmetric}  are similar to ones
used in  the proof of
\cite[Theorem ~3.17]{GKKL}; and they are used both there and here in 
essentially the same manner.  However, the situation in that paper was more
delicate, due to length considerations: the preceding  presentation does not even
have bit-length
$O(\log n)$.  We will deal with this   requirement when we use this lemma
in the proof of Theorems~\ref{Alt(n) and Sym(n) 100} and \ref{Alt(n) and
Sym(n) 100 down to 7}.

\begin{lemma} 
\label{Alt m+n-k}
A presentation of $A_{m+n -k} $  is obtained
  as  in the preceding lemma by  
replacing symmetric groups by alternating groups throughout   
{\rm\eqn{general symmetric}} and assuming that
$m,n>k\ge l+3\ge 6$.

\end{lemma} 

\Proof 
Once again, the restrictions on $m,n,k$ and $l$ are designed to guarantee
that the desired permutations exist.   
The previous picture can again be used.  
 As in the proof of the
preceding lemma, \eqn{general symmetric} implies that  \eqn{commuting
symmetric} holds.

 We will use the presentation (\ref{Carmichael An}) with the union of
\vspace{1pt}
the the   generating sets 
 $x_i:=(1, 2,i ) $ for 
$G$ and $x_j:=\overline {(1, 2,  j) }$ for $\bar G$,
where 
\vspace{1pt}
   $3\le i\le n $  and $-m+1+k\le j\le k$ but $j\ne 1,2$.  As above,
\eqn{general symmetric} implies that
$(1, 2,i ) =\overline {(1, 2,  i) }$ if $3\le i\le k$.

Then all required relations in  (\ref{Carmichael An})   are obvious except
for 
$$
\mbox{$ \big((1, 2, i )\overline {(1, 2,   j)} \big)^2 = 1$  with  $ k
< i\le n$ and $-m+1+k\le j \le 0.$}
$$
Recall that $\< b, d\>  = \Alt(\{l + 1, \dots  , n\}) <G 
 =\Alt(\{1,\dots , n\} )$,
where $6\le l + 3\le k<i\le  n$.
Then some  $g \in  \<b,d\>$   sends $k$ to $i$ and fixes 1 and 2.
By \eqn{commuting symmetric}, $g$ commutes with 
$\overline{(1, 2, j  )}  \in \<  \bar a, \bar e\>$.  Consequently,
$$\big[\big((1, 2, i) \overline{(1, 2, j  )}\big)^2\big]^g =
\big((1, 2, k) \overline{(1, 2, j )}\big)^2 =  \big(\overline{(1, 2,
k)}\,\overline{(1, 2,  j )}\big)^2 = 1 ,$$   as required  in
(\ref{Carmichael An}). 
\qed 

\begin{corollary}
\label{An 100 with y}
  If $A_m$ has a presentation  with $M$
relations and if $m >  k \ge  6,$ then $ A_{2m-k}$  has a presentation
with
$M + 4$ relations. The same holds for the corresponding symmetric groups
using the weaker assumption $m >  k \ge  4$.
\end{corollary} 

\Proof  
 Let $G=\< X  \mid R\>$ be a presentation for $
A_m$  with
$M$ relations.  In  Lemma~\ref{Alt m+n-k} we  use $m=n$  and   $l=3$, but this
time we   introduce an additional generator $y$ corresponding to an even 
permutation sending  $\{1,\dots, m\} \to
\{-m + 1 + k, \dots , k\}$ and inducing the identity on $\{1,\dots, k\} $.

Consider the group
\begin{equation}
\label{8- relations}
J:= \< X ,y \mid R,  a =  a ^y, c =  c ^y, [d,  e ^y] =1\>,
\end{equation}
with  $a,b,c,d  ,\bar a:=a^y,\bar b:=b^y ,\bar c:=c^y, \bar
e:=e^y $ playing the same roles as in Lemma~\
\ref{Alt m+n-k}. 
By that lemma  with $\bar X:=X^y$ and $\bar R:=R^y$, $J$ has a subgroup
$K:=\<X,X^y\>\cong  A_{2 m-k}$. 

Finally, we add an extra relation to ensure that our generator
$y$  is in $K$  and that, as an element of $ A_{2 m-k}$, the action of  $y $ 
 on  $\{ -m+k+1, \dots ,m \}$ is as described above.

The group $S_{2m-k}$ is dealt with in the same manner. 
\qed

\Remark\rm
\label{gluing Sns}
We have just glued two subgroups $A_m$ in order to obtain a group
 $ A_{2m-k}$,  or two subgroups $S_m$ in order to obtain a group
 $ S_{2m-k}$, in each case with suitable restrictions on $m$ and $k$. 
There is a variation on this process that glues two subgroups $S_m$ in
order to obtain a group
 $ A_{2m-k}$ (view $S_m$ as lying in $ A_{m+2}$, as was done on occasion
in Section~\ref{Using 2-transitive groups}).

\subsection{All alternating and symmetric groups}

\label{All alternating and symmetric groups} 
Before proving  Theorem~\ref{C}, we begin  with a weaker result:

\begin{prop} 
\label{Alt(n) and Sym(n) 100}  
For all $n\ge 5,$ $A_n$ and $S_n$  have  
presentations with
$3$ generators  and  $10$  relations.
\end{prop} 

\proof
If $n\le 10 $ then 
$n=p+2 $  or $p+3$ for a prime $p$, and we have already
obtained a presentation with fewer relations than required.
By Ramanujan's version of Bertrand's Postulate \cite{Ramanujan},   
in all other cases  we can write  
$n = 2p + 4 - k $  for  a  prime  $p$ and an integer $k $ such that  
$m:=p+2 > k\ge 6$, and then  use  Corollaries~\ref{Ap+2} and  \ref{Alt
m+n-k}. (A related use of Bertrand's Postulate appears in \cite[Theorem
3.9]{GKKL}.)  \qed

\medskip
This  proposition   is weaker than Theorem~\ref{C} in two significant ways:  the number of relations is larger than in that  theorem, and  bit-length is not mentioned.  We deal with the second of these as follows:
\begin{Lemma} 
\label{Alt(n) and Sym(n) bit-length} 
In {\rm Corollary~\ref{An 100 with y}}  and {\rm Proposition~\ref{Alt(n) and Sym(n) 100},} it is possible to choose  $a,c,d,e$ such that   $a$ is   a $3$-cycle and
the resulting  presentation has  bit-length  $O(\log n)$ with  a bounded number of exponents$,$
each of which is at most $n,$  if either
\begin{itemize}
\item[\rm (i)]
the group is $A_n,$ or

\item[\rm (ii)]
the group is $S_n$ and we used a prime $p\equiv 11$ \rm  (mod $12$). 

\end{itemize}
\end{Lemma} 

\Proof
In Lemma~\ref{Sym m+n-k} and    Corollary~\ref{An 100 with y}   we can choose each of the
elements  $a,c,d,e$ to be a product of a cycle of the form $(i,\dots,j)$ with
$i<j$ and a permutation having bounded support (in fact, in 
Section~\ref{An explicit presentation for $S_n$} each  will   be chosen to be a cycle).  We require
$a$ to be a 3-cycle   (this is needed in Theorem~\ref{Alt(n) and Sym(n) 100 down to 7}); and then
 in the symmetric group case we will need $c$ and $e$ to be odd permutations
 (cf. the hypotheses of Lemma~\ref{Sym m+n-k}).

By Remark~\ref{cycles for p+2},  when we glue two copies of $A_{p+2}$
using  relations   of bit-length $O(\log n)$, the   bit-length and  exponents are as
required, except perhaps  for the crucial  additional relation  
expressing  $y$ as word in  $ X \cup X^y$. 

 We now consider
$y$; there is a reasonable amount of flexibility in the choice of $y$ in
Corollary~\ref{An 100 with y} (recall that  $m=p+2$). In that corollary we
were permuting the
$n=2p+4-k$ points 
\begin{equation}
\label{list of n points}
-p - 1+k, -p +k,\dots -1,0 ;~ 1,\dots ,k ; ~k+1,k+2, \dots ,p+1 ,p+2,
\end{equation}
 where we have
alternating or symmetric groups on the first and last $p+2$ points, with
  an $A_k$  or  $S_k$  on the overlap.

If $p-k$ is even  then we choose  $y$ to be the following  product of $p+2-k$
transpositions:
\begin{equation}
\label{def of y}
y: =(-p - 1+k, p+2)(-p +k,p+1)  \cdots (-1,k+2)(0,k+1).
\end{equation}
We use the following additional permutations:
 \vspace{2pt}
 
 \begin{itemize} 
\item  $x : = (-1,k+2)(0,k+1) = [(1,k+2)(2,k+1)]^{(-1,1)(0,2)},$ which we have
written using permutations from the two alternating groups, and
 \vspace{2pt}
\item   $u^{-1}: =(1,\dots,k,k+1,\dots,p+2)  (1,\dots,k,0,-1,\dots,-p+k-1)   
$

\hspace {19pt} $=(1,\dots,k,k+1,\dots,p+2)  (1,\dots,k,k+1,\dots,p+2)  ^y  $.
 \vspace{2pt}
\end{itemize}  
 By  Remarks~\ref{cycles for p+2} and \ref{cycles for p+2 again}, $u$ and
hence also $s$ can be expressed as a word of bit-length $O(\log p)$ in
$X\cup X^y$  using a bounded number of exponents;
then so can 
\begin{equation}
\label{y }
y =xx^{u^2}x^{u^4}\cdots x^{u^{p -k }}  = (xu^{-2})^{(p -k)/2} x u^{p -k} ,
\end{equation}
using (\ref{Horner}).

If $p-k$ is odd   let $v := (-p - 1+k, p+2)(-p+ k,p+1)  \cdots 
(-3,k+4)$ and use 
\begin{equation}
\label{def of y odd case}
\begin{array}{lll}
y:= &   \hspace{-6pt}   v (-2,k+3)  (-1,k+2,0,k+1) \vspace{2pt} \\
& \hspace{-18pt} = \hspace{1pt} v  [( 2,k+3)  (
1,k+2,k,k+1)]^{(2,-2)(1,-1)(0,k)}. 
\end{array}
\end{equation} 
Then $v$   can be expressed as a
word of bit-length $O(\log p)$ in
$X\cup X^y$   using a calculation similar to  \eqn{y }, and the
final term in $y$ is a product of permutations from the two copies of $A_{p+2}$.
Another application of  
  Remarks~\ref{cycles for p+2} and \ref{cycles for p+2 again} completes the
proof of (i).

Now  (ii) follows from  Remark ~\ref{cycles for p+2 again}.  \qed

\smallskip 
\Remark \rm
\label{small bit-length}
Using Remark~\ref{cycles for p+2}  we see that  
{\em   every cycle $(k,k+1,\dots,l) $  and every element with bounded
support  in $A_n$  has bit-length $O(\log n)$ in our generators. } 

\medskip  
With a bit more number theory, together with Table~\ref{Small degrees},
we obtain an improvement of  Proposition~\ref{Alt(n) and Sym(n) 100}  that is needed
for Theorem~\ref{C}:

\begin{prop} 
\label{Alt(n) and Sym(n) 100 down to 8}  
If  $n\ge 5$ then 
$S_n$  and 
$A_n$   have presentations  with 
$3$ generators$,$  $8$ relations and bit-length $O(\log n)$.
Moreover$,$  these presentations  use  a bounded number of exponents$,$
each of which is at most $n$.
\end{prop}  
 
\proof We refine the argument in  Proposition~\ref{Alt(n) and Sym(n) 100}. First consider 
$ S_n$.   Here  we need  to write 
$n = 2p + 4 - k$  for  a  prime  $p\equiv 2$  (mod 3) such that   $m:=p+2 >
k\ge 4$, so that we can   use  Corollaries~ \ref{Sp+2 4
relations}  and \ref{An 100 with y}, and then continue  
 as in the proof of  Proposition~\ref{Alt(n) and Sym(n) 100}. 
In view of the requirements on bit-length and exponents, we also   {\em  require  that }  $p\equiv 11$  (mod 12) {\em if $n\ge 50$},   so that Lemma~\ref{Alt(n) and Sym(n) bit-length} will  complete
the proof  for $ S_n$.

According to Dirichlet's Theorem, for large $x$
there are approximately $ x/2\log x$ primes $\le x$ of the stated sort, and
subtraction yields a prime $p$ in our situation.  However, we need a more
precise (and effective) result of this type. This is provided in  
\cite{Mor} (updating \cite{Bre, Er, Mol}  with
more  precise estimates):    
if $n\ge 50 $ then there is such a prime $p\equiv 11$  (mod 12).  A straightforward examination
of the cases  $n <50$  leaves
$n=11,12$   or  $13$ to be dealt with.  
See \cite[p.~54]{Ar} for  presentations of $S_{11}$, $S_{12}$ and  $S_{13}$
with 2 generators and 6, 7 and 7 relations, respectively
(cf. \cite[p.~64]{CoMo}).  (Note that Table~\ref{Small degrees} contains the cases $n=11,12$, while Corollary~\ref{Sp+2 4 relations}(i)   handles the case $n=13$.)
 
For $ A_n$ we need  to write 
$n = 2p + 4 - k$  for  a  prime  $p\equiv 11$~(mod~12) such that   $m:=p+2 >
k\ge 6$,  then  use  Corollaries~ \ref{Ap+2 4
relations}  and \ref{An 100 with y}, and again  finish   as in the proof of 
Proposition~\ref{Alt(n) and Sym(n) 100},
using Lemma~\ref{Alt(n) and Sym(n) bit-length}.
Once again, by  
\cite{Mol, Mor},    if $n\ge 50$ then there is such a prime $p$.    Another
straightforward examination leaves  the cases  $n\le13$ and $n=21, 22, 23, 24,
25, 45, 46, 47, 48$  or $49$  to be dealt with.  The cases in which $n=p+2$ or
$p+3$ for some prime $p$  are handled using examples  described earlier, and
Table~\ref{Small degrees} handles the remaining cases.  (For presentations of
$A_{11}$,  $A_{12}$ and  $A_{13}$ using 2 generators and 6,  7 and 7
relations, respectively, see \cite[p.~67]{CoMo}.)~\qed 

\medskip

For the next theorem we will use a presentation in the preceding proposition  
that is valid for most
$n$.   
If  $n\ge50$  (or, more precisely, if $n$ is not one of the exceptions
mentioned in the above proof) then 
\eqn{8- relations}  together with one further relation  is such a
presentation:
\begin{equation}
\label{presentation for use}
\mbox{$A_n$  or $S_n=
\< X ,y \mid R,  a =  a ^y, c =  c ^y, [d,  e ^y] =1, y=w\>,$}
\end{equation}
with $\< X \mid R \>$
in (\ref{Ap+2}) or (\ref{Sp+2})  for $A_n$  or $S_n$,  respectively,
 $a,c,d ,e  $ as in Lemma~ \ref{Alt(n) and Sym(n) bit-length}},
and  a suitable    word $w$  in $X\cup X^y$ 
as in   \eqn{def of y}--\eqn{def of y odd case}.  (The  properties required of
$y$ and $w$ are described at the end of the proof of  Corollary~\ref{An 100
with y} and, in gory detail, in the proof of Lemma~\ref{Alt(n) and Sym(n) bit-length}.) 

We are now able to prove Theorem~\ref{C}: 
\begin{theorem} 
\label{Alt(n) and Sym(n) 100 down to 7}

If  $n\ge 5$ then 
$S_n$  and 
$A_n$   have presentations  with
$3$ generators$,$  $7$ relations and bit-length $O(\log n)$.
Moreover$,$  these presentations  use  a bounded number of exponents$,$
each of which is at most $n$.
\end{theorem}

\Proof  
Let $n=2m-k$ with $m=p+2 >k\ge 6$ (cf. the preceding proposition; below we will
discuss the existence of a suitable prime $p$).

Let $G=\<X\mid R\>$ be the    presentation in 
Corollary~\ref{Ap+2 4 relations}(i) or \ref{Sp+2 4 relations}(ii), so  that
$|X|=2$,  $|R|=3$ and one of the following holds: 
$$
\begin{array}{llll}
\mbox{$A_{n} $  case: }&T=\<{\bf a},{\bf b}\>=\AGL(1,p)^{(2)}~
&  \mbox{$p\equiv11$  (mod 12)}~ & G\cong A_m\times T  \vspace{2pt}\\
\mbox{$S_{n} $ case: }&T=\<{\bf a},{\bf b}\>=\AGL(1,p) 
&  \mbox{$p\equiv 2$  (mod 3)} & G/(1\times T)\cong S_m. \\
\end{array}
$$ 
(In the $A_n$ case $T$ has index 2 in $\AGL(1,p)$; 
in the $S_n$ case $G$ has index  2 in $ S_m\times T$.)
We also require that  $p\equiv 11$  (mod 12) in the   $S_n$ case when  $ n\ge 50$.

Let $t\in T  $ be such that  $T\cap A_{m}=\<(t^3)^T\>$ (for example,  $t={\bf
b}^2$ works since {$p\equiv 2 $~(mod~3)   and the order of $ \bf b $ is odd
in the $A_n$  case)}.%

Then  the presentation 
(\ref{presentation for use}) of $A_{n} $ or $S_{n} $  can be rewritten
\begin{equation}
\label{presentation with t cubed}
 \<X,y \mid R, t^3, a^y = a , c^y=c , [d^y,e ]=1, y=w\>.
\end{equation}
There was a great deal of freedom in our choice of the elements $a,b,c,d,e$
in the proofs of Lemmas~\ref{Sym m+n-k}  and \ref{Alt m+n-k} (and
Corollary~\ref{An 100 with y}).  
   As in Lemma~\ref{Alt(n) and Sym(n) bit-length}, we now  {\em choose}
$a$ so that its image  in the alternating or  symmetric group  $G/(1\times
T)$  is  a 3-cycle  (for the   $S_n$ case, this requires   
$c$  to be    an odd permutation   and  $l\ge3$ in the proof of Lemma
3.15).

 The presentation \eqn{presentation with t cubed}  has 8 relations and 
bit-length $O(\log n)$.  

We use the following additional ingredients:
\smallskip
\begin{itemize}
\item  
Write  $a\in G $  as $a=(a_1,*)\in A_m\times T$
with  $a_1 $ of order 3.
\item Let  $\bar{a }$  and $\hat{a }$ be   words in $X$ such that 
$ \bar{a } = (a_1, 1 ) $ and 
$\hat{a } = (a_1, t ) $  when evaluated in $G$.
\end{itemize} 
\smallskip
 
{\em  We claim that the $7$-relator group
\begin{equation}
\label{8 to 7}
J:= \<X,y \mid R,  \bar{a } ^y = \hat{a }, c^y=c , [d^y,e ]=1,
y= w\>
\end{equation}  
is isomorphic to the group in} (\ref{presentation with t cubed}). 

For, we can view $G=\<X\>\le J$.   Then $a,\bar{a },\hat{a }, c,
d, e \in J$.

Since both $\bar{a } ^y=\hat{a}= (  a_1, t) $ and $ a_1$   have  order
$3$ we have $(1,t)^3=1$ in $J$.  Thus,  $J$ satisfies the
presentation   (\ref{presentation with t cubed}) and hence is as claimed.

{\em Bit-length}: As in the proof of Theorem~\ref{Alt(n) and Sym(n) 100}, 
$a,c,d$ and $e$ can be expressed as words in $X$ of  bit-length $O(\log n)$
mod~$T$.  Since  $T=\<\bf a\>\<\bf b\>$ in Corollaries~\ref{Ap+2 4 relations}
and \ref{Sp+2 4 relations}, 
$\bar{a }$  and $\hat{a}$ have  bit-length $O(\log
n)$  as well.  

Finally, we need to discuss whether we have  handled all groups $A_n$ and
$S_n$; or, what amounts to the same thing, for which $n$  a prime $p$ can
be found satisfying  all of the conditions we have imposed.  

As in Proposition ~\ref{Alt(n) and Sym(n) 100 down to 8},  (\ref{8 to 7})
takes care of $A_n$ except for the cases $n\le13$ and $n=21, 22, 23, 24, 25, 45,
46, 47, 48, 49$; and these are handled exactly
as in that proposition.

For the $S_n$  case   we have  imposed a further condition beyond   what was used in 
Proposition ~\ref{Alt(n) and Sym(n) 100 down to 8}:  
we need  to write 
$n = 2p + 4 - k$  for  a  prime  $p\equiv 2$  (mod 3) such that     
 $m=p+2 > k\ge l+2\ge5$,  and  
  $p\equiv 11$  (mod 12)   if $n\ge 50$.
 (The conditions in Corollary~\ref{An 100 with y}  were $m =p+2 >
k\ge l+2\ge4$, but here we need to be able to find a 3-cycle $a$ in $A_l$.)
 Once again these requirements can be met for all $n$ except if $n<6$ or $n=9,
10, 11$, and  those cases can be handled as before.
\qed
\medskip 

None of  the  presentations in this or the preceding section   has length
$O(\log n)$.  

 By   Remarks~\ref{cycles for p+2} and \ref{cycles for p+2 again},  the
exponents in \eqn{y }  are all less than $n$; there is also an  exponent 
$p+1-k$ used to write
$u$ in that remark.  
As already noted, these  presentations have {\em bounded expo-length} (cf.
Section~\ref{Preliminaries}).   See Remark~\ref{Concluding remarks
expo-length} in Section~\ref{Concluding remarks} for comments concerning the
boundedness of expo-length for  other families of almost simple  groups.
 
Before continuing, we note the following simple improvement of
\cite[Lemma 2.1]{GKKL}. 
 
\begin{lemma}
\label{d generators}
Let $G=  \langle D  \rangle$ be a finite group  having  a presentation
$\langle X \mid  R \rangle  ;$ let $\pi \colon F_X\to G$ be the
natural map from the free group $F_X$ on $X$.  Then
$G$ also has   a presentation
$ \langle D \mid R'\rangle $  such that  $|R'|= |D|+|R|-|\pi(X)\cap D|.$
\end{lemma}

\Proof We recall the simple idea used in the proof of
\cite[Lemma~2.1]{GKKL}.   
Write each $x\in X$ as a word $v_x(D)$ in $D,$ and each $d\in D$ as a
word
$w_d(X)$ in $X$; and let  $V(D)=\{v_x(D) \mid x\in X\}$. 
According to the proof of \cite[Lemma~2.1]{GKKL}, we then obtain another
presentation  for $G$: 
$$G=\langle D  
\mid   d =w_d(V(D )),
\, r(V(D )) =  1, d\in D, r\in R \rangle.
$$
For each $d\in \pi(X)\cap D$,  one of the above relations can be taken to
be 
$d=d$, and hence can be deleted.~\qed
 \medskip

In view of the desire in \cite{Con} for specific generators
(namely, $(1,2)$ and $(1,2,\dots,n)$),
we note the following consequence of the preceding theorem and lemma:

\begin{corollary}
\label{any generators of Sn}  
Let $G=A_n$ or $S_n ,$ $n\ge 5.$ 
\begin{itemize} 
\item[\rm(i)]
If $a$ and $b$ are  any  generators of $G,$ then  
there is a presentation of $G$ using $2$ generators that map onto 
$a$ and $b,$ and    $9$ relations. 

\item[\rm(ii)]
There is a presentation of $G$ using $2$ generators and  $8$
relations.
\end{itemize}
\end{corollary}

We do not have information 
 concerning the bit-length of any of the
resulting presentations.

\proof
Part (i) follows from 
Theorem~\ref{Alt(n) and Sym(n) 100 down to 7} and
the preceding lemma, using $|\pi(X)\cap D|\ge 0$. 

For (ii), note that
we have provided a presentation $\langle X \mid  R \rangle  $ for  $G$
such that~some element of $X$    projects onto  an element   $a\in G$ that
is  either a 3-cycle ($z$ in Lemma~\ref{using Carmichael}) or has a power
that is a 3-cycle (such as
$g$ in Corollary~\ref{Ap+2 4 relations}(ii) or  Proposition~\ref{general
number of relations}).  Let $b$ be  any element of
$G$ such that 
$G=\langle a,b \rangle  $.  Now  use $D=\{ a,b\}$ in the preceding
corollary (compare Section~\ref{Concluding remarks}, Remark~\ref{one less
generator}).
 \qed
 
\subsection {An explicit presentation for $S_n$} 
\label {An explicit presentation for $S_n$}
The presentations in Sections~\ref {Using 2-transitive groups}  and    \ref{Small n}  are not difficult to understand, and  they visibly encode information concerning various  alternating and symmetric groups.
However,  the presentations in
Theorem~\ref{Alt(n) and Sym(n) 100 down to 7} are  not as explicit as one
might wish.  
Therefore, we will provide  a presentation
of $S_n$ for $n$ odd (see Remark~\ref{even n}  for even $n$). Although this presentation is in no sense elegant or informative, it does have the significant advantage of using only 7 relations.

Find a prime  
\mbox{\rm $p\equiv 11$ (mod 12)}  such that  
$n-1\ge p\ge (n+2)/2$.    (This places a mild restriction on $n$, as seen in the proof of Theorem~\ref{Alt(n) and Sym(n) 100 down to 7} .   For $n\ge 50$ there is always such a prime.)

Let $k=  2p+4-n$, so that $p+2>k\ge 6$.
Then $k\equiv n  \equiv 1$ (mod 2),  so that  $p-k$ is even.

{\em The  desired presentation is} 
$$\begin{array}{lll}
 S_n= \< {\bf a} , {\bf g} , {\bf y}    \mid  &\hspace{-8pt} 
{\bf a} ^p = ({\bf g} ^{3})^{p-1 }, 
({\bf a} ^s)^{{\bf g} ^{3}} 
\!=\!{\bf a} ^{s-1} ,   \big ({\bf g} ^{ p-1 } ({\bf g} ^{ p-1})^{\bf a} \big)^2  
\hspace{-2pt}  
 \hspace{-2pt}=1 ,
\vspace{2pt} 
\\
\vspace{2pt}
  &\hspace{-8pt}   {a } ^{\bf y}  = \hat{a }, c^{\bf y} =c , [d^{\bf y} ,e ]=1,
{\bf y} = w
\>  ,
\end{array}   
$$ 
for words $a,c,d,e,\hat a,  w$   defined below and integers 
 $r$ and $s$   such that   $s(r-1)\equiv -1$ (mod $p$)   and $ \F_p^* =\<r\>$.

   Notes:  
In order to conform with the notation in Section~\ref {Using 2-transitive groups}, we view $\AGL(1,p)$ as acting on 
$\{1, \dots,p\}$  with
   ${\bf a} \equiv  (1,\dots,p)\in \AGL(1,p)$,  but we    
use $p$ in place of 1 in Lemma~\ref{using Carmichael}.  
   Finally,   we  use   ${\bf a} $ in place of $a$  since the latter plays a prominent role  in Section~\ref{Gluing alternating and symmetric groups}.  
   
   The permutations in $S_n$ indicated below are not part of the presentation, but are provided in order to help keep  track of the map $\< {\bf a} , {\bf g} , {\bf y}   \>\to S_n$ into the symmetric group on the  $n$ points 
   (\ref{list of n points}).
The notation used  here should  {\bf not} be viewed mod~$p$.%
\vspace{2pt}
    
\begin{enumerate}
\item
$ {\bf z} \hspace{1pt} :={\bf g} ^{p-1} \equiv  (p,p+1,p+2)$,
\vspace{2pt}

\noindent
${\bf  b} : = {\bf g} ^3 $  (so that $\< {\bf  a} ,{\bf  b}\> =\AGL(1,p)$),
\vspace{2pt}

\item
$ {\bf z}(i) := {\bf z} ^{{\bf a}^{-i}} \equiv (i,p,p+1)$  for $1\le i\le  p$,
\vspace{2pt}
 
$c({i, j}):=    {\bf z}(-i)  {\bf z}(-j)^{-1}   {\bf z}(-i)  \equiv (i,j)(p+1,p+2)$
for $1\le i <  j\le p$,
\vspace{2pt}
 
\noindent
$c_i:= \big (c({1,2}){\bf a}  \big ){}^{i-2}  c(1,2) {\bf a} ^{-(i- 2 )} \equiv (1,2,\dots,i) $
 with $i$  odd  and   $3\le i  \le p  ~$   (cf. (\ref{Horner})).

\item  (Constructing a transposition)

\noindent
$b_2:={\bf b} ^{(p-1)/2 }\equiv  (1,p-1)(2,p-2)\cdots \big ((p-1)/2,p-(p-1)/2 \big)$,
\vspace{2pt}


\noindent
$c_{\bullet} : = c_{(p-1)/2} ^{~^{~}}    c_{(p-1)/2}^{-{\bf a}^{(p+1)/2}}  \equiv  \big(1,2,\dots,(p-1)/2 \big) \big(p-1,p-2,\dots, p-(p-1)/2 \big)$,%
\vspace{2pt}
  
$v : =\big (c(1,p-1) c_{\bullet}  ^{-1} \big ){}^{(p-1)/2}c(1,p-1)c_{\bullet}^{(p-1)/2}   $%
\vspace{2pt} 

\hspace{8pt} $ \equiv 
(1,p-1)(2,p-2)\cdots \big ((p-1)/2,p-(p-1)/2\big ) (p+1,p+2)~$   (cf.
(\ref{Horner})),
\vspace{2pt} 

$t  \hspace{1.5pt} : = v b_2 c(1,2) \equiv (1,2)$.
\vspace{2pt}

\item
$a:= {\bf z}(3) ^  {   {\bf z}(1) {\bf z}(2)     }   \equiv (1,2,3)$,
\vspace{2pt}
 
\noindent
$  c \hspace{1pt}  :=  t c_k    \equiv(2,\dots ,k) ~$  (an odd permutation),
\vspace{2pt}

\noindent
$d:=c_3^{-1}{\bf az}  \equiv (3,\dots ,p+2)$,
\vspace{2pt}
 
\noindent
$ e \hspace{.5pt}   :=   {\bf a} c_k^{{-1}  }  t {\bf z}   \equiv (1,2,k+1,\dots,p+2)$,
so that $ e^{\bf y}\equiv  (1,2,0,-1,\dots,-p  + k-1) ~$    (also  odd permutations).
\vspace{2pt} 

\item 
$\hat a :  = a \big ({\bf b} ^2z(1)z(-1 ) \big )^{(p+1)/2}   \in \{ a\} \times T  $.
\vspace{2pt} 
 
\item
$x : =  [c( 1,k+2 )c(2,k+1 )] ^ {  c(1,k+2 )^{\bf y}c(2,k+1 ) ^{\bf y}}  $

\hspace{8pt} $  
 \equiv   (-1,k+2)(0,k+1) = [(1,k+2)(2,k+1)]^{(-1,1)(0,2)}$,
    
\vspace{2pt}
\noindent
 $u^{-1}: = ({\bf a} {\bf z}) ({\bf a} {\bf z})^{\bf y}$
 
\hspace{19pt} $ \equiv (1,\dots,k ,k+1,k+2,\dots,p,p+1,p+2)( 1,\dots,k ,0,-1, \dots, 
- p+k-1   )$,%
\vspace{2pt}

\noindent
$w:=(x  u^{2})^{(p -k)/2} x 
u^{-(p -k)}  $

\hspace{9.3pt} $  
\equiv (-p - 1+k, p+2)(-p +k,p+1)  \cdots (-1,k+2)(0,k+1)~$
(cf.   (\ref{Horner}). 
\end{enumerate}

\Remarks\rm
\label{even n}
1. In the isomorphism given by  (\ref{def phi}), $\bf z$ maps to an element of   $A_{p+2}\times \{1\}$.  Hence, the element $a$  defined above also maps into  $A_{p+2}\times \{1\}$,
 so that the element $\bar a$ used in (\ref{8 to 7}) is just our $a$.  The remainder of the presentation given above is just a straightforward translation from Section~\ref{All alternating and symmetric groups}.
 
  \vspace{2pt}
  
  2. 
 A presentation for  the alternating groups   is similar but slightly simpler:  only  even  permutations are involved. 
 \vspace{2pt}

3.   
 \emph{Changes needed  when  $n$  and $k$  are even}:
\vspace{2pt}

\begin{itemize}
\item [(4$'$)]
 $c:=c_{k-1}  t \equiv  (1,\dots,k) ~ $  (an odd permutation),
\vspace{2pt}
 
 $e:= {\bf a} t  c_{k-1} ^{-1} { \bf   z}  \equiv  (1,k+1,\dots,p+2)~$    (another  odd permutation).
\vspace{2pt}
  
\item [(6$'$)] 
$w:=   t^{\bf a z}    t^{\bf a z}   {}^y   t^{\bf a z} (x  u^{2})^{(p -k-1)/2} x 
u^{-(p -k-1)}  $
\vspace{2pt}

\hspace{10pt} $  
\equiv (-p - 1+k, p+2)(-p +k,p+1)  \cdots (-1,k+2)(0,k+1)~$  as before.
\end{itemize}

\subsection {Weyl groups} 
\label {Weyl groups}
It is easy to use Theorem~\ref{Alt(n) and Sym(n) 100 down to 7}   to obtain 
presentations for the Weyl groups of types $B_n$ or $D_n$.  However, we leave
this to the reader, dealing instead with a subgroup $W_n$ of those Weyl groups
 that is  needed later (in Section~\ref {More presentations of
classical groups}).

Let $W_n :=
\Z_2^{n-1}
\semi A_n$ be  the subgroup of the monomial group of $\R^n$ such that 
$ \Z_2^{n-1}$ consists of all $\pm1$ diagonal matrices of determinant 1,
and the alternating group  $A_n$ permutes the standard basis vectors
 of $\R^n$.  

\begin{prop}
\label{Wn}
If $n\ge 4$ then   $W_n$ has  a presentation with  $4$ generators$,$
$11$ relations and bit-length $O(\log n)$. 
If $n=4$  or $5$ then   $W_n$ also  has  a presentation with  $3$ generators
and
$7$ relations. 
\end{prop}

\proof
We first consider the case $n\ge 5$.   
Let $\langle  X \mid R\rangle$  be a presentation for  $A=A_n$.  Let  
$\sigma  = (1,2,3) \in A$, choose a 2-element generating set of  $ H: =
A_{\{1,2\}}$, and consider the group
$J$ with the following presentation:

\smallskip
 {\noindent \bf Generators:} $X, s$ (where $s$ represents
$\diag(-1,-1,1,\dots,1 )$).

\smallskip
{\noindent\bf Relations:}
\begin{itemize}
\item []
\begin{enumerate}
\item $R$.
\item $   s^2 = 1$.
\item $ [s,H] = 1.$
\item $ss^\sigma  s ^{\sigma^2}= 1$. 
\end{enumerate}
\end{itemize}
\smallskip
There is an obvious surjection $\pi\colon J\to W_n$.  
We can view $A=\langle X\rangle   \le J.$  By (3), 
${ n\choose 2}\ge |s^A|\ge |\pi(s^A) |={n\choose 2}$, so that $s^A$ can
be identified with the 2-sets in $I=\{1,\dots,n\}$.  Thus, there are
well-defined elements  
$s_{ij} = s_{ji} \in s^A$  for all distinct $i,j\in I$.

By (4),  $s_{1j}s_{jk}s_{k1}=1$ whenever $1 ,j,k$  are  distinct.  Since all 
 $s_{ij}$ are involutions, it follows that  $s=s_{12}$ commutes with   all  
 $s_{1j}$, and hence also with  all $s_{jk}$,
so that $N:=\langle s^A\rangle$ is elementary abelian.  Then $J=AN$ has order
$|W_n|$. 

Now  Theorem~\ref{Alt(n) and Sym(n) 100 down to 7}  yields the stated
numbers of generators and relations for $n\ge4$.  For $n=4$ or 5 we  
instead  use (\ref{A4 A5}).~\qed   



\section{Rank 1 groups}
\label{Rank 1 groups }

\subsection{Steinberg presentation} 
\label{Steinberg presentation}
Each rank 1 group $G$ we   consider has a
Borel subgroup $B=U\semi \langle h \rangle $, with   $U$ a
$p$-group.  There is an involution  $t$ (mod $Z(G)$ in the case
$\SL(2,q)$ with $q$ odd) such that $h^t=h^{-1} $
(or $h^{-q }$ in the unitary case).  The {\em Steinberg presentation}  for
these groups \cite[Sec.~4]{St2} consists of    
\begin{itemize}
\item a presentation for $B$,
\item a presentation for  $ \langle h ,t \rangle $, and 
\item $|U|-1$ relations of the
form
\begin{equation}
\label{Steinberg relation}
u_0^t=u_1h_0tu_2,
\end{equation}
with $u_0,u_1, u_2$ nontrivial
elements of
$U$ and
$h_0\in \langle h \rangle $ (one relation for each choice of $u_0$).  
\end{itemize}

\subsection{Polynomial notation}
\label{Polynomial notation}

Our groups will always come equipped with various elements having  names
such as $u$ or $h$.
 For any polynomial $g(x) = \sum_0^e g_i x^i\in
\Z [x]$,   $0\le g_i<p$,  
define    powers as follows:
\begin{equation}
\label{how polynomials act}
[[u^{g (x)  }]]_h =  (u^{g_0})  (u^{g_1})^{h ^{ 1}}
\cdots (u^{g_e})^{h ^{ e}} ,
\end{equation}
so that
\begin{equation}
\label{Horner for polynomials}
[[u^{g (x)  }]]_h =
u^{g_0} h^{-1}u^{g_1} h^{-1} u^{g_2} \cdots h^{-1} u^{g_e}
h^{e}
\end{equation}
by  ``Horner's Rule''  \cite[(4.14)]{GKKL}
(compare (\ref{Horner})).

As in \cite[Sec.~ 4.3]{GKKL}, we need to be careful about rearranging the
terms  in (\ref{how polynomials act}) when   not all of the indicated 
conjugates of
$u$  commute.

\subsection{\boldmath $ \SL(2,q)$ }
\label{$SL(2,q)$} 
\addcontentsline{toc}{subsection}{\protect\tocsubsection{}{\thesubsection}{$
\SL(2,q)$}}
\addtocontents{toc}{\SkipTocEntry}

In \cite{CRW}  there is  a presentation for  $ \PSL(2,q)$  with at
most 13 relations; and it follows readily from that  presentation  that
$ \SL(2,q)$  has one with at most $ 17$ relations.  
We now provide    presentations  having fewer relations, 
based on the  matrices 
\begin{equation}
\label{generators}
\mbox{$u=
   \begin{pmatrix} 
        1 & 1 \\
        0 & 1 \\
     \end{pmatrix}
   $,
   $~t=
   \begin{pmatrix}
       \,\,\,\,0 & 1 \\
        -1 & 0 \\
     \end{pmatrix}
~ $
   and $~h =
   \begin{pmatrix} 
        \z^{-1} & 0 \\
     \!\!\!   0 & \z \\
     \end{pmatrix}
   $. }
\end{equation}  
\begin{theorem} 
\label{SL2 10 relations}
  $\SL(2,q)$  and $ \PSL(2,q)$   have   presentations 
with  $3$ generators$,$
$9$ relations  and bit-length $O( \log q)$. 
When $q$ is even $\PSL(2,q)$ has a presentation with $3$ generators$,$ $5$
relations   and bit-length $O( \log q)$.

\end{theorem}
\begin{proof}
By \cite[pp.~137-138]{CoMo}, if $q\le 9 $ then  $\SL(2,q)$  and $\PSL(2,q)$
have  presentations with 2 generators and at most 4 relations.  
 Assume that  $q>9$,  let $\zeta$ be a generator of $\F_q^*$, and
let $k,l \in\Z$ be such that
$\z^{2k} = \z^{2l}+1$ and $\F_q = \F_p[\z^{2k}]$
(as in 
\cite[Section~3.5.1]{GKKL}). 

Set $d=\gcd(k,l)$. 
Then $\F_q = \F_p[\z^{2d}]$. 

Let  $m(x) \in \F_p[x]$ be the  minimal
polynomial of
$\zeta^{2d}$. If $\c\in \F_q$, let $g_\c(x)\in \F_p[x]$ satisfy   
$g_\c(\zeta^{2d})=\c$  and $\deg  g_\c< \deg   m$.

We will show that $ G =\SL(2,q)$  is isomorphic to the group $J$ having the
following presentation. 

\smallskip
{\noindent\bf
Generators:} 
$u,t,h$.
\smallskip

  {\noindent\bf
Relations: }
\begin{itemize}
\item 
[]
\begin{enumerate}
\item $u^p=1$.
\item $u^{h^k} = u u^{h^l} = u^{h^l} u$.
\item $[[u^{m(x)}]]_{h^d}=1$  in the notation of
 $\eqn{how polynomials act} $. \vspace{2pt}
\item  $u^h=[[u^{g_ {\z^2}(x)}]]_{ h^d} $. 
 \vspace{2pt}
\item $[t^2,u]=1$ (or $t^2=1$  in the case $\PSL(2,q)$ with
$q$ odd).
\item $h^t=h^{-1}$.
\item $  t  = u u^t u $.  	
\item  $h t = [[ u^{  g_{   \z^{-1}  }(x) }    ]]_{ h^d}  \  [[ u^
{g_{\z}(x)} ]]_{ h^d}^{\hspace{1pt} t} \ [[u^{ g_{   \z^{-1}
}(x)}]]_{ h^d}    $.

    \end{enumerate}
\end{itemize} 
\smallskip

Matrix calculations using \eqn{generators} easily show that there is a
surjection  $J\to G$.  
By (1), (2)  and  \cite[
Lemma~4.1]{GKKL} (compare \cite{Bau,CR1,CRW}), 
$U:=\<u^{\langle h^k,h^l \rangle}\>$ is elementary abelian;  since $d=\gcd(k,l)$ we have   
 $U=\<u^{\< h^d \>} \>$.
By (1) and (3)  we can identify $U$ with the additive group of $\F_q$
in such a way that  
$h^d$ acts as multiplication by $\zeta^{2 d}$. 
By (4), 
 $h$ acts on $U$ as a transformation of order $(q-1)/(2,q-1)$. 
(Note that $U$ is defined  using  $h^k$ and $ h^l $  rather than $h$ 
so that 
\cite[ Lemma~4.1]{GKKL} can be used.)

By (7) and (8), $J=\<U,U^t\>$, and  $J$ is perfect since $U=[U,h]$.
Moreover, using (6) we see that  $z:=h^{(q-1)/(2,q-1)}$ is inverted by $t$,
centralizes
$U$ and
$U^t$, and hence is an element of  $Z(J)$ having order 1 or 2.

Thus, 
$\<u,h\> /\<z\>$ is isomorphic to a  Borel subgroup of
$\PSL(2,q)$.
By (6), $\<h,t\>/\<z\>$ is   dihedral  of order
$2(q-1)/(2,q-1)$.

We already know that 
$\<h\>$ acts on   the
nontrivial elements of $U$ with at  most 2 orbits, with orbit
representatives $u^1$ and $[[ u^ {g_{   \z^{}}(x)} ]]_{{ h^d}}$ 
if $q$ is odd. 
As in the proof of \cite[Sec.~4.4.1]{GKKL},   (7) and (8) provide the 
relations
\eqn{Steinberg relation} required to let us deduce that $J/\<z\>\cong
\PSL(2,q)$.

Now $J$ is a perfect central extension of  $\PSL(2,q)$, and hence is
$\SL(2,q)$ or  $\PSL(2,q)$.  Finally, (5) distinguishes between these
groups when $q$ is odd.
The bit-length of the presentation  is clear from (\ref{Horner for
polynomials}).

Finally, if $q$ is even there are significant simplifications.  
We may assume that $k=d=1$, so that  $h^d=h$ acts on $U=\<u^{\langle h^k,
h^l \rangle} \>$; the induced automorphism has order   
 $q-1$ by  (3), and (4) can be deleted.     
Relation (5) can be deleted  since  (1) and   (7)
imply that
$t^2=1$;  and (8) is not needed since
$\langle h
\rangle$ has only one  orbit on the nontrivial elements of $U$.
\end{proof}

\Remark 
\label{lengths in SL2}
Every element of $\SL(2,q)$ has bit-length $O(\log q)$ in our generators.  
\rm  For, this is true of all elements of $U$  by (\ref{Horner for
polynomials}), while $\SL(2,q)=UU^tU$.

\begin{Remark}  \rm
If $ d =1$ then relation (4) can be removed since then  
$h^d=h$ acts as multiplication by $\zeta^{2  }$, by (3). 

 In 
\cite[Section~3.5.1]{GKKL} it was observed that we can choose $k=1,$
$l=1$  or $k=2$.  Thus, $d\le 2$ for some choice of $k$ and $l$.  If
$q$ is even then $d=1$.  If $q\equiv 3$ (mod 4) and $k=2$ we can change
$\z
$ to
$-\z^2$ in order to obtain $k=1  $ and hence  $d=1. $
We can also prove that there are  choices for $\z, k, l$ that yield 
$d=1  $  when $q\equiv 5$ (mod 8), but we do not know how to
obtain such choices   in general.

\end{Remark}

 \Remark  \rm
\label{degree q+2}
Now that we have the notation  in  (\ref{how polynomials act}),
we can give  more examples along the lines of 
Examples~\ref{more AGL1 examples}(\ref{AGL(1,9)}) and
(\ref{AGL(1,9)/2}).  Let $q$ be a power of an odd prime $p$ such that 
$(3,q-1)=1$; we may assume that $q>5$.

(1$'$)  {\em   $S_{q+2}$ has a presentation
with  $2$ generators and $6$ relations.}  For, if $\zeta$ is as above,
 $\zeta+1=\zeta^s$, and  $g(x)$ is the minimal polynomial of $\zeta$
over $\F_p$,  then 
$\AGL(1,q)\cong \< u,h\mid  u^p=h^{q-1}, u^{h^s}=u u^h=u^h u,
[[u^{g(x)}]]_h=1\>$;   this is proved
as above, using \cite[
Lemma~4.1]{GKKL}.  Then Proposition~\ref{general number of relations}
provides the stated presentation.
\smallskip

 (3$'$)  {\em  $A_{q+2}$ has a presentation
with  $2$ generators and $6$ relations -- only $5$ relations if} 
$q\equiv 3$ (mod 4).   This
time let
$k$,  $ l$  and $m(x)$ be as in the  proof of the preceding theorem, 
and obtain
$\AGL(1,q)^{(2)}\cong \<  u,h \mid u^p=h^{(q-1)/2}, 
u^{h^k} = u u^{h^l} = u^{h^l} u,
[[u^{m(x)}]]_h=1 \>$, after which 
Proposition~\ref{general number of relations} provides the stated
presentation (with $\rho = (2,(q-1)/2)).$
 

\subsection{Unitary  groups}
\label{Unitary  groups}
\label{Other rank 1 groups}

We will obtain presentations for 3-dimensional unitary\break
 groups by  taking   the  presentations in
\cite[Sec.~4.4.2]{GKKL}  and deleting the portions that were needed to
produce   short presentations.
We use    matrices of the form 
\begin{equation}
\label{unitary matrices}
u= \begin{pmatrix}
        1& \a & \b\\
          0 & 1 &  \!\!\!\!-\bar \a \\
   0 & 0 &1\\
     \end{pmatrix}\!\!, 
w= \begin{pmatrix}
        1&  0 & \c\\
          0 & 1 &  0 \\
   0 & 0 &1\\
     \end{pmatrix}\!\!, 
t=\begin{pmatrix}
        0& \,0 & 1 \\
          0 &\!\!\! -1   & 0\\
   1 & \,0 & 0\\
     \end{pmatrix} \!\! ,
h= \begin{pmatrix}
        \bar \z^{-1} & 0 & 0\\
          0 &  \bar \z/\z    &0 \\
          0 & 0 &  \z\\
     \end{pmatrix} \! \!  
\end{equation}
with  $\a,\b,\c\in \F_{q^2}$ arbitrary such that $\b+\bar \b = - \a \bar \a \ne0
$ and
$\c=-\bar
\c\ne0$.

\begin{theorem} 
\label{SU3 10 relations}
 $\SU(3,q)$  and  $\PSU(3,q)$   have   presentations 
with  $3$ generators$,$ $23$ and 
$24$ relations$,$ respectively$,$   and bit-length $O( \log q)$. 
\end{theorem} 

\proof
Let
  $\z$ be a generator of  $\F_{q^2} ^*$.   
As in \cite[Sec.~4.4.2]{GKKL}, we will assume that $q\ne 2, 3, 5$ and use
elements 
 $a=\z^k,b=\z^l $, where   $1\le k,l<q^2$, 
such that
$$
\begin{array}{lllll}
a^{2q-1}+ b^{2q-1} =1~ \mbox{ and } ~
 a^{q+1}+ b^{q+1} =1, \vspace{2pt}
\\
\mbox{$\F_q=\F_p[a^{q+1}] $, and }\vspace{3pt}
\\
\mbox{$\F_{q^2}=\F_p[a^{2 q-1 }]$ 
if $q$ is odd while 
$\F_q=\F_p[a^{2 q-1 }]$ if  $q$ is even. }
\end{array}
$$

Let $d=\gcd(k,l)$. 

 If $\c\in \F_{q^2}$  write $\c':=\c^{q+1}$  and  $  \c'':=\c^{2q-1}$, 
and also let   $m_\gamma (x)$  denote its minimal polynomial over
$\F_p$. If $\delta \in \F_p[\gamma]$, let
$f_{\delta;\gamma}(x)\in\F_p[x]$ with $f _{\delta;\gamma}(\gamma)=\delta$
and $\deg f _{\delta;\gamma}  < \deg m  _{ \gamma} $
(compare \cite[Sec.~4.4.2]{GKKL}).
\smallskip

The required presentation is as follows:

\smallskip

 {\noindent \bf Generators:} $u,  
h ,t$.
\smallskip

{\noindent\bf Relations:}
\begin{itemize}
\item []
\begin{enumerate} 
\item $w= w^{h^k} w^{h^l} =w^{h^l} w^{h^k} ,$
where $w$ is defined by 
$u = u^{h^k}u^{h^l}w$.
\item $w^p=1$.
\item $[[w^{m_{a'}(x) } ]]_{h^d }=1$.
\item $[[w^{f_{\z';a'}(x) } ]]_{h^d }=w^{h}$.
\item $u = 
u^{h^l}u^{h^k} w_1 $.
\item  
$[u,w] = [u^{h^k},w] =1.$

\item $u^p=w_2$.
\item $[[u^{m_{a''}(x) } ]]_{ h^d}=w_3$.
%

\item[(9$'$)]
$[[u^{f_{\z'';a''}(x) } ]]_{{ h^d} }=u^{}w_4$  if $q$ is odd.
\item[(9$''$)]
$
([[u^{f_{\alpha ;a''}(x) } ]]_{{ h^d} })^{h }
[[u^{f_{\beta;a''}(x) } ]]_{{ h^d}}
=u^{h ^2}w_5$
if $q$ is even and $\z''$ satisfies $\z''{}^2 = \alpha \z''{} +\beta $ for
$\alpha ,\beta \in \F_q$. 
\item[(10)]
$[u,u^{ h }]= w _6 $  and  $[u^{ h^k},u^{ h }]= w_7$  if $q$ is even.

\item [(11)]  $ t^2 = 1$.  
\item [(12)]  $h ^t =  h ^{-q}  $.
\item [(13)]  $ u_i^t=   u_{i1}  h_{i}   t   u_{i2}
 $ for
$1\le i\le 7$,   relations due to Hulpke and Seress
\cite{HS}.
\item [(14)] $1 =  v_{ 1 1}  v_{ 12}^t
 v_{ 13}  t  $   and $h  =  v_{ 2 1}  v_{2 2}^t
 v_{ 23}  t  $.

\end{enumerate}
Here,  the elements $  w_i$ are specific
words in of $ w^{\langle  {h^d}\rangle}$,  the elements         
 $u_{ ij},v_{ i j }$ are specific words in  $ u^{\langle  {h^d}\rangle}$,
and the elements
$h_i$  are specific powers of $ h $; all depend on the initial 
choice of    $u $ and  $\z$.

\end{itemize} 

\smallskip

Note that $ \langle w^{\langle {h^d}\rangle}\rangle$ is defined using
${h^d}$ rather than $h$ so that  \cite[Sec.~4.1]{GKKL} can be used
together  with (1) and (5). 
 See \cite[Sec.~4.4.2]{GKKL} for a proof that this is, indeed,  a
presentation of $\SU(3, q) $.  
As in \cite[Sec.~4.4.2]{GKKL}, one further relation of bit-length
$O(\log q)$ produces  a
presentation for
$\PSU(3,q)$.   \qed

\Remark
\label{SU3 lengths}\rm
Let $U:=\langle u^{\langle {h^d}\rangle}\rangle$.
By applying (\ref{Horner for
polynomials})  to both $Z(U)=\langle w^{\langle {h^d}\rangle}\rangle$
and $U/Z(U)$, we find that all elements of $U$ have bit-length $O(\log q)$,
and hence the same holds for $\SU(3,q)=UU^t UU^t U$.

\subsection{Suzuki  groups}
\label{Suzuki  groups}
 
The short presentation in 
\cite[Section~4.4.3]{GKKL} uses 7 generators and 43  relations for 
$\Sz(q)$.  This  can be used for  $^2\!F_4(q)$
in  Section~\ref{Rank 2 groups}.
Nevertheless, we note that there is an easy modification similar to what
occurred for $\SL(2,q)$ and $\SU(3,q)$: three generators are powers of a
fourth, allowing us to decrease  to 4 generators and 31 relations. 


\section{{$\SL(3,q)$}}
\label{SL3 100} 

The groups $\SL(3,q)$   will reappear
 more often in the rest of the proof than any other rank 2 groups: we will 
use  $\SL(3,q)$  to obtain first all higher-dimensional groups $\SL(n,q)$, and
 then all higher-dimensional classical groups.
Therefore we will be more explicit with these groups than the other rank
2 groups  (cf. Section~\ref{Rank 2 groups}).

\begin{theorem}
\label{SL3 for 100}
$\SL(3,q)$   has a  presentation 
  with  $4$ generators$,$ $1 4$  relations and bit-length $O(\log q)$. 

\end{theorem}
\begin{proof}
 When   $q\le 9$ there are presentations 
 with 2 generators and at most 10
 relations 
 \cite{CMY,CR3}.  
(These cases also can be handled directly by a slight
variation on the following  
approach; compare the end  of Section~\ref{most
rank 2}.)
Hence,     
 we will assume that
$q>9$.

 Let
$\SL(2,q)\cong L=\langle X\mid R \rangle$, with  $X =\{u,t,h\}$  and
$\<\z\>=\F^*_q$ as in the proof of  Theorem~\ref{SL2 10 relations}.
We view the elements of  $\SL(3,q)$ as  matrices,
with $L$ consisting of the matrices 
$\bigl( \begin{smallmatrix}
* &0 \\
0&1
\end{smallmatrix} \bigr)$.

We will show that $ G=\SL(3,q)$  is isomorphic to the group $J$ having the
following presentation. 

\smallskip
{\noindent\bf
Generators:} \rm
$X$ and  $c$ (corresponding to the permutation matrix 
$ \Bigl( \begin{smallmatrix}
0 &1 &0 \\
0 &0 &1 \\
1 &0 &0 \\
\end{smallmatrix} \Bigr)$
acting as $(3,2,1)$). 

\smallskip
  {\noindent\bf
Relations: }
\begin{itemize}
\item 
[] 

\begin{enumerate}
\item $R$.
\item $c^t=t^2c^2$.
\item $h h^c h^{c^2} =1$.
\item $u^{h^c} = u^{\diag(1,\zeta^{-1})} ,$ written as a word in $X$. 
(The matrix  
$ \diag(1,\zeta^{-1}) $
 is not in $L$  if $q$ is odd, but  it  can be viewed as  inducing an automorphism of
$L$.)

\item $[u^{c},u ] = (u^t)^{c^2}$.
\item $[u^{c},u^t]=1$.
\end{enumerate}

\end{itemize}
\smallskip

Matrix calculations easily show that there is a
 surjection   $J\to G $.   There is a
subgroup of $J$ we identify with $L=\<X\>$.  We separate the argument into four
steps.
\smallskip

1.  \emph{Computations in   $\<t,c\>$.}
The relations $t^4=1$ and (2)  imply that
\begin{equation}
\label{tct}
tct=c^2, \quad t^{-1}c^2 t^{-1}=c,\quad t^{-1}c^3t=t^{-1}c^2 t^{-1}tct=cc^2,
\end{equation}
and hence
\begin{equation}
\label{$t^c$}
t^c =  c^{-1}tctt^{-1}  \! = c^{-1} c^2 t^{-1}  \! = ct^{-1},
~
t^{c^2} = (c t^{-1})^c  = t^{-1}c.
\end{equation}
It follows that  
\begin{equation}
\label{c cubed ts}
t^2 (t^2)^{c^2}(t^2)^{c}= 
t^2 \cdot t^{-1} c t^{-1} c\cdot  c t^{-1} c t^{-1} =
t c (t^{-1} c^2 t^{-1}) c t^{-1}=
t c c c t^{-1} = c^{3}.
\end{equation}


\quad
  
2. \emph{The action of $h,h^c,h^{c^2},  c^3$ on $L$.}  
By \eqn{tct} and  (3), 
$
(h^c)^t =    h^{t^{-1} c^2}  
= (h^{-1})^{c^2} = 
hh^c.
$
 Consequently, conjugating both sides of (4) by $t$ gives
\vspace{4pt}

\hspace{-20pt}
$
\begin{array}{rllll}
 \left(u^{h^c} \right)^t \hspace{-8pt}
& \! = \left(u^t\right)^{h^{ct}} =
\left(u^t\right)^{hh^c} = 
\left(u^{th}\right)^{h^c}
\vspace{2pt}
\\
\left(u^{\diag(1,\zeta^{-1})}\right)^t \hspace{-8pt} 
& \! = 
\left(u^t\right)^{\diag(1,\zeta^{-1})^t}  \! \!=\hspace{-.5pt}
\left(u^t\right)^{\diag(\zeta^{-1},1)} \!=\hspace{-.5pt}
\left(u^t\right)^{h\, \diag(1,\zeta^{-1})}   \!=\hspace{-.5pt}
\left(u^{th}\right)^{\diag(1,\zeta^{-1})} \!.%
\end{array}
$

\noindent
(Recall that $h=\diag(\zeta^{-1},\zeta)$ in Section~\ref{$SL(2,q)$}.)
Thus, by (4),  $h^c$ acts on $L=\langle u, u^{th} \rangle $ as conjugation by
${\diag(1,\zeta^{-1})}$.  
We know how $h$  acts on $L$ since $h\in L$. 
Now (3)
implies that
$h^{c^2}$ acts on $L$ as conjugation by ${\diag(\zeta^{-1},1)}$.

If $q$ is odd then 
  $t^2=h^{(q-1)/2}$  is in  $Z( L)$.  It follows  that $t^2$, 
$(t^2)^c =(h^c)^{(q-1)/2}  $ and~%
$(t^2)^{c^2}=$ $ (h^{c^2})^{(q-1)/2}  $ 
act on $L$ as they should:~%
 as conjugation by 1,   $ \!\diag(1,-1)$ and~%
  $\diag(1,-1)$, respectively.
Now  \eqn{c cubed ts} implies that   $c^3 = t^2 (t^2)^{c^2} (t^2)^{c} $ acts
trivially on
$L=\< X \>$ and hence on  $J=\< X,c\>$.
This also holds trivially if $q$ is even. 
 
\medskip

3. \emph{The elements  $e_{ij}(\lambda) $.} 
For all integers $m$ and all $\lambda\in \F_q$, write
$$
\begin{array}{llllllll}
e_{12}(\z^m) :=  u^{(h^m)^c}\!, ~ e_{12}(0):=1 \qquad
&
e_{21}(\lambda) :=   e_{12}(-\lambda)  ^t
& \vspace{2pt}
\\
e_{23}(\lambda):= e_{12}(\lambda)^ c 
&
e_{32}(\lambda):= e_{21}(\lambda)^c \vspace{2pt}
\\ 
e_{31}(\lambda):= e_{12}(\lambda)^{c^2}
&
e_{13}(\lambda):= e_{21}(\lambda)^{c^2}.
\end{array}
$$ 
Then $e_{12}(1)=u$, $e_{23}(1)=u^c$, and $e_{12}( \F_q) $ is an elementary
abelian subgroup of $ L$ by  (4).  Then we also have  
$e_{2 1}(\F_q) <L$.   By (\ref{$t^c$}), $c=tt^c\in \<X\cup X^c\>$, so
that 
$J$   is  generated by the  elements
$e_{ij}(\lambda)$.

Clearly, $\< c\>$ acts on the set of subgroups $e_{ij}(\F_q)$;  in fact  
  $\<t,c\>$ acts   as the symmetric group $S_3$ on subscripts (this is the
Weyl  group of $G$).  For example,  
by \eqn{tct}, 
 $e_{23}(\F_q) ^t=e_{12}(\F_q) ^{ct}
=e_{12}(\F_q) ^{t^{-1}c^ 2}=e_{13}(\F_q)  $; and (as we have seen) $t^2$ acts
correctly on
$L^{c^{-1}}=L^{c^2}$ and hence also on $e_{12}(\F_q)^{c^2}=e_{31}(\F_q)$.
Similarly, $t$ acts correctly on each $e_{ij}(\F_q)$.

\medskip

4. \emph{Verification of the Steinberg relations} 
(see  \cite[Section~5.1 or
5.2]{GKKL}). The relations  
$$
e_{ij}(\lambda)  e_{ij}(\mu) = e_{ij}(\lambda+\mu)
$$
follow from the corresponding relation in $L$ (with $\{i,j\}= \{1,2\} $)  by
conjugating with
$t$ and 
$c$. We will deduce the remaining relations from (5) and (6) by conjugating with 
  $t$,  $c$, $h$, $h^c$ and $h^{c^2}$; 
we have seen that these act on the set of subgroups $e_{ij} (\F_q)$ as they
do in
$G$.

We have $[e_{23}(1),e_{21}(1)]=[u^{c},u^t]=1$ by (6).
Conjugating by $h^x (h^c)^y$ we obtain  
$$
[e_{23}(\zeta^{2y-x}),e_{21}(\zeta^{-2x+y})] = 1.
$$
This does not cover all relations of the form
$
[e_{23}(\lambda),e_{21}(\mu)] = 1
$
since $\det\bigl( \begin{smallmatrix}
-1&2\\ -2&1
\end{smallmatrix} \bigr)=3$,
but this does imply  that 
$$
[e_{23}(1),e_{21}(\mu^3)] = 1\quad \mbox{ for all } \mu \in \F_q.
$$
The additive subgroup of $\F_q $ generated  by the cubes 
in $\F_q ^*$ is closed under
multiplication, and so is a subfield of size $\ge 1+(q-1)/3$ and
hence is all of $\F_q$  since $ q\ne 4$.  It follows that  
$$
[e_{23}(1),e_{21}(\mu)] = 1\quad \mbox{for all } \mu\in \F_q.
$$
Conjugating this by all $h^x$ yields  all relations of the form
$$
[e_{23}(\lambda),e_{21}(\mu)] = 1 \quad \mbox{ for all } \lambda ,\mu\in \F_q.
$$
Conjugating by $\<t,c\>$, we obtain
\begin{equation}
\label{commuting}
[e_{kl}(\lambda),e_{km}(\mu)]=[e_{km}(\lambda),e_{lm}(\mu)]=1
\quad \mbox{for all distinct } k,l,m \mbox{ and all  } \lambda, \mu.
\end{equation}

Similarly, (5) implies  that 
$$
[e_{23}(1),e_{12}(1)] =[u^{c},u]=
(u^t)^{c^2}=e_{13}(-1).
$$
Since $t$ acts correctly on each $e_{ij}(\F_q)$, conjugating by $t$ gives
$
[e_{13}(-1),e_{21}(1)]   
=e_{23}(1).
$
Conjugating by $h^x (h^c)^y$ we obtain 
$$
[e_{23}(\zeta^{2y-x}),e_{12}(\zeta^{2x-y})] =
e_{13}(-\zeta^{x+y}),
\quad
[e_{13}(\zeta^{x+y}),e_{32}(\zeta^{x-2y})] =
e_{12}(\zeta^{2x-y}).
$$
Once again these do  not cover all relations of the form
$$
[e_{23}(\lambda),e_{12}(\mu)]=e_{13}(-\lambda\mu)
\quad \mbox{and} \quad
[e_{13}(\lambda),e_{32}(\mu)]=e_{12}(\lambda\mu).
$$
However, using the standard identity $[x, ab ] = [x, b][x, a]^b$, together
with  
\eqn{commuting} and the fact that the cubes in $\F_q^* $ generate $\F_q
$ under addition, we deduce all such relations. 

Conjugating by $\<t,c\>$ yields all remaining Steinberg relations. 
Then $J$ is a homomorphic image  of $G=\SL(3,q)$, and hence   $J\cong G$. 
By Theorem~\ref{SL2 10 relations}  and Remark~\ref{lengths in SL2}, our
presentation has the required bit-length;
note that $c$ is the product of an element of $L$ and an element of $L^c$.  
The presentation    uses $|R|+5=1 4$ relations.%
\end{proof}

In order to obtain a presentation for  $\PSL(3,q)$ when $m:=(q-1)/3$ is
an  integer, add the relation 
$h^m(h^{2m})^c  =1.$  

Note that the  computations in the above proof  would be
 considerably simpler if we added the relation
$c^3=1$;  this is a relation we will have available  in the proof of
Theorem~\ref{SLn 100}.

\Remark 
\label{lengths in SL3}
\rm Using Remark~\ref{lengths in SL2} we see that  {\em each element 
of
$\SL(3,q)$ has bit-length
$O(\log q)$ in our generators.}

\Remark
\label{generating SL3}
 \rm For future use  we will need to know that \em
$\SL(3,q)=\langle a,a^{(1,2,3)}\rangle$ for some $a\in L $ 
of order $q+1$.  
\rm 
 For,   $\langle a,a^{(1,2,3)}\rangle$ is an irreducible subgroup
of $\SL(3,q)$.  Using the order of $a$ and the list of possible subgroups
\cite{Mitchell,Hartley} proves the claim.



\section{$\SL(n,q)$}
\label{PSL(n,q)} 

We now turn to the general case of Theorem~B for the  groups
$\PSL(n,q)$,  using a variation on the approach in Section~\ref{SL3 100}.

\begin{theorem}
\label{SLn 100}
Let $n\ge4.$
\begin{itemize}
\item[(a)] 
$\SL(n,q)$   has a  presentation with  $7$ generators$,$ $25$ relations
and bit-length\break $O(\log n
 +\log q)$. 
\item[(b)] $\PSL(n,q ) \!$ has~a presentation with $7$ generators$,$  $26$
relations  and bit-length $O(\log n +\log q).$%

  \item[(c)]  $\SL(4,q) $ has a    presentation  with  $6$
 generators$,$
 $2 1$ relations and bit-length\break
  $O( \log q)$.%
 \item[(d)]  $\SL(5,q) $ has a  presentation  with  $6$
 generators$, $
 $2 2$ relations and bit-length\break
 $O( \log q)$.%
    
\end{itemize}
\end{theorem}


\Proof 
We use two presentations:

\smallskip
\begin{itemize}
 \item
The  presentation  $\langle X \mid R \rangle$ for
$F=\SL(3,q)
$ in Theorem~\ref{SL3 for 100}.  
We view  $F$ as the group of matrices 
$\bigl( \begin{smallmatrix}
* &0 \\
0&I
\end{smallmatrix} \bigr)$ in $G=\SL(n ,q)$ 
and only write the upper left $3\times 3$ block. 
\smallskip 

\item
  The  presentation $\< Y\mid S\>$  for  $ T=A_n$ in
Theorem~\ref{Alt(n) and Sym(n) 100 down to 7}, where
$T$ acts  on
$\{1, \ldots, n\}$   (with $X$ and $Y$ disjoint).  We view $T$ as
permutation  matrices. 
\end {itemize}  
\smallskip

We will also use the  subgroup 
$L = \SL(2,q)$ of $F$ consisting of matrices  in the upper left $2\times 2 $
block.   We   use the following elements:

\smallskip
\smallskip
\begin{itemize}
\item [] $c =  
\begin{pmatrix}
0&0&1 \\
1&0&0 \\
0&1&0
\end{pmatrix}\! , \,f = 
\begin{pmatrix}
0&1&\,0 \\
 1&0&\,0 \\
0&0&\!\!\!-1
\end{pmatrix}\in F$, 
\smallskip

\item  [] $a \in L$ such that $\langle a, a^f \rangle = L$, 

\item  []
$(1,2,3),(1,3)(2,4)\in T$, 
and

 
\item  []$\sigma$  and $\tau =(1,2)(3,4)  $  in $ T$   interchanging $1$ and
$2$ and generating   the set-stabilizer~$T_{\{1,2\}}$  of $\{1,2\}$~in~$T$.

\smallskip
\hspace{-14pt}
{\em Bit-length}: $c,$  $ f$ and $a$ have bit-length $O(\log q)$ using 
Remark~\ref{lengths in SL2}.  
We may assume that     $\sigma$ is a
cycle of length $n-2$ or $n-3$ on $\{3,\dots,n\}$; both $\sigma$  and
$\tau$   can be viewed as 
words in  $Y$  of  bit-length 
$ O(\log n)$ (by Remark~\ref{small bit-length}). 

\end{itemize}
\smallskip

We will show that $ G$  is isomorphic to the group $J$ having the
following presentation. 

\smallskip
{\noindent \bf Generators:} $X,Y$. 
\smallskip

{\noindent\bf Relations:}
\begin{itemize}
\item []
\begin{enumerate} 
\item $R$.
\item $S$.
\item $c=(1,2,3)$.
\item $a^\sigma = a^f$.
\item $a^\tau = a^f$.
\item $(a^f)^\sigma = a$.
\item $[a, a^{(1,3)(2,4)}]=1$.
\item $[a^f, a^{(1,3)(2,4)}]=1$ (needed only when $n$ is 4 or 5). 
\end{enumerate}
\end{itemize}
\smallskip
 
As usual,  there is a  surjection $\pi\colon J\to G$, and $J$  has subgroups
we will identify with $ F =\<X\>$ and $T=\<Y\>$. Since $\tau$ has order 2, (5)
implies that $(a^f)^\tau = a$.
Hence, by (4)--(6), $\<\sigma,\tau \>$ normalizes  
$L$, inducing the same automorphism group as $\<f\>$ on $L$.
In particular,  elements of $\<\sigma,\tau \>  $ that fix 1 and 2 must centralize
$L$, while   elements interchanging 1 and 2 act as $f$. 

It follows that  
$|L^T |\le  \binom{n}{2}$; as usual, we use $\pi$ to obtain 
equality. Then $ L^T $  can be identified
with the set of all 2-sets  of $\{1, \ldots, n\}$.  Its subset 
$L^{\< c \>}$ consists of 3 subgroups corresponding to the 2-sets
in $\{1,2,3 \}$.
Consequently, any two distinct members of 
$L^T$  can be conjugated by a single element of $T$ to one of the pairs
$L,L^{ (1,2,3)}$ or $L,L^{ (1,3)(2,4)}$.   
Here $\<L,L^{ (1,2,3)} \>=  \<L,L^{ c} \>=  F$ by (3).   

We will use (7)--(8) to show that  
$[L,L^{ (1,3)(2,4)}]=1$. 
By our comment about elements of $\<\sigma,\tau \>$, we have 
$a^{(1,2)(5,6)}=a^f$ and 
$a^{(3,4)(5,6)}=a$.
Then
\begin{equation}
\label{only for n<6}
\begin{split}
1&=[a,a^{(1,3)(2,4)}]^{(1,2)(5,6)} \hspace{4.8pt}
  =   
[a^f,
(a^{(3,4)(5,6)})^{(1,3)(2,4)}]  
=
[a^f, a^{(1,3)(2,4)}] 
\\
1&=[a,a^{(1,3)(2,4)}]^{(1,2)(3,4)} \hspace{4.9pt} 
  =   [a^f, (a^f)^{(1,3)(2,4)}]
\\
1&=[a^f,a^{(1,3)(2,4)}]^{(1,2)(3,4)} 
=[a , (a^f)^{(1,3)(2,4)}],
\end{split}
\end{equation} 
where the first equations explain the comment in  (8).
 
Thus, any
two distinct members of 
$L^T$  either generate  a conjugate of 
$F=\SL(3,q)$  or  commute.
%
Consequently,  $N:=\langle L^T  \rangle\cong G$ 
(see  
\cite[Sections~5.1  or 5.2]{GKKL}  for the Steinberg presentation).
Moreover, $N\unlhd 	J$, and
$J /N$ is a homomorphic image of $\langle  Y \rangle \cong A_{n}$ in which $ c
$  is sent to 1. Thus, $J/N=1$.

The bit-length follows easily from those of  $\<X\mid R\>$ and  $\<Y\mid S\>$.

We still need to count the number of relations.  If $n\ge 6$ then we have
used 
$14 +7+5$ relations by Theorems~\ref{Alt(n) and Sym(n) 100 down to 7} and
 \ref{SL3 for 100}.
However, {\em we can remove one generator and one relation as follows}.  
In  Theorem~\ref{SL3 for 100}
we used a generator ``$c$''  (corresponding to the permutation
matrix acting as
$(3,2,1)$).  Hence we replace that  element $c$   by a word in $Y$
representing $(1,2,3)\in T$   in the relations appearing in that theorem.  Then
(3) is implied by the {\em new} presentation  
 $\<X\mid R\>$.

If $n=4$ or 5  we use  the presentation (\ref{A4 A5})
 with only 2 generators and 3 relations.   
Then the preceding presentation requires  only 
$14+4+6 -1  $   relations.
Moreover, when $n=4$ we can delete $\sigma$
entirely, saving two further relations (4) and (6); and 
when $n=5$ we can
choose $\sigma$ of order 2  and  delete (6).

Finally, we need to add one further relation in order to obtain  $\PSL(n,q)$. 
Let $h_{i,j}$ be the matrix with $\z$ and $\z^{-1}$ in positions $i$ and $j$, and
1 elsewhere.  Let $m=(q-1)/(d,q-1)$. If $n$ is odd, use 
 Remark~\ref{cycles for
p+2} to obtain the $(n-2)$-cycle $(2,\dots,n-1)$, and then the additional
relation 
$h_{ 1,n }^{m}\big( h_{ 2,n }^m {(2,\dots,n-1)} \big)^{n-2}=1$
 produces  $\PSL(n,q)$ with the required bit-length. The case
$n$   even is similar.
\qed

\section{Remaining rank 2 groups}
\label{Rank 2 groups}
In this section we will provide  
presentations  required in Theorem~\ref{B} for some of the rank 2 groups of Lie
type.
Since ${\rm
P}\Omega^-(6,q)\cong\PSU(4,q)$, and we will handle all unitary groups in a
different manner in Theorem~\ref{unitary 100}, we only need to consider the
groups $\Sp(4,q)$, $G_2(q)$, $^3\! D_4(q)$ and $^2\! F_4(q)$.
Note that the last three groups   do  not appear inductively as Levi factors of 
any higher rank groups of Lie type.

 \subsection{\boldmath $\Sp(4,q)$, $G_2(q)$ and $^3 \! D_4(q)$ }  
 \label{most rank 2}  

\addcontentsline{toc}{subsection}{\protect\tocsubsection{}{\thesubsection}{$\Sp(4,q)$,
$G_2(q)$ and $^3 \! D_4(q)$ }}
\addtocontents{toc}{\SkipTocEntry}

Here the Weyl group
is dihedral of  order $2m=8$ or  12. 

\begin{theorem}
\label{rank 2s}

\begin{itemize}
\item [ (a)]  $\Sp(4,q)$ has a  presentation  
with   $6$
generators$,$ $3 5$ relations and  bit-length $O(\log q)$.

\item [ (b)]  $\PSp(4,q) \cong \Omega(5,q)$  
has a  presentation with   $6$
generators$,$ $3 6$ relations and  bit-length $O(\log q)$.

\item [ (c)]  $G_2(q)$ and $^3\! D_4(q)$ 
  have   presentations 
with   $6$
generators$,$  $4 0$ relations and  bit-length $O(\log q)$.

\end{itemize}
\end{theorem}

 \Proof
(a)  
The   root system $\Phi $ of $G=\Sp(4,q)$  has  $8$ roots, half of them 
long and  half short.  Let  $\Pi=\{  \a_1,\a_2\}$ be a set of  fundamental
 roots  with $\a_1$ long; there are   corresponding rank 1 groups 
$L_{\a_i}\cong \SL(2,q) $.

We assume that $ q>9$ until the end of this proof,  and use the
presentation   
$\<X_i\mid R_i\>$ of $L_{\a_i}$   in Theorem~\ref{SL2 10 relations},
with 
\begin{equation}
\label{Xi}
X_i = \{ u_{\a_i},r_i,h_i\}\!,\, i=1,2  
\end{equation}
(here we are using $ r_i$ instead of $t_i$ in order to  approximate 
standard Lie notation; we assume that $X_1$ and $X_2$ are disjoint).
 The action of $h_i$ on
$u_{\a_i}$ is given in $R_i$, and 
$ U_{\a_i}:=\< u_{\a_i  } ^{\<h_i\>} \>$ has order
$q$.  Let $
U_{\a_i}^\#:=U_{\a_i}\backslash\{1\}.$
There are $(2,q-1)$  orbits
of $\<h_i \>$ on   $ U_{\a_i}^\#$, with orbit
representatives 
$u_{\a_i ,1}:=u_{\a_i}$ and $u_{\a_i , 2}$ (where  
$u_{\a_i , 2} : =u_{\a_i }$ if $q$ is even).   If $q$ is odd then
$\<h_1,h_2\>$ has 2 orbits on both 
$ U_{\a_1 }^\# \times  U_{ \a_2 }^\# $ and 
 $U_{\a_1 +\a_2 }^\#\times  U_{ \a_2}^\# $,
with respective orbit representatives
$(u_{\a_1  , a}, u_{  \a_2})$, $a=1,2$, and 
$(u_{\a_1 +\a_2 , a}, u_{ \a_2 ,a})$, $a=1,2$ 
(see  \cite[Lemma~5.3]{GKKL}).

The root groups $U_\a ,\a\in \Phi$, will be built into our  
presentation.

We will show that $ G$  is isomorphic to the group $J$ having the
following presentation. 

\smallskip
{\noindent\bf
Generators:} \rm
$X_1\cup X_2$.

\smallskip
  {\noindent\bf
Relations: }
\begin{itemize}
\item 
[] 

\begin{enumerate}
\item $R_1\cup R_2$. 
\item $h_1^{r_2}$, $h_2^{r_1}$  and $w^4$ are explicit words in
$\{h_1,h_2\}$,  where 
$w:=r_1r_2$. 
 
%

%

 \vspace{2pt}
\item  
 $(u_{\a_i  })^{h_j^{w^n}}$  
is an explicit word in  $ u_{\a_i}^{\< h_i\>}$   for $i\ne j$ and
$0\le n<4$.

 \vspace{2pt}
Notation: $ U_{\a_i}:=\<  u_{\a_i  } ^{\<h_i\>} \>  $ for
$i=1,2$,  
$H:=\<h_1,h_2\>$, $N:=\<r_1,r_2, H\>$, and $W:=N/H$ (which is isomorphic
to $D_8$).
 \vspace{2pt}

 Label the subgroups $( U_{\a_i} ) ^{w^n}$, $0\le n<4$, in the usual way
as
$$\mbox{
$U_\a$, $\a\in  \Phi:=\{\pm\a_1,\pm\a_1, \pm(\a_1+a_2), \pm(\a_1+2a_2)  
\}$ }
$$
\cite[p.~46]{Ca}, so that  $\Phi =\a_1^{\{w^n\mid 0\le n
<4\}}\cup
\a_2^{\{w^n\mid 0\le n <4\}}$.
 \vspace{2pt}

For $i=1,2$ let $u_{\a_i, 1}$ be another name for $u_{\a_i }$, and  let
$u_{\a_i, 2}$ stand  for an explicit  word in $ u_{\a_i  } ^{\<h_i\>}$
corresponding to an element of $ G$ described above.

If   $\a=\a_i^{w^n}$ for (unique) $i\in \{1,2\}$ and 
$0\le n<4$, write  $u_{\a ,a}:= ( u_{{\a_i, a}} ) ^{w^n}$ for $a=1,2$.

\vspace{2pt}

%




 \vspace{2pt}

\item 
\begin{enumerate}
\item  
 $[u_{\a_1}, u_{\a_1+2\a_2 }] = [u_{\a_1 }, u_{\a_1+ \a_2 }] =1$.
\item 
 $[u_{\a_1 ,a}, u_{ \a_2 }]$ is an explicit word in
$ u_{\a_1}^{\< h_1\>w^3} \cup u_{\a_2}^{\< h_2\>w^3}  \subset
U_{\a_1+2\a_2 } \cup  U_{\a_1+
\a_2 } $ for $a=1,2$.
\item 
 $[u_{\a_2 ,a}, u_{ \a_1 +\a_2,a}]$ is an explicit word in
$ u_{\a_1}^{\< h_1\>w^3} \subset U_{\a_1+2\a_2 }$  for $a=1,2$.
  
\end{enumerate}

\end{enumerate}

\end{itemize}
\smallskip

First of all, the explicit words  mentioned in these relations are 
obtained in 
$G$, and   have bit-length $O(\log q)$ in the
generators by Remark~\ref{lengths in SL2}.   It follows that  this
presentation has bit-length $O(\log q)$, and that there is a surjection   
$\pi \colon J\to G$.
Moreover, this presentation has  
$|X_1|+|X_2|=6 $ generators and  at most  
$|R_1|+|R_2| +3 + 8   + 6 = 3 5 $   relations.

By (1) there is a subgroup $L_{\a_i} \cong \SL(2,q)$ of $J$ we can
identify with
$\<X_i\>$. 
 
 By (2),
 $H\unlhd N$  and 
$W\cong  D_8$, as stated in (3).

Using  $L_{\a_i} $ we see that   $ U_{\a_i} ^{{r_i}}=U_{-\a_i} $.
By (3), $H$ normalizes each subgroup $U_\a=\<u_\a^H\>$.  It follows that
$W$ acts on $\Phi$  in the natural manner:
 there are  $8$ root groups $U_{\a}$,
$\a\in\Phi$, labeled as in (3) and permuted by $N$.

Fix roots $\a$ and $\b\ne \pm\a$. Then  
$H$ acts by conjugation on the set of all 
Chevalley commutator relations  \cite[p.~47]{GLS}, which we write as 
\begin{equation}
\label {commutator rels}
\mbox{$[y_{\a },y_{\b }]=\prod_\g
v_\g~$ with $~  y_\a\in U_\a ^\#$, $  y_\b\in U_\b^\#$, $v_\c\in
U_\c$},
\end{equation}
where $\c$ runs through all  positive
integral combinations  of $\a$ and $\b$. 
Clearly  (4) provides us with instances of \eqn {commutator rels}
corresponding to each of the four $W$-orbits on  pairs $\{  \a,  \b
\},$  $\b\ne
\pm\a$. 

If $[y_{\a },y_{\b }]=1$ with  $\{y_{\a },y_{\b }\}=\{u_{\a_1},
u_{\a_1+2\a_2 }\}$ or  $\{u_{\a_1 }, u_{\a_1+ \a_2 } \}$ in (4a), 
 then $[U_{\a },U_{\b }]=  [\<y_{\a }^H\>,\<y_{\b }^H\>]=1$.
(For, the Chevalley
commutator relations for $G$ \cite[p.~47]{GLS} show that  the
corresponding fundamental subgroups of
$G$ commute. In particular, already in $\pi( J ) =G$    suitable
elements of 
$h_1^N\cup h_2^N$  act independently on $U_{\a }$  and  $U_{\b }$.)

If $[y_{\a },y_{\b }]\ne 1$ for the pairs $\{y_{\a },y_{\b }\}$  in (4b,c)
then, by  
\cite[Lemma~5.3]{GKKL}, we obtain all relations \eqn {commutator rels} by
conjugating the   ones in (4b,c)  by elements of $H$.

At this point we have verified the 
Steinberg-Curtis-Tits relations mentioned in 
Section~\ref{Preliminaries}.  Thus, $J$ is a homomorphic image of  the
universal group of Lie type for  $G$, and hence is $G$
\cite[pp. 312-313]{GLS}. 
\smallskip

(b) When $q$ is odd, we just need one  further relation to kill the
center of
$\Sp(4,q)$.%

(c) These groups  are handled in a manner  similar to what was
done in (a), replacing the number 4 by 6 at suitable places in order to
obtain the Weyl group
$D_{12}$.  We sketch this very briefly since these groups do not arise
 in higher rank settings.

\smallskip

Once again we have fundamental subgroups $L_{\a_1}\cong \SL(2,q) $
and  $L_{\a_2}\cong \SL(2,q) $ or  $\SL(2,q^3) $, where the latter occurs
for short $\a_2$ in  $^3\!D_4(q)$.  We use versions of (1)--(4),
labeling the roots as in \cite[p.~46]{Ca}.  However,
in order to avoid any use of $(3,q-1)$  elements  $u_{\a_1 ,a}$ we proceed
as follows in (4).  The long root groups $U_\a$ of $G$ generate a subgroup
$\SL(3,q)$. Relations (5) and (6) from our presentation in
Theorem~\ref{SL3 for 100}, with $c:=w^2$, are instances of \eqn{commutator
rels}; hence these are the only relations of this sort needed for long
$\a,\b$. By
\cite[Lemma~5.3]{GKKL}, this is the only situation where we might have
needed to use several elements $u_{\a_1 ,a}$. 

This time $W=N/H$ has seven  orbits  
on   pairs $\{  \a,  \b \}$,  $\b\ne \pm\a$.  Counting, we see that there
are only $4 0$ relations. Note the  significant savings due to 
not needing  more  than one  element of any root group
$U_{\a}$.

This proves the theorem when $q>9$. 

  It remains to make some straightforward
remarks about the remaining cases.  For each of these, if $L_{\a_i}$ is a
central extension of $\PSL(2,q)$ then it has a presentation  $\<Y_i\mid S_i\>$
using 2 generators and at most 3 relations
\cite{Sun,CMY,CR3}.  We
do not have to be concerned about lengths of words in these bounded
situations, so we can assume that each of the elements 
$u_{\a_i },$  $u_{\a_i,  a}, r_i,h_i$ is replaced by a word in $Y_i$. 
Once this has been done,  we can again write the relations   (1)--(4) in
our new generating set, and then  our  previous argument goes through
without difficulty.~\qed

   \subsection{\boldmath $^2\hspace{-1pt} F_4(q)$}
\label{^2F_4(q)}  
\addcontentsline{toc}{subsection}{\protect\tocsubsection{}{\thesubsection}{$^2\hspace{-1pt}
F_4(q)$}}
\addtocontents{toc}{\SkipTocEntry}

Here $q=2^{2e+1}>2$. 
There  is no root system  in the classical sense, but there are 16
``root groups'' $U_i$, $1\le i\le 16$.  There are  rank 1 groups 
$L_1=\Sz(q)$ and $L_2=\SL(2,q)$, and we use   the  presentation   
$\<X_i\mid R_i\>$ for  $L_{i}$    in 
Section~\ref{Suzuki  groups} or Theorem~\ref{SL2 10 relations}.  
If $i=2$ then
\eqn{Xi} holds, and  $U_2:=\<  u_2 ^{\<h_2\>} \>$ is
elementary abelian of order $q$.
On the other hand, $X_1$ has size 7 and contains elements
$u_1,r_1,h_1$ behaving essentially as before, except that this time
$U_1 =\<  u_1 ^{\<h_1\>}\>$  is  nonabelian of order
$q^2$ with 
 $Z(U_1)=\Phi(U_1)=\< (u_1^2)^{\<h_1\>} \>$ (where $\Phi(U_1)$ denotes
the Frattini subgroup of $U_1$).

We use the following presentation. 

\smallskip
{\noindent\bf
Generators:} \rm
$X_1\cup X_2$.

\smallskip
  {\noindent\bf
Relations: }
\begin{itemize}
\item 
[] 

\begin{enumerate}
\item $R_1\cup R_2$. 
\item $w^8=1$, where  $w:=r_1 r_2$.  
\item $h_1^{r_2}$ and $h_2^{r_1}$   are explicit words in
$\{h_1,h_2\}$.

 \vspace{2pt}
\item  $(u_{\a_i})^{h_j^{w^n}}$  is an explicit word in  $
u_{\a_i}^{\< h_i\>}$ for $\{i,j\}=\{1,2\}$ and  $0\le n<8$.
 \vspace{2pt}

Notation:
Let 
$u_{i+n}:= ( u_{ i} ) ^{w^n}$ for $i=1,2$ and $1\le n<8$.

\item  
$[u_{i}, u_j]$ is an explicit product  of words in  $u_k^{\<h_1, h_2\>}$,
$i<k<j$, for the pairs $(i,j)$  with $i=1$ and $ 2\le j\le 8$, or $i=2$
and $j=4,$ $ 6$ or 8. 

\end{enumerate}

\end{itemize}
\smallskip

In the presented group  $J$ there are again subgroups $L_i$ we can
identify with $\<X_i\>$, $i=1,2$.  Also, $\<r_1,r_2\>$ is dihedral of
order 16 by (2), and normalizes $H:=\<h_1, h_2\>$ by (3).
It then follows from  (4) that $H$ normalizes each subgroup 
$U_{i+n}:= ( U_{ i} ) ^{w^n}$ for $i=1,2$ and $1\le n<8$.

 As in the case of the other rank 2 groups,  the known actions  of
$  h_1^{w^n}$ and $ h_2^{w^n} $ in (4)  allow us to deduce from (5) an
additional relation   analogous to \eqn {commutator rels}
for each pair  of nontrivial cosets  of the form
$y_i \Phi(U_i),$ $y_j  \Phi(U_j) $, for $i, j$ as in (5)
and $y_i\in  U_i ,$ $y_j \in  U_j $ (compare
\cite[Lemma~5.3]{GKKL}).
 By using 
 the elementary identity $[x,uv]= [x,v] [x,u]^v$, we see that  these
conjugates of the relations
(5) imply all analogues of   \eqn {commutator rels} for these $i,j$.
Finally, we obtain further
relations  analogous to \eqn {commutator rels}
by conjugating by elements  of  $\<w\>$. 
It is now easy to see that we have all relations required for a 
presentation of $G$  (\cite [p.~412]{Gri}, \cite [p.~105]{BGKLP} and
\cite[p.~48]{GLS}
  give the 10 formulas implicitly contained  in (5)).

There are $1 6$ relations (4) and $10  $  relations (5), for a total  of 
$7+3$  generators and     $43+9 +1+ 2  + 16+10 =8 1$ relations.
As noted in Section~\ref{Suzuki  groups}, these numbers  
can easily be decreased to $4+3$ and $31+9 +1+ 2  + 16+10  $.


   \section{Unitary groups }  

\label{Presentations of unitary groups } 
Since the commutator relations for the odd-dimensional unitary groups are
especially complicated  (see, for example, \cite[Theorem~2.4.5(c)]{GLS}),
we will deal with these groups separately.  In fact, when combined with 
Theorems~\ref{Alt(n) and Sym(n) 100 down to 7} and  \ref{SU3 10 relations},
recent presentations in 
 \cite{BeS} allow us to use surprisingly few generators and relations
(cf. Theorem~\ref{unitary 100}).

\subsection{Phan style presentations}
\label{Phan style presentations}
We will outline the   presentation of $G = \SU(n,q ),$  $n\ge 4$,  in
\cite{BeS}, based on one in 
\cite{Ph}.
In \cite{BeS},  subgroups $ U_1 , U_2\cong \SU(2,q )$ of $\SU(3,q )$  are called 
a {\em standard pair} if  $U_1$ and $U_2$ are the respective stabilizers   in
$\SU(3,q )$ of  perpendicular nonsingular vectors.

Using an orthonormal basis, it is easy to see that  $G $ 
 has subgroups 
$U_i\cong \SU(2,q),$ $ 1\le i\le n-1,$ and 
$U_{i,j},$ $ 1\le i<j \le n-1,$ satisfying the following conditions.
\smallskip
\begin{itemize}
\item [(P1)]  If $|j-i|>1$ then $U_{i,j}$  is a central product of $U_{i }$ 
and $U_{ j}$. 
\item [ (P2)] For $1\le i<n-1$, $U_{i,i+1}\cong \SU(3,q)$, and $U_i, U_{i+1}$
is a standard pair in
$U_{i,i+1 }$.

\item [(P3)]  $G=\<U_{i,j}\mid 1\le i<j \le n-1   \>$.
\end{itemize}
\smallskip

The presentation we will use is the following analogue of the
Curtis-Steinberg-Tits presentation mentioned in
 Section~\ref{Preliminaries}.
\begin{theorem}
\label{Phan theorem}
{\rm\cite{Ph,BeS} }  If {\rm(P1)--(P3)}  hold in a group $G,$  then $G$ is
isomorphic to a factor group of $\SU(n , q) $ in each of the following
situations.

\begin{itemize}
\item [ (a)]   $q >3 $  and  $n\ge 4$.
\item [ (b)]  $q =2$ or $3,$ $n\ge5$  and the following hold$:$
\begin{itemize}
\item [ (1)]  
$\< U_{i,i+1 }, U_{ i+1,i+2 } \>\cong \SU(4,q) $ whenever    $1\le i\le n-3
;$ and

\item [ (2)]  if  $q = 2$ then
\begin{itemize}
\item [ (i)] $[U_i, U_{j,j+1}]=1 $ whenever  $1\le i\le n-1,$ $1\le j\le
n-2 
$  and
$i\ne  j-1,j,j+1,j+2  $$;$ and 
\item [ (ii)] $[U_{i,i+1}, U _{j,j+1}]=1$
whenever  $1\le i\le n-2,$    $1\le j\le n-2$ and
$i\ne  j-2,j-1,j,j+1,j+2  $.

 \end{itemize} 
 \end{itemize}
 \end{itemize}

\end{theorem}

 In
\cite{BeS}, it is remarked that (P1)--(P3)  do not provide a  presentation
of $\SU(n , 2) $.  It is also noted that a standard pair in
$\SU(3,2)$ does not generate that group.

\subsection{Some specific presentations}
\label{Some specific presentations of unitary groups}
~
When $q=2$  or $3$ we need to deal with
some small groups. 
The computer package  MAGMA \cite{CP} contains the following presentation: 
$$
\begin{array}{llll}
\hspace{-12pt} \SU(4,3) = \<x,y\mid x^3& \hspace{-6pt}=y^8=[y^4,x] =(x 
y)^7=[x,(x y^{-1})^7] &  \vspace{2pt}\\&
\hspace{-6pt}=[y^2,x y^2 x y^2 x y^2]=x x^{ {x y x^{-1} y^{-1} x y}}=
(x y x^{-1}y^2)^8\> =1\>.
\end{array}
$$    
We  also need other small cases \cite{Bray}:
$$
\begin{array}{llll}
\SU(6,2) = \< a, b \mid a^2 = b^7  & \hspace{-6pt}
= (ab^3)^{11} = [a,
b]^2  \vspace{2pt}
\\
& \hspace{-6pt} = [a, b^2]^3 = [a, b^3]^3 =(ab)^{33} =
(abab^2ab^3ab^{-3})^2  = 1 \>.
\vspace{2pt}
\\
\SU(5,2)=\<  a, b \mid a^2 = b^5&\hspace{-6pt} = (ab)^{11} = [a, b]^3 =
[a, b^2]^3 = [a, bab]^3 = [a, bab^2]^3 = 1 \>.%
 \vspace{2pt}
\\ 
\SU(4,2)=\<
 a, b \mid a^2 = b^5 &\hspace{-6pt}= (ab)^9 = [a, b]^3 = [a,bab]^2 = 1 
\>.
\end{array}
$$  
The last of these is  \cite[(10.8)]{CMY}.

\subsection{Presentations of unitary groups}
\label{Presentations of unitary groups}
In this section we will prove the following

\begin{theorem}
\label{unitary 100}
Let  $n\ge  4$.

\begin{itemize}
\item[(a)]
 $\SU(n,q)$ has a presentation with   $6$
generators$,$  
$39$ relations  and bit-length\break
 $O(\log n+\log q)$.

\item[(b)]
 $\PSU(n,q)$ has a presentation with   $6$
generators$,$   $40$  relations  and bit-length $O(\log n+\log q)$.

\item[(c)] 
 $\SU(4,q)$ and  $\SU(5,q)$ have presentations with  $6$
generators$,$  
$35$ relations  and bit-length $O( \log q)$.

\item[(d)] 
 $\PSU(4,q)$ and  $\Omega^-(6,q)$ have presentations with  $6$
generators$,$  
$36$ relations  and bit-length $O( \log q)$.

\end{itemize}

\end{theorem}

\Proof
Let 
\smallskip
\begin{itemize} 
\item   $F:=\SU(m,q)$, with $m=3,4 $ or 6 and $m\le n$, and~
\item   $A:=A_n$, acting on 
$\{1,\dots,n\}$.
\smallskip
\end{itemize} 
We view both of these groups as lying in
$G=\SU(n,q)$, using an orthonormal basis of the underlying vector space: $F$
consists of the matrices 
$\bigl( \begin{smallmatrix}
* &0 \\
0&I
\end{smallmatrix} \bigr)$   with an $m\times m$ block in the upper left corner,
and $A$ consists of permutation matrices.

For each $m$ we assume that we have the following additional information:
\smallskip

\begin{itemize} 
\item The presentation    
$\<X\mid R\>$   of $ F$  in 
 Theorem~\ref{SU3 10 relations},
 except in the case of
  the pairs
$(m, q)=(4,2),(4,3)$ or $(6,2)$, in which case
  $\<X\mid R\>$  is  given  in Section~\ref{Some
specific presentations of unitary groups}.

\item The presentation   $ \<Y\mid S\>$  of $A$  in
Theorem~\ref{Alt(n) and Sym(n) 100 down to 7} (where $X$ and $Y$ are
disjoint).

\item   Two generators $g,h$ for $W:=\SU(2,q)<F$, where $W$
consists of the matrices 
$\bigl( \begin{smallmatrix}
* &0 \\
0&I
\end{smallmatrix} \bigr)$  with a $2\times2$ block in the upper left corner,
and  where  $g$ and $h$ can be 
 viewed as  words in $X$ of   
bit-length $O(\log q)$  (using  Remark~\ref{SU3 lengths}).

\item  A word in $X$  of   
bit-length $O(\log q)$ representing the element $c_{(1,2,3)}\in F$ that acts
as the 3-cycle
$(1,2,3)$ on the orthonormal basis (cf.  Remark~\ref{SU3 lengths}).

\item  A word in  $Y$ of   bit-length $O(\log n)$  representing  
$(1,2,3)$  (cf.  Remark~\ref{small
bit-length}).

\item   Permutations  $\sigma= (1,2)(3,4)$  and  $ \tau $ in $A$ that
interchange 1 and 2 and generate the set-stabilizer
$S_{n-2}$ of $\{1 , 2\}$ in $A$, where  $\sigma$  and $\tau$ 
   can be viewed as 
words in $Y$ of   bit-length $O(\log n)$ (using  Remark~\ref{small
bit-length}).

\end{itemize}

\smallskip 
Parts (b) and (d) of the theorem are handled at the end of the proof.

\smallskip 
{\noindent \bf Case  
$q > 3$}:  Here we let   $m=3$.   
We will show that $ G$  is isomorphic to the group $J$ having the
following presentation. 
  In view of the preceding remarks,
this presentation has the desired bit-length.

\smallskip 
{\noindent \bf Generators:} $X\cup Y.$
\smallskip

{\noindent\bf Relations:}
\begin{itemize}
\item []
\begin{enumerate}
\item $R\cup  S $.
\item  $c_{(1,2,3)}=(1,2,3)$.      
\item $g^\sigma=g^\tau= g^{\!^{\bigl( \begin{smallmatrix}
0&1 \\
1 &0
\end{smallmatrix} \bigr)}}$ and
$h^\sigma=h^\tau  = h^{\!^{\bigl( \begin{smallmatrix}
0&1 \\
1 &0
\end{smallmatrix} \bigr)}} \!   $.  
(The matrix  
$ {\bigl( \begin{smallmatrix}
0&1 \\
1 &0
\end{smallmatrix} \bigr)}  $
 is not in $W$  if $q$ is odd, but can be viewed as  inducing an automorphism of
$W$.)
\item  $[W, W^{(13)(24)} ]=1.$ 
 
\end{enumerate}
\end{itemize}
\smallskip

Since there is a surjection $\pi\colon J\to G$,  there are subgroups of $J$ we
can identify with 
$F=\<X\>$ and $A=\<Y\>$.

By (3),
$|W^{A}|\le
\binom{n}{2}$, and using $\pi$ we see that   
$W^{A }$  consists of
$\binom{n}{2}$ subgroups  we will call
$W_{i,j}=W_{j,i} $, $1\le i < j\le n$,  
where $W_{i,j} \le F$ for $1\le i < j\le 3$.  If $i,j,k,l$ are distinct, then
(4) and   the 4-transitivity of
$A$ imply that  
$[W_{i,j},W_{k,l}]=1 $   and  
$\<W_{i,j},W_{i,k}\>\cong \<W_{1,2},W_{1,3}\> = F $.

Let  $U_i:=W_{i,i+1}$ and $U_{i,j} := \<U_i,U_j\>.$  These subgroups 
satisfy (P1)--(P2) and hence, by  Theorem~\ref{Phan theorem},   
 $N :=\<  U_{i,j} \mid 1\le i< j\le n-1\>$   is a homomorphic image of 
$G$. 

We claim that  $N\unlhd G$.  For,   $W_{1,3}\le \<W_{1,2},W_{2,3} \>$. By the  
3-transitivity of $A$, it follows that 
$W_{i,k}\le \<W_{i,j},W_{j,k} \>$ for all distinct $ i,j,k$.  
By induction, if $i<j +1$   then 
$$
\begin{array}{llll}
W_{i, j} &\hspace{-6pt}\le \<W_{i,j-1}, W_{j-1,j}  \> \vspace{2pt}
\\ & \hspace{-6pt}\le \<W_{i,i+1},
\dots ,W_{j-2,j-1}, W_{j-1,j}\> = \<U_{i}, \dots , U_{j-2}, U_{j-1}\>  \le
N.
\end{array}
$$
 Thus, $N=\< F^A\> \unlhd J.$
By (2),
  $J/N$ is a quotient of $A_n$ in which $(1,2,3)$ is mapped to 1.  Thus,   
$J/N=1$. 

Total: $|X|+|Y| = 3+3$ generators and $|R|+|S|+9  = 23+7+9   $ relations.

This proves (a) when $q>3$.  For (c), we 
use  the presentation (\ref{A4 A5})
 with only 2 generators and 3 relations.    
Hence, the preceding presentation requires only 
$23 + 3+  9$   relations.
\smallskip

{\noindent \bf Case  
$q =3$}:  This time  we let  $m=4$ and   use the presentation $\<X\mid R\>$
for
$\SU(4,3)$ given  in Section~\ref{Some specific presentations of unitary
groups}.  We assume  that $n\ge5$, and  that we   have 
\smallskip
\begin{itemize}
\item     Words in $X$ representing the elements   $c_{(1,2,3)}, c_{( 2,3
,4)}\in F$  that act as the indicated permutations on the orthonormal
basis.
\end{itemize}

\smallskip  
We will show that $ G$  is isomorphic to the group $J$ having the
following presentation. 
\smallskip  
 
{\noindent \bf Generators:} $X\cup Y.$

\smallskip 
{\noindent\bf Relations:}
\begin{itemize}
\item []
\begin{enumerate}
\item   [(1$'$)] $R\cup  S $.
 
\item   [(2$'$)] $c_{(1,2,3)} = (1,2,3),$ $  c_{(2,3,4)}=(2,3,4)$.  
 
 \item   [(3$'$)] $g^\sigma= g^\tau= g^{\!^{\bigl( \begin{smallmatrix}
0&1 \\
1 &0
\end{smallmatrix} \bigr)}}$ and
$ h^\sigma=h^\tau  = h^{\!^{\bigl( \begin{smallmatrix}
0&1 \\
1 &0
\end{smallmatrix} \bigr)}} \!   $.        
 
\end{enumerate}
\end{itemize}
\smallskip

By  (2$'$) and  (3$'$), 
 $g^\sigma=g^\tau \in W$ and
$h^\sigma=h^\tau  \in W $.  As before, it follows that 
$W^{A}$  consists   of  $\binom{n}{2}$ subgroups  
$W_{i,j}=W_{j,i} $, $1\le i < j\le n$,  
where $W_{i,j} \le F$ for $1\le i < j\le 4$.

The previous relation    (4)   follows  from   the corresponding
relation  in  $F=\SU(4,q)$.

If $i,j,k,l$ are distinct, then (4)   and   the 4-transitivity of
$A$ give 
$[W_{i,j},W_{k,l}]=1 $,    
$\<W_{i,j},W_{i,k}\>\cong \<W_{1,2},W_{1,3}\> \cong \SU(3,q) $
and
$\<W_{i,j},W_{i,k},W_{i,l}\>\cong \<W_{1,2},W_{1,3},W_{1,4}\> = F  $.

Once again, the subgroups  $U_i:=W_{i,i+1}$ and $U_{i,j}: = \<U_i,U_j\> $
of  $J$  satisfy (P1)--(P2).  They also behave as in  Theorem~\ref{Phan
theorem}(b1) since $\< U_{i,i+1}, U_{i+1,i+2} \>\cong \< U_{1,2}, U_{2,3}
\>=F$.

Once again, the subgroup $N$   generated by all $U_{i,j}$  is isomorphic to $G$. 
As before, $N$  is normal in $J=\<X , Y\>$ and hence is $J$.

Total: $|X|+|Y|= 2+3$ generators and $|R|+|S|+6 =8 +7 + 6 $ relations.
 
\smallskip

{\noindent \bf Case  
$q = 2$}:  This time  we let  $m=6$.  Using the 
presentations  in Section~\ref{Some specific presentations of unitary groups} 
we may assume  that $n\ge 7$, and  that we   have 
\smallskip
\begin{itemize}
\item Generators $g',h'$ for $V: =\SU(3,2)$ as words in $X$; 
\item     Words in $X$ representing elements  $c_{(1,2,3)} , c_{( 2,3
,4,5,6)}\in F$ that act as the indicated permutations on the orthonormal basis;
and 

\item   Two  permutations  $\sigma'= (1,2)(4,5)$  and  $ \tau'$ that generate the
set-stabilizer  of $\{1 , 2,3\}$ in $A$ and  can be viewed as 
words in $Y$ of   bit-length $O(\log n)$ (cf. Remark~\ref{small
bit-length}).    
\end{itemize}   

\smallskip

We will show that $ G$  is isomorphic to the group $J$ having the
following presentation. 
\smallskip

{\noindent \bf Generators:} $X\cup Y.$
\smallskip

{\noindent\bf Relations:}
\begin{itemize}
\item []
\begin{enumerate}
 \item  [(1$''$)] $R\cup  S $.
\item [(2$''$)] $c_{(1,2,3)} = (1,2,3),$ $ c_{(2,3,4,5,6)}  = {(2,3,4,5,6)}
$. 
%
 \item  [(3$''$)]  $ g^\sigma=g^\tau= g^{\!^{\bigl( \begin{smallmatrix}
0&1 \\
1 &0
\end{smallmatrix} \bigr)}}$ and
$g^\sigma= h^\tau  = h^{\!^{\bigl( \begin{smallmatrix}
0&1 \\
1 &0
\end{smallmatrix} \bigr)}} \!   $.

\item [(4$''$)]    
$g'{}^{\tau'} =g'{}^{c[\tau']}$ and $h'{}^{\tau'}  =g'{}^{c[\tau']}$,
where $c[\tau'] $ denotes the automorphism of $V$ induced by the permutation
matrix for the restriction of $\tau'$ to the first 3 coordinates.

\end{enumerate}
\end{itemize}
\smallskip 

   As before, it follows from (2$''$) and  (3$''$)    that 
$W^{A}$  consists of  $\binom{n}{2}$ subgroups   $W_{i,j}=W_{j,i} $, $1\le i <
j\le n$,   where $W_{i,j} <F$ for $1\le i < j\le 6$. The previous relation  
(4) follows from  the corresponding relation  in  $F=\SU(6,2)$.

If $i,j,k,l$ are distinct,  then
 $W_{i,j} $ commutes with $W_{k,l}  $  by (4) and the 4-transitivity of $A$.

By (2$''$), $\sigma'\in\<c_{(1,2,3)},c_{(2,3,4,5,6)}\><F$, so that 
$V^{\sigma'}=V$.
As before, (4$''$) then implies that  $V^A$  consists of $\binom {n}{3}$
subgroups 
$V_{i,j,k}$ for distinct $i,j,k\in \{1,\dots,n\}$, where $V_{i,j,k} \le F$ for
$1\le i,j,k\le 6$  and  $W_{i,j} \le   V_{i,j,k}$.

In view of the  transitivity properties  of $A$,  if $i,j,k,l,r,s$ are
distinct then\break
$\<V_{i,j,k},V_{i,j,l}
\>
\cong 
 \<V_{1,2,3},V_{1,2,4} \> \cong \SU(4,2)$   and
 $[V_{i,j,k} , V_{l,r,s}]=1  $.

Let $U_i:=W_{i,i+1}$,  $U_{i,i+1}: = V_{i,i+1,i+2}$  and $U_{i,j}:
=\<U_i,U_j\>$  iff $|i- j|>1 .$

The  subgroups  $U_{i,j}$ satisfy  (P1)--(P2).   
They also behave as in Theorem~\ref{Phan theorem}(b1)
since $\< U_{i,i+1}, U_{i+1,i+2} \>\cong \< U_{1,2}, U_{2,3}
\>=\SU(4,2)$.~The conditions  in  Theorem~\ref{Phan theorem}(b2) also 
hold  since they hold for the subgroups $U_i, U_{i,i+1}, U_{j,j+1}  $ 
that lie in
$F$.

 Hence, the subgroup $N$   generated by all $U_{i,j}$  is isomorphic to $G$.  As
before, $N$ is normal in $J=\<X , Y\>$ and hence is $J$.

Total: $|X|+|Y| =2+3$ generators and $|R|+|S|+8 =  8 +7 + 8 $ relations.

\smallskip
This proves (a) and (c)  for all $q$.  For (b) and (d) we note that    a
presentation of each group in question is obtained as  in the proof of
Theorem~\ref{SLn 100}  by adding one new relation  of bit-length $O(\log
n+\log q)$ to a presentation in (a) or (c).
\qed


\section{General case}
\label{General case}
We now complete the proofs of Theorems~\ref{A} and  \ref{B}.
\begin{theorem} \label{classical theorem}
All finite simple groups        of
Lie type and rank $n \ge 3 $ over $\F_q$  have presentations with at most  
$1 4$ generators$,$  $7 6$ relations and bit-length $O(\log n+
\log q)$.

\end{theorem}
\Proof
We   use a   variation on the methods in 
\cite[Section~6.2]{GKKL}.  By Theorems~\ref{SL3 for
100}, \ref{SLn 100} and
\ref{unitary 100}, we may assume that
$G$ is neither a special linear nor  unitary group.

Until the end of the proof we assume that our simple group  $G$ is the
simply connected group of the given type. 
As usual,  $\Pi=\{ \alpha_1, \ldots, \alpha_n\}$ is the set of fundamental 
roots  of $G$, and for each $i$ there are  root groups $U_{\pm
\a_i}$.

\smallskip

  {\noindent\bf Case 1: $G$ is a classical group.}  Number $\Pi$ as in
\cite[Section~6.2]{GKKL}: 
the subsystem $\{ \alpha_1, \ldots, \alpha_{n-1}\}$
is of type $A_{n-1}$, $\alpha_n$ is an end node root and is connected
to only one  root $\alpha_j$ in the Dynkin diagram
(here $j=n-1$ except for
type $D_n$, when $j=n-2$).
Let 
\begin{equation}
\label{G1 G2}
\mbox{$G_1= \< U_{\pm \alpha_i}\mid  1 \le i < n \>$, 
$G_2 =  \< U_{\pm \alpha_n}, U_{\pm \alpha_j}  \> $,
$L_2 = \< U_{\pm \alpha_n}\>$,   $L= \< U_{\pm \alpha_j} \>$.}
\end{equation}
Then $G_1$ has type $A_{n-1}$ and $G_2$ is a rank 2 classical group.

Let $L_1$ be the subgroup of $G_1$ generated by the root groups
that commute with $L_2$.  Then  $L_1$ is of type $A_{ n-2}$  unless
$G$ has type $D_n$, in which case $L_1$ is of type $A_{ 1} \times
A_{n-3}$.

We will use the following presentations  (with $X$ and $Y$ disjoint):
\begin{itemize}

\item The presentation $\langle X \mid  R \rangle$  for  $G_1$
in Theorem~\ref{SLn 100}; and 
 
\item The presentation $\langle Y \mid   S  \rangle$ for  $G_2$ in 
Theorems~\ref{rank 2s}  and
 \ref{unitary 100}(d).~(The latter presentation is needed for the universal
cover of $\Omega^-(6,q)$ rather than for general unitary groups.)

\end{itemize}

 Choose two generators
for $L_1$, two for $L_2$ and two for $L$,  all
viewed as words in $X$ or $Y$.

{\em Bit-length}: Generators for the subgroup $  \< U_{\pm \alpha_1} \> \cong
\SL(2,q)$ are used in the presentation 
 $\langle X\mid  R \rangle$  (cf. Section~\ref{PSL(n,q)}).  Thus, by
Remark~\ref{lengths in SL2},  each element of   $ \< U_{\pm\alpha_1}\>$  
has  bit-length $O( \log q)$ in $X$.  

Recall that the presentation in Theorem~\ref{Alt(n) and Sym(n) 100 down to 7} 
was used in the proof of   Theorem~\ref{SLn 100}.  Hence,  the cycles
in  Remark~\ref{small bit-length} have bit-length 
$O(  \log n)$ in $X$.
Since $L_2 = \< U_{\pm \alpha_n}\>$ and   $L= \< U_{\pm \alpha_j} \>$  can
each be obtained by conjugating $ \< U_{\pm\alpha_1}\>$   using one of those
cycles, it follows that our generators for these groups have bit-length
$O( \log n+\log q)$ in   $X$.
Similarly, in view of  Remark~\ref{lengths in SL2} the proof  of
Theorem~\ref{rank 2s} shows that  our two generators for $L$ have bit-length 
$O(  \log q)$ in terms of $Y$.  

We need to be more careful with our choice of generators of $L_1$ in order to
achieve the desired bit-length.    
We will assume that $G$ does not have  type
$D_n$; the omitted case is very similar.
 If $n-1=2$, we can use any pair of generators by
Remark~\ref{lengths in SL2}.  Assume that
$n-1\ge 3$.  We temporarily write elements of  $L_1\cong \SL_{n-1}$ 
using $(n-1)\times (n-1)$ matrices. View 
$\< U_{\pm \alpha_1}\>=\bigl(
\begin{smallmatrix} * &0 \\
0&I
\end{smallmatrix} \bigr)$ with $2\times 2 $ blocks~$*$, inside 
$\SL(3,q)=\< U_{\pm \alpha_1}, U_{\pm \alpha_2}\>=\bigl( \begin{smallmatrix}
* &0 \\ 0&I
\end{smallmatrix} \bigr)$  with $3\times 3 $ blocks $*$.
Choose $a\in \< U_{\pm \alpha_1}\>$ such that $\SL(3,q)= \langle a,
a^{(1,2,3)}\rangle$ (cf. Remark~\ref{generating SL3}).
The monomial transformation  
$g:=(1,2,\dots,n-1)$ or 
$(1,2,\dots,n-1, -1, -2, \dots)=  r_{12}\,\, (2,\dots,n-1) $ is in
$L _1$;  and $g$   has bit-length $O(\log n + \log q)$ in $X$  by
the preceding paragraph. Note that $a^{(1,2,3)}= a^{g}$: since
$n-1\ge 3$, $(1,2,3)$ and
$g$ agree on the first two coordinates.
Then $  L_1 \ge \langle a,  {g}   \rangle\ge 
\langle \langle a, a^{g} \rangle ^{\langle  {g}   \rangle} \rangle 
= \langle\langle a, a^{(1,2,3)} \rangle  ^{\langle  {g}   \rangle} \rangle
=\langle\SL(3,q)  ^{\langle  {g}   \rangle} \rangle
= L_1$, where
$a$  and  $ a^{g}
$   have  bit-length $O(\log n + \log q)$   in $X$.

\smallskip

We will show that $ G$  is isomorphic to the group $J$ having the
following presentation. 
  In view of the preceding remarks, this presentation has the
desired bit-length.%

\smallskip
{\noindent \bf Generators:} $X,Y$. 
\smallskip

{\noindent\bf Relations:}
\begin{itemize}
\item []
\begin{enumerate} 
\item $R \cup S$. 

\item $[L_1,L_2]=1$.
 
More precisely,  impose the
four obvious commutation relations using  pairs of words in $X$ or $Y$ that
map onto the chosen  generators for   $L_1$  or  $L_2$.

\item Identify
the copies of $L$ in $G_1$ and $G_2$.  

More precisely,  take the two generators for $L$, viewed as words in   $X$ and
also in 
$Y$, and impose the equality of the corresponding words.
\end{enumerate}
\end{itemize}
\smallskip

We claim that $J \cong  G$.  For, clearly  $J$ surjects
onto
$G$, and hence  we may assume that  $G_1=\<X\>$ and $G_2=\<Y\> $, $L$, $L_1$ and 
$L_2$  are subgroups of
$J$. Clearly $J$ is generated by the
fundamental  root groups contained in $G_1$ or $G_2$.  Any two
of these root  groups satisfy the Curtis-Steinberg-Tits relations
(see the references in Section~\ref{Preliminaries}): 
either they are both in $G_1$ or $G_2$, or they commute since  $[L_1,
L_2]=1$.   Thus, $J$ is a homomorphic image of the universal finite group
of Lie type  of the given type, which proves the claim. 

By  Theorems~\ref{SL3 for 100}, \ref{SLn 100}, \ref{rank 2s}  and
 \ref{unitary 100}(c), if the type is not $D_n$ then the number of generators is
$|X|+|Y|\le 7 + 6  $  and the number of relations is
$|R|+|S| + 6
\le 25 +  3 6  +6 =67  $ (since
$4$ relations are required to ensure that $[L_1,L_2]=1$ and $2$  to identify
the copies of $L$).  For type $D_n$ these numbers become 
$7+4$ and $25 + 14 + 6 $.

As in \cite{BGKLP,GKKL},  in each case  we can kill the center of $G$
with at most one 
additional relation  of bit-length $O(\log n + \log q)$,
except for some of the groups $D_n$, where two relations may be needed.
Thus, in all cases we use at most $1 3$ generators and at most 
$69$ relations for all simple classical groups.

This proves the theorem for the classical groups.  Howeover,  for use in
$F_4(q)$, when  
 $G\cong \Sp(6,q)$ we can  use the presentation  in Theorem~\ref{SL3 for
100} instead of the one in Theorem~\ref{SLn 100}.  This time  
$|X|+|Y|= 4+  6$ and  $|R|+|S| + 6 =   14 + 3 6 +6 = 56$.
 
\smallskip

  {\noindent\bf Case 2: $G$ is an exceptional group.}   
We modify the above argument   slightly.
If $G$ is the universal cover of  $E_n(q)$ with $6 \le n \le 8$, for a suitable
numbering of
$\Pi$ we can define
$G_1,G_2,L_1,L_2, L$ essentially as in \eqn{G1 G2}:   $G_1$ still has type
$A_{n-1}$,
$G_2$ has type $A_2$, and this time $L_1$ has type  $A_2 \times A_{n-4}$.  Since
$L_1$   is still generated by $2$ elements, precisely as above we obtain 
$|X|+|Y|= 7 + 4
$ and 
$|R|+|S|+6= 25 + 14+6  $.

  If $G$ is  $F_4(q)$, then $G_1=\Sp(6,q)$,  and both $G_2 $  and  $L_1$  have 
type $A_2$.  We just saw that  $G_1$ has a presentation with 10 generators and
$ 5 7$ relations.    By Theorem~\ref{SL3 for 100}  we obtain    
 $|X|+|Y|=10 + 4$ and  
$|R|+|S|+6=  56+ 14+ 6=76$.

Similarly, if $G$ is the universal cover of $ {^2}\!E_6(q)$, then
$G_1=\SU(6,q)$ and  $G_2 $  and 
$L_1$ both have  type
$A_2$.  Using Theorems~\ref{unitary 100}(a)  and
\ref{SL3 for 100},  we obtain
 $|X|+|Y|=8+4$ and 
$|R|+|S|+6= 43+14+6 =63$.
 
Since $\langle X \mid  R \rangle$ and  $\langle Y \mid   S  \rangle$  have
bit-length $O( \log q)$, the same is true of our presentation  for $G$.

Once again we can kill the center of $G$ 
with at most one  additional relation  of bit-length $O(  \log q)$.~\qed

\medskip
\noindent
{\bf Proof of Theorems~\ref{A} and \ref{B}.}  As pointed out in  
Section~\ref{Introduction}, Theorem~\ref{A}  follows from Theorem~\ref{B} 
and \cite[Lemma~2.1]{GKKL}.  Theorem~\ref{B} is contained in
Theorems~\ref{Alt(n) and Sym(n) 100 down to 7}, 
\ref{SL2 10 relations}, \ref{SU3 10 relations}, \ref{SL3 for 100}, \ref{SLn
100},
\ref{unitary 100} and \ref{classical theorem}, together with Sections 
\ref{Suzuki  groups} and \ref{^2F_4(q)}.~\qed

\section {Additional presentations of classical groups}
\label {More presentations of classical groups}

We now provide an alternative to the approach in  the preceding section for
 presentations  of  classical groups.
We will use Section~\ref {Weyl groups} and its notation
concerning the group $W=W_n =\Z_2^{n-1}
\semi A_n$,    consisting of $n\times n$ real
monomial matrices with respect to an orthonormal basis $v_1 ,\dots,v_n$
 of $\R^n$. 
 
\subsection {\boldmath Groups of type $D_n$}  
\label {Groups of type $D_n$}  
\addcontentsline{toc}{subsection}{\protect\tocsubsection{}{\thesubsection}{Groups
of type $D_n$}}
\addtocontents{toc}{\SkipTocEntry}
By Theorem~\ref{SLn 100}, 
since $\PO^+(6,q)\cong\PSL(4, q)$
 we only need to consider the case
$n\ge 4$.  We will use Proposition~\ref{Wn} to imitate the argument in
Theorem~\ref{SLn 100}.

\begin{Theorem}
 \label {$D_n$ again}
The group $\Omega^+(2n,q)$ has  a presentation
with 
\begin{itemize}
\item [\rm (a)]   
 $8$ generators$,$
$31$ relations and bit-length $O(\log n + \log q)$ if $n\ge 4,$ and 

\item [\rm (b)] $8$ generators$,$
$29 $ relations and bit-length $O(  \log q)$  if $n=4$ or $5$. 
  
\end{itemize} 
 At  most one additional relation 
of bit-length $O( \log n+ \log q)$   
is needed in order to obtain a presentation for $\PO^+(2n,q)$.
\end{Theorem}

\proof

There is a hyperbolic basis $e_1, f_1,\dots,e_n , f_n$  of
$V=\F_q^{2n}$ associated with
$G=\Omega^+ (2n,q)$.  Then $W  $ consists of isometries.  Moreover, $W$  lies
in
$G$  and permutes the pairs $\{e_i,f_{i } \}$,  $1\le i\le n$: if
$n\ge5$ then $W$  is perfect, and for $n=4$ we can see this by restricting from the
group
$ \Omega^+(10,q)$.

Each element of $W$ can be viewed using two different vector spaces: 
  $\R^n$ and
$ V $.  In the action on $V$, we write elements in terms of
the  standard orthonormal basis $v_1,\dots,v_n$.  The resulting diagonal 
elements of
$W$ are the elements leaving each pair  $\{e_i,f_{i } \}$ invariant.  Since
these two views are potentially confusing (especially when $q$ is even), we
will initially write elements in both manners.  

We  digress in order to observe that, \emph{when $q$ is odd$,$  $W$ does not
lift to an isomorphic copy inside the universal cover $\hat G$ of $G$.}  
For, recall Steinberg's criterion
\cite[Corollary~7.6]{St2}: an involution in $\Omega^+(2n,q)$ lifts to an element  
of order 4 in the spin group if and only if the dimension of its $-1$ 
eigenspace  on $V$  is $\equiv 2$~(mod~4). 
 Apply this to $\diag(-1,-1,1,\dots,1 ) = (e_1,f_{1 } ) (e_2,f_{2 } )   \in
W$   in order to obtain an element of order 4 in $\hat G$, which  proves our
claim.
 
We view $F=\SL(3,q) $ as the subgroup of $G$ preserving the subspaces 
$\langle e_{1} ,e_{2},e_{3} \rangle $  and $\langle f_{1 } ,f_{2 },f_{3 }
\rangle $ while fixing the remaining basis vectors.
 We  use 
\smallskip
\begin{itemize} 
\item the   presentation  
$\langle X
\mid R
\rangle$ for
$F $  in Theorem~\ref{SL3 for 100} and
 
\item  the  presentation
$\langle Y\mid S\rangle $ for  $W=W_n  $ in Proposition~\ref{Wn}
(where $X$ and $Y$ are disjoint), 

\smallskip
\end{itemize}  
together with  the following elements:
\smallskip 

\begin{itemize} 
 
\item  $c = 
\begin{pmatrix}
0&1&\,0 \\
 0&0&1 \\
1&0&0
\end{pmatrix} \! , \,f =  
\begin{pmatrix}
0&1&0 \\
1&0&0 \\
0&0&\!\!-1
\end{pmatrix}\!
\in F$;
 
\smallskip
\item  $a
\in F$ such that  $L:=\langle a, a^f\rangle   \cong \SL(2,q)$  
consists of all matrices 
$ {\bigl( \begin{smallmatrix}
* &0 \\
0&1
\end{smallmatrix} \bigr)}  $
 in $F$%
;  
\smallskip
\item $(3,2,1)= (e_3,e_{2} ,e_{1}) (f_{3} ,f_{2 } ,f_{1 }) ,~ (1,3)(2,
4)= (e_1, e_{3})(e_2, e_{4})(f_{1 } , f_{3 })(f_{2 } , f_{4 })   \in   W$ and
$s=\diag(-1,-1,1,\dots,1 )= (e_1,f_{1 } ) (e_2,f_{2 } )\in Y$;  and 
\smallskip
\item  $\sigma$ and $\tau =(1,2)(3,4) = (e_1, e_{2})(e_3, e_{4})(f_{1 } ,
f_{2 })(f_{3 }, f_{4 })\in W$  that generate the subgroup of $W$  fixing
$\langle v_1- v_2\rangle$ (within $\R^n$) and   send  $v_1- v_2$ to $v_2-
v_1$.

\smallskip
\hspace{-14pt}
{\em Bit-length}: $c,$  $ f$ and $a$ have bit-length $O(\log q)$ using  
Remark~\ref{lengths in SL2}. 
We may assume that     $\sigma$ is the  product of
$\diag(-1,-1,-1,-1,1,\dots)$  and a
cycle of odd length $n-2$ or $n-3$ on $\{3,\dots,n\}$;  
both $\sigma$  and
$\tau$   can then be viewed as 
words in  $Y$  of  bit-length 
$ O(\log n)$ (by Remark~\ref{small bit-length}).

\end{itemize}

\smallskip 

We will show that $ G$  is isomorphic to the group $J$ having the
following presentation. 
\smallskip

 {\noindent \bf Generators:} $X, Y$.
\smallskip

{\noindent\bf Relations:}
\begin{itemize}
\item []
\begin{enumerate}
\item $R$.
\item $S$.
\item  $c = (3,2,1)$.
\item   $a^\sigma  = a^f,$ $ a^\tau  = a^f$.
\item   $(a^f)^\sigma   = a $.
\item   $[a, a^{(1,3)(2,4)}] = 1$.
\item   $[a^f, a^{(1,3)(2,4)}] = 1$ if $ n= 4$ or $5$.
\item  $[a, a^{s^{(3,2,1) }}] = 1$.
\item  $[a^f, a^{s^{(3,2,1) }}] = 1$ if $n=4$.
\end{enumerate}
\end{itemize}
\smallskip
 
There is a
surjection
$J 
\to G$, in view of the previous remark concerning the universal cover of $G$
together with the fact that the chosen element $f$ acts on $L$ in the same manner
as any element of $  W $ that interchanges $1$ and $2$.

As usual, we may assume that  
$F=\langle X
\rangle$ and 
$W =\langle Y \rangle $ lie in $J$.  
 By (4) and (5),
 $\langle \sigma,\tau\rangle$
normalizes  $L$ (since  $\tau$ has order 2 we also have 
$(a^f)^\tau   = a $).  As usual, it follows that 
$L^{W}$ can be identified with the set of   $n(n-1)$ pairs $\{\pm\alpha\}$ of
vectors
$\alpha=\pm v_i\pm v_j\in \R^n$, $i\ne j$, in the root system $\Phi$ of type
$D_n$. The groups in $L^{W}$  produce corresponding root groups $X_\alpha$,
$\alpha\in\Phi$.

 Any unordered pair of distinct, non-opposite roots can be moved by $W$ to 
one of the following pairs within $\R^n$:
\vspace{-2pt}
$$
\begin{array}{llll}
\mbox{(i) $v_1- v_2, \,  v_1 + v_2$}
&
\mbox{\quad   (ii)   $v_1- v_2,\,   v_1 - v_3$}
\vspace{4pt}
\\
\mbox{(iii)    $v_1- v_2, \,  v_3 - v_1 $ }
&
\mbox{\quad (iv)   $v_1- v_2, \,  v_3 - v_4$.
}
\vspace{-1pt}
\end{array}$$
Then the corresponding root groups can be moved in the same manner.
\vspace{2pt}

Let  $N:=\langle L^{W } \rangle =\langle X_\alpha\mid \alpha\in \Phi\rangle
=\langle F^{W} \rangle\noreq J $.
 
We need to verify the Steinberg relations for the  root groups $X_\alpha$.  The
pairs (ii) and (iii) already lie in   $F = \SL(3,q)$, so the desired relations are
immediate.  It remains  to consider the cases (i) and (iv).

 As in  (\ref{only for n<6}) (inside the proof of Theorem~\ref{SLn 100}), (6) and
(7)  imply that  the root groups determined by  (iv) commute.   
 
Before considering  (i), we note that (4) and (5) imply that every element  of $W$
that interchanges 
$v_1-v_2$ and $v_2-v_1$ also interchanges $a$ and $a^f$.  Two such elements are 
$s$ and, when $n\ge 5$, also  $(1,2)(4,5)$.  Since 
$t:=s^{(3,2,1)}=\diag(1,-1,-1,1,\dots,1)$  commutes with
$s$, and  $ t^{(1,2)(4,5) }=\diag(-1,1,-1,1,\dots,1) =s t $,
by (8) we have   
$$
 \begin{array} {lll}
1 =[a^s,(a^s)^t] = [a^f, (a^f)^t]  
\vspace{6pt}
 \\
1=  [a^{(1,2)(4,5)}, a^{t(1,2)(4,5)}] = 
[a^f, (a^{(1,2)(4,5) s})^t] = [a^f,a^t].\vspace{-2pt}
 \end{array}  
$$
By (9), the second of these relations also holds when $n=4$.  
  Then also 
$1= [a^f{}^t,a ]$.  Now the root groups $X_{v_1-v_2}<L_{v_1-v_2}=L=\langle a,a^f 
\rangle
$ and 
 $X_{v_1+v_2}<L_{v_1+v_2}=L^t=\langle a^t,a^f{}^t  \rangle $ commute, as required
for  (i).  

Thus,  $N \cong  G$.
Relation (3) pulls   $(1,2,3)$ into $N$,  so that $J/N=1$.

Now use Proposition~\ref{Wn} in order to obtain a presentation for 
$\Omega^+(2n,q)$ having the stated numbers of generators and relations. 
Finally, at most  one further relation 
of bit-length $O( \log n+ \log q)$  is needed to kill the center. 
\qed
 
\smallskip 
\smallskip 

This should be compared to the  presentations for these groups in
Section~\ref{General case} using 11 generators and   46 relations.

\subsection {\boldmath Groups of type $B_n$ and $C_n$} 
\addcontentsline{toc}{subsection}{\protect\tocsubsection{}{\thesubsection}{Groups
of type $B_n$ and $C_n$}}
\addtocontents{toc}{\SkipTocEntry}

This time we will glue  $W_n$ and  a group of type $B_2$ or $C_2$.

\begin{Theorem}
\label{classical again}
The groups $\Sp(2n ,q),$    $\Omega(2n+1,q)$ and 
$\Omega^-(2n+2,q)$ have  presentations with
\begin{itemize}
\item [\rm (a)]   $10$ generators$,$
$58$ relations and bit-length $O( \log n+ \log q)$ if $n\ge 4,  $   

\item [\rm (b)] $9$ generators$,$
$5 7 $ relations and bit-length $O(  \log q)$  if $n=4$ or $5,  $ and
 
\item [\rm (c)] $8$ generators$,$
$ 52$ relations and bit-length $O(  \log q)$  if $n=3$. 
\end{itemize} 
 At  most one additional relation 
of bit-length $O( \log n+ \log q)$   
is needed in order to obtain a presentation for  $\PSp(2n,q)$
or 
$\PO^-(2n+2,q)$.
\end{Theorem}

\proof 
The root system $\Phi$ of type
$C_n$ or $B_n$ for
 $G=\Sp(2n ,q)$,
$\Omega(2n+1,q)$   or  $\Omega^-(2n+2,q)$  consists of the vectors
$\pm v_i \pm v_j$ for $1\le i<j\le n$, and all
$\pm  2v_i$ or  $\pm v_i$,  respectively. 
We may assume that a fundamental system is 
$$
\Pi=\{ \alpha_1, ~ \alpha_j=v_{j+1} -v_{j }\mid 2\le j\le n-1 \},
\vspace{-2pt}
$$
where $\alpha_1= 2v_1$ or $v_1$.

We will use the subgroup
 $L_{12} \cong \Sp(4,q)$, $\Omega(5,q)$   or
$
\Omega^-(6,q)$ corresponding to the root subsystem 
$\Phi_{12}$  generated by $\alpha_1$ and $\alpha_2$, and its rank 1 subgroups 
$L_{1}$ and $L_{2}$ determined by 
$\pm \alpha_1$ and $\pm \alpha_2$, respectively. We will also need the subgroup 
$L_{23}$ corresponding to the root subsystem generated by $\alpha_2$
and $\alpha_3$.

There is a hyperbolic basis $e_1, f_1,\dots,e_n , f_n$  associated with $G$ (with  
additional basis vectors $v$, or $v$ and $v'$,  perpendicular  to all of these in
the orthogonal cases).  We may assume that
$W=W_n $ permutes these $2n$ vectors, with its normal subgroup  $\Z_2^{n-1}$ fixing
each pair $\{e_i, f_i  \}$.
We may assume  that the support of $L_{12}$ is $ \langle e_1, f_1, e_2, f_2
\rangle$ (or $ \langle e_1, f_1, e_2, f_2,v \rangle$  or
 $ \langle e_1, f_1, e_2, f_2,v,v' \rangle$ in the orthogonal cases); and
if  $n\ge4$ then  the support of $L_2^{(1,3)(2,4)}$ is $ \langle e_3,  
e_4 ,   f_3,   f_4
\rangle$,  so that 
\begin{equation}
\label{L's commute}
[L_{12} , L_2^{(1,3)(2,4)}] =1.
\end{equation}
 Similarly,
\begin{equation}
\label{L1 L2 commute}
[L_1,L_2^{c^2}]=1.
\end{equation}

We will use the following elements,  writing matrices for $L_2$ using the
vectors  $e_1,   e_2$ (and, implicitly, also  their ``duals''  $f_1,   f_2$),
and writing elements  of $W$ using $\R^n$:
\smallskip

\begin{itemize} 

\item $c=(1,2,3) =(e_1,e_{2} ,e_{3}) (f_{1} ,f_{2 } ,f_{3 })\in W$ and
$s=\diag(-1,-1,1,\dots,1 )= (e_1,f_{1 } ) (e_2,f_{2 } )\in Y$; 

\smallskip 
\item $ (2,3,4)=(e_2,e_{3} ,e_{4}) (f_{2} ,f_{3 } ,f_{4 }) 
, ~ (1,3)(2, 4)=(e_1,e_{3} ) (e_2,e_{4} ) (f_{1} ,f_{3 })(f_{2} ,f_{4 }) \in
W$ if $n\ge4$;

\smallskip 
\item  $\sigma,\tau \hspace{-1pt} \in\! W\!$   generating the stabilizer   
$W_{\{\pm v_1,\pm  v_2 \}}\! \hspace{-.5pt}=  \!W_{\{\{e_1,f_1\},\{ e_2,f_2\} \}}\!\hspace{-.5pt}  \cong
\Z_2^{n-1} \hspace{-.8pt}  \semi S_{n-2} $;%
\smallskip
\item  $h =\diag(\zeta^{-1},\zeta)$  
generating  the torus that  normalizes the root subgroups $X_{\pm  \alpha_2  }$
of $ L_2 :=L_{\alpha_2}$;

\smallskip 

\item  $u =\bigl( \begin{smallmatrix}
1&1 \\
0&1
\end{smallmatrix} \bigr)\in  X_{ \alpha_2  }$;

\smallskip
\item  $a \in L_{2}$ such that  $L_{2} =\langle a, a^s\rangle
$;   
\smallskip 
\item 
$a_2\in  L_{ 2}$  such that $L_{2} =\langle a_2,a_2^{r_2h} \rangle $;
\smallskip 
\item  $b \in L_{12}$ such that  $L_{12}=\langle b, b^s\rangle $;    and
 \smallskip 
\item two generators for $L_1 $.

\smallskip
\hspace{-14pt}
{\em Bit-length}: Once again, by  
Remarks~\ref{lengths in SL2}  and  \ref{small bit-length}
  the stated elements of 
$L_1 $, $L_2$   and $L_{12}$
have bit-length
$O(\log q)$, while  $\sigma$  and
$\tau$ have bit-length 
$ O(\log n)$. 
\end{itemize}
\smallskip

The required elements  $a,a_2$  and  $b$ exist.  For $b$, use
the fact that 
$\PSp(4,q) \cong \PO(5,q)$  and then, in the orthogonal 5- or 6-dimensional group
$L_{12}$,  choose an element $b$  of order $(q^2+1)/(2,q-1)$ such that 
$\langle b, b^s\rangle $ has no proper invariant subspace of dimension $>1$.
(This means that $\langle b, b^s\rangle $ is irreducible, except in the case
$\PO(5,q)$, $q$ even.)

We will use the following presentations:
\smallskip
\begin{itemize}

\item  a presentation   $\langle X \mid  R \rangle $  for  $L_{12} \cong
\Sp(4,q)$,
$\Omega(5,q)$  or
$  \Omega^-(6,q)$
when  $G= \Sp(2n ,q),$    $\Omega(2n+1,q)$ or 
$\Omega^-(2n+2,q)$, respectively;  and
 
\item a presentation   $\langle Y \mid  S \rangle $   for
$W=W_n $    (where $X$ and $Y$  are disjoint).

\end{itemize} 
\smallskip
 We have corresponding root groups $X_\alpha$
whenever $\alpha\in \Phi_{12}$.  

As in Section~\ref {Groups of type $D_n$},   in the orthogonal
cases with $q$  odd  $W$ does not lift to an isomorphic subgroup of the universal
cover.

We will show that $ G$  is isomorphic to the group $J$ having the
following presentation. 
\smallskip

 {\noindent \bf Generators:} $X, Y$.

\smallskip
{\noindent\bf Relations:}
\begin{itemize}
\item []
\begin{enumerate}
\item $R$.
\item $S$.
\item $x'{}^\sigma $ and $  x'{}^\tau $ written as words in $X$, for each $x'\in
\{b, b^s\}$.
\item $ (2,3,4) $  commutes with   $L_1$ if $n\ge4$.  
\item $c^{r_2 }=  r_2  ^2c^{2}$.   
\item $h h^c h^{c^2} =1$.
\item $a_2^{h^c} = a_2^{\diag(1,\zeta^{-1})}$  written as a word in $X$.
\item $[u^c,u ] = (u^{r_2 })^{c^{2}}$. 
\item $[u^c,u^{r_2 }]=1$.
%
%
\item $[b,a^{(1,3)(2,4)}]=1$ if $n\ge4$.    

\item $[L_1,L_2^{c^2}]=1$  if $n=3$.

\item $s $  written as a word in $L_{ 1 } \cup L_{ 2  } $  if $n=3$.
\end{enumerate}
\end{itemize}  
\smallskip

Note that the relations (2)--(6)  in Theorem~\ref{SL3 for 100} are the
present relations (5)--(9) for
 $L_{23}:=\langle L_2,c\rangle \cong \SL(3,q)$.   Also,  the
present relations (10)
and (11) follow  from (\ref{L's commute}) and  (\ref {L1 L2
commute}).   Therefore, as usual, there is a surjection $J  \to G$,
 and  we may assume that  
$L_{12}=\langle X\rangle$,
$W =\langle Y \rangle $ and  $L_{23} =\langle L_2,c\rangle$  lie in
$J$.

\smallskip
{\noindent \bf Case $n\ge4$:} 
 Relations (3) and (4) imply that 
$\langle W_{\{\pm v_1,\pm  v_2 \}} , (2,3,4)\rangle =W_{\{\pm v_1\}  }$
normalizes 
$
L_{ 1}$.  As usual, it follows that
$|L_{ 1}^W|=n$.   Similarly,
by (3), $|L_{12}^W|=n(n-1)/2$  and each element of $W_{\{\pm v_1,\pm  v_2 \}}$ acts
correctly on $L_{12}$.
Then the normalizer of $L_2$ in $W_{\{\pm v_1,\pm  v_2 \}}$ has index 2, so that 
$|L_{  2}^W|  \le  n(n-1) $.  As usual, it follows  that $|L_{  2}^W| = n(n-1) $.
Starting with  the root subgroups $X_{\alpha_1}$  and  $ X_{\alpha_2}$  of $L_{12}$,
in  $J$ we obtain the ``correct" set $X_{\alpha_1}^W \cup X_{\alpha_2}^W
=\{X_\alpha\mid \alpha\in \Phi\}$ of
root subgroups.  

We will verify that the Steinberg relations hold for 
$N:=\langle X_{\alpha_1}^W \cup X_{\alpha_2}^W\rangle =\langle L_{12}^{W}
\rangle\noreq J $, after which we will have $J/N =1$ by (10). 
Many of the required relations are already available for $L_{12}$
and $L_{23}.$

Any unordered pair $\alpha,\beta$ of distinct, non-opposite roots  
can be moved by $W $ to  one of 
 \vspace{-4pt}
$$
\mbox{(i) $\alpha_1,   \pm \alpha_2$,
\quad   (ii)  $\alpha_1,  \pm  \alpha_1\pm  \alpha_2$,
\quad  (iii)   $ \alpha_1,  \alpha_4, $  
\quad (iv)  $\alpha_2,   \alpha_3,~$ or
\quad (v)  $ \alpha_2,   \alpha_4$.  
}\vspace{-4pt}$$
Only pairs (iii) and (v)  are not inside $L_{12}$  or $L_{23}$. 

Since each element of $W_{\{\pm v_1,\pm  v_2 \}}$ acts correctly on $L_{12}$,
$s':=s^{(1,3)(2,4)}$ \emph{commutes with}  $L_{12} $  (cf. (\ref{L's commute})).  
  By (10), 
$$
 \begin{array} {lllll}
1=[b,a^{(1,3)(2,4)  }] ^{s }
\hspace{6.9pt}  =
 [b^{s },a^{s'(1,3)(2,4)}] 
 =
[b^s, a ^{(1,3)(2,4) }]  
\vspace{3pt}
 \\
1=[b^{},a^{(1,3)(2,4)}]^{s' } 
\hspace{3.6pt}  =
[b,(a^{s}){}^{(1,3)(2,4) }]  
\vspace{3pt}
\\
1=[b,a^{(1,3)(2,4) }]^{ss'} =
  [ b^{s'}{}^s,(a^{s's}){}^{(1,3)(2,4) }]
 =
  [ b ^s,(a^ s)^{(1,3)(2,4) }]. 
 \end{array}  
$$
Now  $[L_{12},L_{ 2}^{(1,3)(2,4)}]=[ \langle b, b^s\rangle ,  \langle a,
 a^s \rangle^{(1,3)(2,4)} ] =1$, which  takes care of the pairs  (iii) and (v). 

Thus, $J$ is a central extension of $G$.  We 
already noted that $J$ cannot be the universal cover in the orthogonal cases. 
Hence, $J=N \cong G$.

For the counts  in (a) and (b)   use Theorems~\ref{rank
2s} and \ref{unitary 100} together with Proposition~\ref{Wn}. 

\smallskip
{\noindent \bf Case $n=3$:}   Once again, $|L_{ 1}^W|=3$ and  $|L_{  2}^W|  =
6$, and we obtain  root groups $X_\alpha$, $\alpha  \in \Phi$. 

This time we consider the subgroup $M:=\langle L_{12},L_{23} \rangle$ of 
$N:=\langle X_{\alpha_1}^W \cup X_{\alpha_2}^W\rangle =\langle L_{12}^{W}
\rangle\noreq J $.
By the Curtis-Tits-Steinberg presentation   mentioned in 
Section~\ref{Preliminaries}, in order to prove that $G\cong M$ we only need
to consider pairs 
$X_\alpha, X_\beta$  of root groups  with $ \alpha,  \beta$  lying in the subsystem
of $\Phi$ determined by one of the pairs  $\{\alpha_1, \alpha_2\}$,  
$\{\alpha_2, \alpha_3\}$   or $\{\alpha_1, \alpha_3\}$. In the first two cases the
desired relations already hold in $L_{12}$ or $L_{23}$.  The last case is covered
by  (11). 

It follows that  $M\cong G$.  As in (\ref{$t^c$}), from (12) we obtain
$W=\langle s,c \rangle \le M$. Then $J=N=M \cong G$.
This time we use the presentation for $W_3\cong A_4$ in
(\ref{A4 A5}).  Hence,  in   (c) this presentation for $J$ uses only $6+2$ 
generators and at most $ 36+ 3+ 7+4$ relations.~\qed

\smallskip 
\smallskip 

Recall that we already had presentations for these groups
 in Section~\ref{General case} having at most  13 generators and at most 
$67$ relations.

\section{Concluding remarks}
\label{Concluding remarks}

\begin{comments}\rm
\label{Concluding remarks CGT}
Short and bounded presentations are goals of one aspect of 
Computational Group Theory (\cite[pp.~290-291]{Sim} and 
\cite[p.~184]{HEO}).   Such  presentations have various applications, such
as  in \cite{LG,KS2} for gluing together presentations in
a normal series  in order to obtain a
presentation for a given matrix group.
The  presentations in the present paper  are not short in the  
sense of length used in  
\cite{Sim,HEO,GKKL}.  However, they  have the
potential advantage that they are simpler than those in \cite{GKKL}, at least
in the sense of requiring fewer relations.  We hope that both types of
presentations will turn out  to be useful in Computational Group Theory.
\end{comments}
 
\begin{comments}\rm
\label{Concluding remarks awkward}
 The presentations for $S_n $ and  $ A_n$  in
Section~\ref{Alternating groups} that are related to prime numbers appear to
be practical.  The ones in 
Section~\ref{All alternating and symmetric groups} for general
$n$  have one unusual and awkward relation  $y=w$, expressing $y$ 
as a word  in $X\cup X^y$; see \eqn{8 to 7} and the description of this word 
in the proof of Theorem~\ref{Alt(n) and Sym(n) 100}.   Experimentation appears
to be needed in order  to find a ``nice'' additional relation of this sort.   That is,
the presentation in Section~\ref{An explicit presentation for $S_n$} is among the easiest to describe 
of the presentations obtained using our methods, but it may not be the best in practice.
 
\end{comments}

\begin{comments}\rm
\label{Small alternating groups}
We used the presentations \eqn{A4 A5} of $A_4$ and $A_5$  in the proofs of
Proposition~\ref{Wn} and Theorems~\ref{SLn 100}, \ref{unitary 100} and
\ref{classical again}. 
If we had instead used the corresponding presentations of $\SL(2,3)$ or 
$\SL(2,5)$, with 2 generators and 2 relations,  we could have saved an
additional relation.

There are further small-rank cases where we could have proceeded  in the same
manner.  Presentations are known for the universal central extension of 
$A_n$,
$n\le9$,  with 2 generators and 2 relations
\cite{CHRR,CHRR2}; and  for $A_{10}$ using 2 generators and 3 relations
\cite{Hav} (cf. Example~\ref{more AGL1 examples}(\ref{A10})).   For these small $n$    such presentations can be used in our 
presentations in order to  save  several relations. 

\end{comments}

\begin{comments} 
\label{one less generator}\rm
As in the  proof of  Corollary~\ref{any generators of Sn}, 
Lemma~\ref{d generators}  can be used  in
order to \emph{decrease the number of relations by one} in  
Theorem~\ref{A}.   The easiest way to see this is to use the
``$3/2$-generation" of all finite simple groups \cite{GK}, according to
which any one of  our generators $a$ is a member of a generating pair
$\{a,b\}$; and then proceed as in Corollary~\ref{any generators of Sn}. 
However, this uses an unnecessary amount of machinery, since each of our
generating sets contains a member $a$ for which a suitable $b$ can be found
without much difficulty. 

\end{comments}

\begin{comments}\rm
\label{Concluding remarks rank 2}
For the groups $S_n, A_n$ and $\SL(n,q)$ 
we were careful to make the number of relations small.  For the other groups
we were somewhat less careful.  It is likely that one can obtain
presentations of these groups with fewer relations using  ideas
provided here.

In particular, the presentations of rank 2 groups in Section~\ref{Rank 2
groups} undoubtedly could be improved using the more careful approach in
Section~\ref{SL3 100}.  For example, there should be no need for
both $u_{\a_i,1}$  and $u_{\a_i,2}$.

The number of relations in Proposition~\ref{Wn}  probably can be
improved somewhat by using  the ideas in Sections~\ref{Using 2-transitive
groups}  and
\ref{All alternating and symmetric groups} in place of  Theorem~\ref{Alt(n)
and Sym(n) 100 down to 7}.

 Phan-type presentations   for  
orthogonal groups \cite{BGHS,GHNS}  should help decrease the
numbers of relations in Sections~\ref{Rank 2 groups},
\ref{General case} and 
\ref{More presentations of classical groups}.
\end{comments}

\begin{comments}\rm
\label{Concluding remarks spin}
Theorem~\ref {$D_n$ again} did not deal with all central extensions of
orthogonal groups, in particular, it did not handle spin groups.  These can
be dealt with in the same manner, by using a double cover of $W_n$. 
\end{comments}

\begin{comments}\rm
\label{Concluding remarks Cayley graph}
There is a small generating subset  of each finite simple group $G$  producing 
a Cayley graph of diameter $O(\log|G|)$ \cite{BKL,Ka}.
Such generators need to be incorporated into   Theorem~\ref{A}  
 in order to obtain somewhat shorter presentations. 
\end{comments}

\begin{comments}\rm
\label{Concluding remarks expo-length}
  As observed in Section~\ref{Alternating groups} we have
constructed presentations for alternating and symmetric groups with bounded
expo-length.  The remaining
presentations in this paper do not have this property, unless we only
consider  groups over {\em bounded} degree extensions of the prime field.
An obstacle to our obtaining presentations with bounded expo-length is that
we do not know sufficiently nice presentations of  $\F_q$
when this field has   large degree  over the prime field $\F_p$.

In Section~\ref{Rank 1 groups }, we started with a presentation of  
$\F_q$ (as an algebra over $\F_p$)
of the form
$
\F_q = \F_p [ x, y]/\!\left( m(x),y-g_{\zeta^2}(x) \right)\!,
$
where $x,y$ map onto $\zeta^{2d}$ and $\zeta^2$, respectively,
and used it to obtain a presentation of  $\SL( 2,q)$.  

Roughly speaking, short (with length $O(\log q ) $) presentations of 
$\F_q$ as an algebra over $\F_p$ yield  presentations of
$\SL( 2,q)$ with short bit-length.
However in order to obtain a presentation 
of $\SL( 2,q)$ with bounded expo-length
using the same method, we would need a  presentation of  $\F_q$ in which
 every relation involves only a bounded number of monomials. 
This is a   computational question concerning ``sparse''
constructions of finite fields about which    little appears to be known
\cite{GaN}. The same remarks  also apply to the  unitary and Suzuki  groups.
\end{comments}

 \smallskip \noindent
\emph{Acknowledgements}:  We are grateful to George Havas for providing the
presentations in Example~\ref{more AGL1 examples}(\ref{PGL(2,9)},
\ref{A10}),
to  John Bray for  providing the
presentations  in  Section~\ref{Some specific
presentations of unitary groups}, and  to   Bob  Griess for pointing out
Steinberg's criterion
\cite[Corollary~7.6]{St2} used in Section~\ref{Groups of type $D_n$}.

\end{document}